\documentclass{amsart}
\usepackage{nicematrix}
\usepackage{tikz}
\usetikzlibrary{fit,patterns,arrows.meta,decorations.pathmorphing}
  
\usepackage[table]{xcolor}

\usepackage{a4wide}
\usepackage{amscd}
\usepackage{amssymb}
\usepackage{cite}
\usepackage{amsmath}
\usepackage{comment}
\usepackage{bookmark}
\usepackage{mathtools}
\usepackage{enumerate}
\usepackage{diagbox}
\usepackage{caption}

  \usepackage{mathrsfs}
\newtheorem{theorem}{Theorem}[section]
\newtheorem{thm}[theorem]{Theorem}

\newtheorem{definition}[theorem]{Definition}
\newtheorem{lemma}[theorem]{Lemma}

\newtheorem{corollary}[theorem]{Corollary}
\newtheorem{quest}[theorem]{Question}

\newtheorem{proposition}[theorem]{Proposition}
\newtheorem{prop}[theorem]{Proposition}
\newtheorem{rem}[theorem]{Remark}

  \usepackage{makecell}
\newcommand{\norm}[1]{\left\lVert#1\right\rVert}

\numberwithin{equation}{section}

\newcommand{\cH}{{\mathcal H}}

\newcommand{\cM}{{\mathcal M}}
\newcommand{\cN}{{\mathcal N}}

\newcommand{\cP}{{\mathcal P}}

\newcommand{\cZ}{{\mathcal Z}}

\usepackage{mathtools}
\usepackage{hyperref}                   % Reference
\hypersetup{colorlinks,%                % Reference color setup
    linkcolor=blue,%
    citecolor=blue}
\begin{document}

\title[On Isomorphic Embeddings of Noncommutative Lorentz Spaces]{%Isomorphic  embedding of $C
%_{p_1,q_1}$ into $C_{p_2,q_2}$, $0<p_1,p_2,q_1,q_2 <\infty $
On Isomorphic Embeddings of Noncommutative Lorentz Spaces
}

\author[C. Gao]{Chenhao Gao}
\address{Institute for  Advanced Study in  Mathematics of HIT, Harbin Institute of Technology, Harbin, 150001, China}
\email{{\color{blue}chenhao.gao@stu.hit.edu.cn}}

\author[J. Huang]{Jinghao Huang}
\address{Institute for  Advanced Study in  Mathematics of HIT, Harbin Institute of Technology, Harbin, 150001, China}
\email{{\color{blue}jinghao.huang@hit.edu.cn}}

\author[F. Sukochev]{Fedor Sukochev}
\address{School of Mathematics and Statistics, University of New South Wales,
Kensington, 2052, Australia}
\email{\color{blue} f.sukochev@unsw.edu.au}

\author[Z. Yu]{Zhizheng Yu}
\address{Institute for  Advanced Study in  Mathematics of HIT, Harbin Institute of Technology, Harbin, 150001, China}
\email{\color{blue}  zhizheng.yu@hit.edu.cn}

\thanks{J. Huang was supported the NNSF of China (No.12031004, 12301160 , 12471134 and 12671159)}


\begin{abstract} 
It is shown that for 
  $0<p_1,p_2,q_1,q_2 <\infty $,  
the Lorentz operator ideal  $C_{p_1,q_1}$ embeds  isomorphically into $C_{p_2,q_2}$ if and only if 
\begin{enumerate}
    \item $p_1=q_1=2$
    \item or $p_1=p_2$ and $q_1=q_2$,
\end{enumerate}
  which is in contrast to the classical result that $\ell_{q}$ embeds into $\ell_{p,q}$ as a (complemented) subspace. 
We also study the isomorphic embeddings between operator ideals generated by the weighted Lorentz sequence spaces $d(w,p)$ ($0<p<\infty$).
We characterize  those (commutative/noncommutative) Lorentz function spaces $L_{W}^p $ that are isomorphic to subspaces of noncommutative $L_p$-spaces for $1\leq p<2$ over some noncommutative probability spaces, which is  a noncommutative version of Sch\"{u}tt's theorem.
%Furthermore, we give a necessary and sufficient condition for  $d(w,p)$ ($0<p<\infty$) with exactly two nonequivalent symmetric basic sequences.
%Thus we establish the embedding property for $C_{d(w,p)}$ analogous to the $C_{p,q}$ - case.
\end{abstract}
\maketitle

\section{Introduction}
Lorentz spaces $L_{p,q}$ introduced by G. G. Lorentz in the 1950s~\cite{Lorentz50,Lorentz51} are important  generalizations of classical $L_p$-spaces.
Their noncommutative counterparts were later developed in the context of (semifinite/infinite) von Neumann algebras\cite{Kosaki, DDDE}.
These spaces play an important role in various branches of analysis, such as harmonic analysis and interpolation theory (see e.g. \cite{BS,Dilworth,CD85} and references therein).
Our theme in this paper is concerned with the isomorphic classification of $L_{p,q}$-spaces, notably, of the `atomic' version, that is ideals $C_{p,q}$ (of compact operators) for different values of $(p,q)$.
For the isomorphic classification of $L_{p,q}$-spaces over different measure spaces/semifinite von Neumann algebras, we refer to \cite{AHS,SS,KS,HS21,HSSZ,HSS,CD88,HSZ25}.
%The geometric structure of Lorentz spaces and their discrete counterparts (i.e. the Lorentz sequence spaces $\ell_{p,q}$) has been a central subject in Banach space theory.

\subsection{A noncommutative Sch\"{u}tt's theorem}
%An important topic in the field of Banach space geometry  is the classification of subspaces.
Let $B(\mathcal{H})$ be the $*$-algebra of all bounded operators on an infinite dimensional separable Hilbert space $\mathcal{H}$.
Fix an orthonormal basis in $\mathcal{H}$ (the particular choice of a basis is inessential) and identify the algebra $\ell_{\infty}$ of all bounded sequences with the subalgebra of all diagonal operators on $\mathcal{H}$ with respect to this basis. 
Given an ideal $C_E\subset B(\mathcal{H})$, the commutative core $E$ corresponding to $C_E$, which is also frequently
called a Calkin space (see, e.g., \cite{Kalton_S,LSZ})  is the collection of all diagonal operators from $B(\mathcal{H})$ belonging to $C_E$ and we denote it by 
\begin{align}\label{commucore}
E:= C_E \cap \ell_{\infty}.
\end{align}
If $\left(C_E, \norm{\cdot}_{C_E}\right)$ is a Banach ideal, then its commutative core $E$ is equipped with the induced norm $\norm{x}_E:=\norm{x}_{C_E}, \enspace x\in E$.
On the other hand, for any symmetric sequence space $E$, one can define an operator ideal $C_E$ (see Section \ref{Preliminaries} below or \cite{Kalton_S}). 
In \cite{JSZ2}, it was shown that if the commutative core $E$ ($\ne \ell_p$) of the ideal 
			$C_E$ admits an isomorphic embedding into $L_p(0,1)$, $1\le p<2$,  then $C_E$ embeds isomorphically into $\mathcal{L}_p(\cM,\sigma)$ for some finite von Neumann algebra $(\cM,\sigma)$. 
            Here, we are concerned with the converse to this implication stated in \cite{JSZ2} as an open question.
\begin{quest}\label{JSZquestion}\cite[Question 8]{JSZ2}
			Suppose that  (an ideal of compact operators) $C_E$ admits an isomorphic embedding into $\mathcal{L}_r(\cM,\sigma)$, $1\le r<2$, for some finite von Neumann algebra $(\cM,\sigma)$. Is it true that the commutative core $E$ of the ideal 
			$C_E$ admits an isomorphic embedding into $L_r(0,1)$? 
		\end{quest} 
%Question \ref{JSZquestion} above is the second motivation of the present paper. 
In this paper, we answer this question in the setting of Lorentz operator ideals $C_{p,q}$ (see Remarks~\ref{quesanws1} and  \ref{quesanws2} below) generated by $\ell_{p,q}$  and $C_{d(w,p)}$ (see Theorem~\ref{ifnv} below) generated by weighted Lorentz sequence spaces $d(w,p)$  (see the definition in section~\ref{Preliminaries} below). 

%Recall that 
% Sch\"{u}tt characterized 
An important contribution related to our theme is Sch\"utt's characterization \cite{Schutt} of those Lorentz function spaces $L_W^p(0,1)$ (see section~\ref{Preliminaries} below) which are subspaces  of $L_1(0,1)$.
 \begin{theorem}\cite[Theorem 1.1]{Schutt}
		Let $1\le p<2$ and $W$ be a nonincreasing, positive, continuous function on $(0,1)$ with 
		$\int_0^1 W(t)dt =1$. The  following are equivalent:
\begin{enumerate}
    \item[(i)] $L_W^p(0,1)$ is isomorphic to a subspace of $L_p(0,1)$.
    \item[(ii)] $L_W^p(0,1)$ has cotype $2$ (is $2$-concave).
    \item[(iii)] there is a constant $C>0$ such that for all $x\in (0,1) $,
			we have 
			$$\int_0^x W(t)t^{-p/2} dt \le C x^{-p/2} \int_0^x W(t)dt . $$  
\end{enumerate}
	\end{theorem}
As an application, it is shown in \cite[Corollary 1.5]{Schutt} (see also \cite{CD91}) that $L_{p,q}(0,1)$   is isomorphic to a subspace of $L_1(0,1)$ if and only if $1\le q\le p<2$ or $p=q=2$. 

Let $\cM$ be an atomless semifinite $\sigma$-finite von Neumann algebra equipped with a semifinite infinite faithful normal trace $\tau$.
	Recently,   it was established  that a symmetric (Banach) function space $E(0,\infty)$ with $E(0,\infty)\cap L_{\infty}(0,\infty)\not =L_p(0,\infty)\cap L_{\infty}(0,\infty)$, $1\leq p<2$, embeds isomorphically into $L_p(\cN,\nu)$ for some noncommutative probability space $(\cN, \nu)$ if and only if its noncommutative analogue $E(\cM,\tau)$ embeds isomorphically into $L_p(\hat{\cN},\hat{\nu})$ for some noncommutative probability space $(\hat{\cN},\hat{\nu})$\cite[Theorem 1]{HJSZ}. 
	This   motivates us to establish a noncommutative analogue of Schütt's theorem.
   % This result connects the embedding properties of noncommutative Lorentz spaces to the geometry of their classical counterparts.
	\begin{theorem}
	    \label{nvos}
		Let $1 \leq r\leq p < 2$ and $ W $ be a nonincreasing, positive, continuous function on $(0, 1)$ with $\int_0^1 W(t)\,dt=1$. Let $\cM$ be an atomless finite von Neumann algebra equipped with a faithful normal trace $\tau$.
		The following are equivalent.
\begin{enumerate}
            \item[(i)] there exists a noncommutative probability space $(\cN, \nu)$ such that $L_W^p(0,1)$ is isomorphic to a subspace of $L_r({\cN,\nu})$.
		\item[(ii)] there exists a noncommutative probability space $(\hat{\cN},\hat{\nu})$ such that $L_W^p(\cM,\tau)$ is isomorphic to a subspace of $L_r(\hat{\cN},\hat{\nu})$.
		\item[(iii)] $L_W^p(0,1)$ is isomorphic to a subspace of $L_p(0,1)$.
		\item[(iv)] $L_{W}^p(0,1)$ is of  cotype $2$ (is $2$-concave). 
		\item[(v)] there is a constant $ C > 0$ such that for all $x \in (0, 1)$, we have
            $$
			\int_0^x W(t) t^{-p/2} \, dt \leq Cx^{-p/2} \int_0^x W(t) \, dt.
		$$
\end{enumerate}
		The implication $(i)\Rightarrow(ii)$ holds true when $\cM$ is an arbitrary (not necessarily atomless) finite von Neumann algebra equipped with a faithful normal trace.
	\end{theorem}
We also establish Schütt's theorem for semifinite $\sigma$-finite infinite von Neumann algebras 
(see Theorem \ref{ifnv} below). 
    %{\color{red}On the other hand, we 
 % characterize those parameters $1\leq q\leq p<2$ and $1\leq r$  for which $L_{p,q}(\cM,\tau)$ embeds into $L_r(\cN,\nu)$ for some noncommutative probability space $(\cN,\nu)$.
%In particular, % we show that that 
	%Descending to the discrete setting, we address the embedding of Lorentz sequence spaces.
%    we show that for $1\leq r<2$, there exists some finite noncommutative probability space $(\cN,\nu)$ such that
%	$\ell_{r,q}\hookrightarrow L_r(\cN,\nu)\mbox{ and } q=r$ 
 %   (whereas if $r>2$, then $\ell_{r,q}\not\hookrightarrow L_r(\cN,\nu)$ for all  $1\le q<\infty$) and $\ell_{p,r}\hookrightarrow L_r(\cN,\nu)$
	%	if and only if
	%	$r\le p<2$  }.  

\subsection{Isomorphic embeddings between Lorentz ideals}
It is shown in \cite{Schutt} that  $L_{p,q}(0,1)$  does not embed isomorphically into $L_1(0,1)$ when $q>p>1$, by a deep result due to Kwapien and Sch\"{u}tt\cite{KS89}  (see also \cite[p. 76]{CD91} for another argument relying on a theorem due to  Guerre and Levy\cite{GL}). 
It is not known   whether $L_{p,q}(\cM,\tau)$ embeds  into the predual of some von Neumann algebra $ \cN $ when $q>p>1$. 
 This naturally motivates  a more general question:
 \begin{quest}
 Let $(\cM,\tau)$ and $(\cN,\nu)$ be two semifinite von Neumann algebras.
     For which $0<p_1,p_2,q_1,q_2<\infty$, $L_{p_1,q_1}(\cM,\tau)\hookrightarrow L_{p_2,q_2}(\cN,\nu)$?
 \end{quest}
The second main result of the present paper addresses this question in the setting when $\cM=\cN =B(H)$. 

 Recall that Arazy's result \cite[Section 3]{Arazy81} yields   $C_{p,q}\not\hookrightarrow C_1$
when  $p,q > 1$ and $(p,q)\ne (2,2 )$. 
%Having characterized the embedding properties in the general continuous framework, a natural question arises: Do these embedding phenomena persist when we restrict our attention to the 'atomic' setting of $B(H)$, specifically the Schatten-Lorentz operator ideals? 
%To answer this, we utilize the theory of (symmetric) basic sequences in unitary matrix ideals $C_E$. 
%A cornerstone is the work of Arazy \cite{Arazy81}, which characterizes embeddings between such spaces in terms of tensor products. 
In the present paper, 
we extend Arazy's result \cite[Corollary 5.3]{Arazy81} from the Banach setting to the quasi-Banach setting, 
 which is the crucial tool in establishing the second main result of the present paper (see Lemma~\ref{arazylemmaforisomorphism} below).
 We denote the standard unit matrices by $\{e_{i,j}\}_{i,j\ge1}$ with $e_{i,j}(k,l)=\delta_{i,k}\cdot \delta_{j,l}$ (see also Definition~\ref{definieij}).
The   proof will be presented in Appendix \ref{Arazy-qbs}.

\begin{lemma}\label{arazylemmaforisomorphism}
Let $E,F$ be two quasi-Banach symmetric sequence spaces such that $c_0,\ell_2\not\hookrightarrow E$ and $\ell_2\not\hookrightarrow F$. If the lower triangular part   $T_F$ of $C_F$ is isomorphic to a subspace of $C_E$, then there exists $b\in C_E$ such that the mapping  $W:C_F\to C_E$ defined by $$We_{i,j} = b\otimes e_{i,j},~i,j\ge 1,$$
is an isomorphism. 
\end{lemma}

It is a direct corollary of the main result in \cite{Arazy81} (see also \cite{AL} and \cite{HSS}) that $C_{p} \not\hookrightarrow C_q$ if $1\le p,q  <\infty, $  $
p\ne q, 2   $. 
It is well-known that $\ell_q$ is a complemented subspace of $\ell_{p,q}$, see e.g. \cite{Dilworth,CD89}.
Leveraging Lemma~\ref{arazylemmaforisomorphism},  we show that the noncommutative counterpart of this result fails, i.e., if $q\ne 2$ and  $p\ne q$,   then $C_q$ does not embed into $C_{p,q}$.
As an application, we fully characterize those $C_{p_1,q_1}$-ideals which can be embedded into $C_{p_2,q_2}$. 

 \begin{theorem}Let $0<p_1,p_2,q_1,q_2 <\infty $. Then, $C_{p_1,q_1}\hookrightarrow C_{p_2,q_2}$ if and only if 
 \begin{enumerate}
 \item either $p_1=q_1=2$,
 \item or $p_1=p_2$ and $q_1=q_2$. 
 \end{enumerate}
 
 \end{theorem}

\subsection{Isomorphic embeddings between weighted Lorentz ideals}
We also study the isomorphic embeddings between  $C_{d(w,p)}$-ideals  generated by  weighted Lorentz sequence spaces $d(w,p)$  for $0<p<\infty $.
  Altshuler, Casazza, and Lin \cite[Theorem 6]{ACL} gave a necessary and sufficient condition for $d(w,p)$ for $1\le p<\infty$ to possess exactly two nonequivalent symmetric basic sequences.
We extend this result to the quasi-Banach setting (i.e., $0<p<\infty$, see Theorem~\ref{fourmainthm} below).
Consequently, we characterize those $C_{d(w_1,p_1)}$-ideals which can be embedded into $C_{d(w_2,p_2)}$ when $d(w_2,p_2)$ has only two nonequivalent symmetric basic sequences (see Theorem~\ref{weighted} below). %This is useful in the study of  isomorphic embeddings between operator ideals $C_{d(w,p)}$, $0<p<\infty$. 

Concerning terminology in the theory of quasi-Banach  spaces,  we
shall follow the standard terminology used in Banach spaces. In particular, the
notions related to Schauder basis, such as ``basic sequence''
and ``symmetric basis'',  are analogous to the corresponding notions in Banach
spaces and can be found in \cite{LT1,LT2,KR07,Popa}.

\section{Preliminaries}\label{Preliminaries}
	
	Recall that two Banach spaces $X$ and $Y$ are said to be {\it isomorphic} (denoted by $X\approx Y$) if there exists a linear bijection $T:X\rightarrow Y$ such that $T$ and $T^{-1}$ are bounded.
	We say that $X$ embeds into $Y$ if there is a linear subspace $Z\subset Y$, such that $Z$ is a Banach space and $Z\approx X$.
	If there exists an isomorphic embedding of $X$ into $Y$, then  we write $X\hookrightarrow Y$.
	If $X$ does not embed into $Y$ isomorphically, then  we write $X\not\hookrightarrow Y$. 

    A sequence $\{e_n\}_{n\ge 1}$ in $X$ is called a basis if for every $x\in X$ there is a unique sequence of scalars $\{\alpha_n\}_{n\ge 1}$ such that $x=\sum_{n=1}^{\infty}\alpha_n e_n$.
    If $\{x_n\}_{n\ge 1}$ and $\{y_n\}_{n\ge 1}$ are bases of $X$ and $Y$, respectively, we say that $\{x_n\}_{n\ge1}$ dominates $\{y_n\}_{n\ge 1}$, and write $\{x_n\}_{n\ge 1}>\{y_n\}_{n\ge1}$, if $\sum_{n=1}^{\infty}\alpha_n x_n$ converges in $X$ implies that $\sum_{n=1}^{\infty} \alpha_n y_n$ converges in $Y$.
    The basis $\{x_n\}_{n\ge 1}$ is equivalent to the basis $\{y_n\}_{n\ge1}$, and we write $\{x_n\}_{n\ge1} \sim \{y_n\}_{n\ge1}$, if $\{x_n\}_{n\ge 1}>\{y_n\}_{n\ge1}$ and $\{y_n\}_{n\ge 1}>\{x_n\}_{n\ge1}$.
    If $\{x_n\}_{n\ge 1}\sim\{y_n\}_{n\ge1}$, then $X\approx Y$.
\subsection{Symmetric spaces}
	Let $I$ be $(0,1)$ or $(0,\infty)$ and let $L(I)$ denote the space of all Lebesgue-measurable functions $f$ on $I.$
	For a Lebesgue-measurable function $f$ on $I,$ we define its {\it distribution function} by the formula
	$$d_f(s)=m(\{t:\ |f(t)|>s\}),\quad s\ge 0,$$
	where $m$ stands for Lebesgue measure. 
	Denote by $S(I)$ the subalgebra of $L(I)$
	consisting of all functions $f$ such that $d_{f}(s) < \infty$  for some 
	$s > 0$.
	
	Two measurable functions $f$ and $g$ are called {\it equimeasurable} (written, $f\sim g$) if their distribution functions $d_f$ and $d_g$ coincide. In particular, for every measurable function $f\in S    (I),$ the function $f$ is equimeasurable with its {\it decreasing rearrangement} $\mu(f)$ defined by the formula
	$$\mu(t;f):=\inf \{s\geq0:\ d_{f}(s)<t \},\quad t>0.$$
	If $ f,g\in S(I),$ then $\mu(f)=\mu(g)$ if and only if $f$ and $g$ are equimeasurable.  % We recall that a function $x$ is said to be symmetrically distributed, if $x$ and $-x$ are equimeasurable.

    \begin{definition}

	A linear subspace $E$ in $S(I)$ equipped with a quasi-Banach norm $\left\|\cdot\right\|_E$ is said to be {\it symmetric space} (see \cite{KPS,LT2})
		if $f\in E$, $g\in S(I)$ with $\mu(g)\leq \mu(f)$ imply that $g\in E$ and $\left\|g\right\|_E\leq \left\|f\right\|_E$.
       \end{definition}
	Without loss of generality, in what follows we always assume that $\left\|\chi_{(0,1)}\right\|_E=1$ for any symmetric function space  $E=E(0,1)$ on $(0,1)$.

\subsection{Noncommutative symmetric spaces}
For a detailed exposition of material in this subsection,  we refer to \cite[Chapter 6.7]{DPS2023}.
	In what follows,  $\cH$ is a  Hilbert space and $B(\cH)$ is the
	$*$-algebra of all bounded linear operators on $\cH$ equipped with the uniform norm $\left\|\cdot\right\|_\infty$, and
	$\mathbf{1}$ is the identity operator on $\cH$.
	Let $\mathcal{M}$ be
	a von Neumann algebra on $\cH$.
	We denote by $\cP(\cM)$ the collection  of all projections in $\cM$, by $\cM'$ the commutant of $\cM$ and by $\cZ(\cM)$ the center of $\cM$.

	A closed, densely defined operator $x:\mathfrak{D}\left( x\right) \rightarrow \cH $ with the domain $\mathfrak{D}\left( x\right) $ is said to be {\it affiliated} with $\mathcal{M}$
	if $yx\subseteq xy$ for all $y\in \mathcal{M}^{\prime }$, where $\mathcal{M}^{\prime }$ is the commutant of $\mathcal{M}$.
	A  closed,
	densely defined
	operator $x:\mathfrak{D}\left( x\right) \rightarrow \cH $ affiliated with $\cM $ is said to be
	{\it measurable}  if  there exists a
	sequence $\left\{ p_n\right\}_{n=1}^{\infty}\subset \cP\left(\mathcal{M}\right)$, such
	that $p_n\uparrow \mathbf{1}$, $p_n(\cH)\subseteq\mathfrak{D}\left(x\right) $
	and $\mathbf{1}-p_n$ is a finite projection (with respect to $\mathcal{M}$)
	for all $n$.
	The collection of all measurable
	operators with respect to $\mathcal{M}$ is denoted by $S\left(
	\mathcal{M} \right) $, which is a unital $\ast $-algebra
	with respect to strong sums and products (denoted simply by $x+y$ and $xy$ for all $x,y\in S\left( \mathcal{M%
	}\right) $).

	From now on, let $\mathcal{M}$ be a
	semifinite von Neumann algebra equipped with a faithful normal
	semifinite trace $\tau$.

	An operator $x\in S\left( \mathcal{M}\right) $ is called $\tau$-measurable if there exists a sequence
	$\left\{p_n\right\}_{n=1}^{\infty}$ in $P\left(\mathcal{M}\right)$ such that
	$p_n\uparrow \mathbf{1}$, $p_n(\cH)\subseteq \mathfrak{D}\left(x\right)$ and
	$\tau(\mathbf{1}-p_n)<\infty $ for all $n$.
	The collection $S\left( \mathcal{M}, \tau\right)
	$ of all $\tau $-measurable
	operators is a unital $\ast $-subalgebra of $S\left(
	\mathcal{M}\right) $.

	\begin{definition}
		Let $x\in
		S(\mathcal{M},\tau)$. The generalized singular value function $\mu(x):t\mapsto  \mu(t;x)$ of
		the operator $x$ is defined by setting
		$$
		\mu(s;x)
		=
		\inf\{\left\|xp\right\|_\infty:\ p\in \cP(\cM)\mbox{ with}\ \tau(\mathbf{1}-p)\leq s\}, ~ \forall s\in (0,\infty ).
		$$
	\end{definition}

	%It is well-known \cite{LSZ,DP2} that if $a\in S(\cM,\tau)$ and $b,c\in \cM$, then
	%\begin{align}\label{ineqfact}
	%\mu(t;bac)\le \left\|b\right\|_\infty \left\|c\right\|_\infty \mu(t;a), ~\mu(t;a^*)=\mu(t;a).
	%\end{align}
	%{\color{red}In particular, $\mu(A) =\mu(|A|)$.}
	%{\color{red}Suppose that $X\in S(\cM,\tau)$. If $0<\alpha \in \mathbb{R}$ and $E = E^{|X|}(\alpha,\infty)$, then
		%\begin{align}\label{mu_shift}
		%\mu(|X|E) = \mu(X)\chi_{[0,\tau(E))}
		%\end{align}
		%and
		%\begin{align}\label{mu_shift2}
		%\mu(t; |X|E^\perp) = \mu(t+\tau(E);X)
		%\end{align}
		%for all $t\ge 0$ whenever $\tau(E)<\infty$.}

	\begin{definition}
		A linear subspace $E$ of $S(\cM,\tau)$ equipped with a complete norm $\norm{\cdot}_E$, is called a symmetric space (of $\tau$-measurable operators) if $x\in S(\cM,\tau)$, $y \in E$ and $\mu(x)\le \mu(y)$ imply that $x\in E$ and $\norm{x}_E \le \norm{y}_E$.
	\end{definition}
	
	It is well-known that any symmetric space $E$ is a normed $\cM$-bimodule, that is, $axb\in E$ for any $x\in E$, $a,b\in \cM$ and $\left\|axb\right\|_E\leq \|a\|_\infty\left\|b\right\|_\infty \left\|x\right\|_E$ \cite{DP2,DPS2023}.
	%A linear subspace $E$ of $S(\cM,\tau)$ equipped with a complete norm $\|\cdot\|_E$, is called \emph{fully symmetric space} (of $\tau$-measurable operators) if $X\in S(\cM,\tau)$, $Y \in E$ and $X\prec\prec Y$ imply that $X\in E$ and $\|X\|_E \le \|Y\|_E$.
	%
	
	%
	%
	%
	%
	%The so-called K\"{o}the dual is identified with an important part of the dual space. If $E  \subset S(\cM,\tau)$ is a symmetric space, then the K\"{o}the dual $E^\times $ of $E$ is defined by setting $$ E^\times =\{   X\in S(\cM,\tau) : \sup_{\|Y\|_E\le 1, Y\in E}\tau (|XY|)   <\infty    \}.$$
	%
	
	A wide class of symmetric operator spaces associated with the von Neumman algebra $\cM$ can be constructed from concrete symmetric function spaces studied extensively in e.g. \cite{KPS}. 
    Let $ E(0,\infty )$ (respectively, $E(0,1)$) be a symmetric function space  and let $(\cM,\tau)$ be a semifinite von Neumann algebra equipped with a semifinite faithful normal trace $\tau$ (respectively, a noncommutative probability space, i.e., $\tau$ is a faithful normal tracial state). 
     Then, the pair given by 
	$$E(\cM,\tau)=\{x\in S(\cM,\tau):\mu(x)\in E(0,\infty) \},\quad \left\|x\right\|_{E(\cM,\tau)}:=\left\|\mu(X)\right\|_{E }$$ is a symmetric space on $\cM$ \cite{Kalton_S} (see also \cite{DPS2023,LSZ}). 
If $ E$ is a symmetric sequence space and $\cM=B(\cH)$ is equipped with the standard trace, then we denote 
\cite{LSZ,GK,DPS2023,Kalton_S}
$$C_E := \{x\in B(\cH): \{ \mu(i; x)\}_{i\ge 0}\in  E \},\quad \left\|x\right\|_{C_E}:=\left\|\{\mu(i;X)\}_{i\ge 0}\right\|_{E } .$$
    
\subsection{Lorentz spaces}

For  $0<p,q<\infty$ and $I\subset (0,\infty)$, the Lorentz space $L_{p,q}(I)$ consists of all measurable functions on $I$ satisfying 
	$$
	\left\|f\right\|_{p,q}:=\left(\int_I \mu(t;f)^q\,dt^{q/p}\right)^{1/q}<\infty,
	$$
	where $\mu(t;f)$ stands for the decreasing rearrangement of $|f|$. 
    It is well known that for $1\leq q\leq p<\infty$, $\left\|\cdot\right\|_{p,q}$ defines a norm under which $L_{p,q}$ is a separable   Banach symmetric function space; if $1< p< q<\infty$, $\left\|\cdot\right\|_{p,q}$ defines a quasi-norm on $L_{p,q}$, which is  known to be equivalent to a norm  \cite{BS}.

    The Lorentz sequence space is the special case $L_{p,q}(\mathbb{N})$ when the measure space $\Omega$ is $\mathbb{N}$ of all natural numbers equipped with the counting measure.
    We denote by $\ell_{p,q}$ consisting of all scalar sequences $x=(x_i)_{i=1}^{\infty}$ for which $\left\|(x_i)\right\|_{p,q}<\infty$, where 
    $$
    \left\|(x_i)\right\|_{p,q}=\left(\sum_{i=0}^{\infty}\mu(i;x)^q(i+1)^{\frac{q}{p}-1}\right)^{1/q}, \quad 0<p,q<\infty
     $$
    \noindent Clearly, $\ell_{p,q}$ is isomorphic  to a sublattice of $L_{p,q}(0,\infty)$.
    For any $p\ge 1$, we have $L_{p,p}=L_p$ and $\ell_{p,p}=\ell_{p}$; in this case, we will simply write $\left\|\cdot\right\|_p$.

 The weighted Lorentz function space $L_W^p(I)$ is the space of all measurable functions $f$ on an interval $I\subset (0,\infty )$
	with norm 
	$$\norm{f}_{L_W^p} =\left( \int_I  \mu(t;f)^p W(t)dt  \right)^{1/p}<\infty , \quad 0<p<\infty$$
	where a nonnegative and measurable function $W:I\rightarrow (0
	,\infty)$ satisfies the conditions
     $$
      \int_{0}^{t}W(s)\,\mathrm{d}s \;<\;\infty
      \quad\text{for all }t\in I,
     $$
     and when \(I=(0,\infty )\),
     $$
     \int_{0}^{\infty}W(s)\,\mathrm{d}s = \infty.
     $$
Note that $L_{p,q}=L_W^q$ with $W(t)=\frac{q}{p}t^{q/p-1}, \enspace t\in I$.

For every $0<p<\infty$ and every nonincreasing sequence of nonnegative numbers $w=\{w_n\}_{n\ge1}$, the weighted Lorentz sequence space $d(w,p)$ is the space of all sequences of scalars $x=(a_1,a_2,\cdots)$ for which 
$$
\norm{x}_{d(w,p)} = \sup_{\pi} \left(\sum_{i=1}^{\infty} |a_{\pi(i)}|^p w_i\right)^{1/p}<\infty,
$$
where $\pi$ denotes all permutations of positive integers.
It is well known that $\norm{\cdot}_{d(w,p)}$ is a $p$-norm for $0<p<1$ and a norm for $1\le p<\infty$ (see \cite{Popa,NO}, \cite[Section 
 7]{Dilworth} or \cite[p. 56, Proposition 3.6]{BS}).
It is easily checked that if $\inf_{n} w_n>0$ then $d(w,p)\approx\ell_p$.
If $\sum_{i=1}^{\infty} w_i<\infty$, then $d(w,p)\approx \ell_{\infty}$  \cite{LT1}.
In the present paper, we assume that $\sum_{i=1}^{\infty} w_i=\infty$.
If $E=d(w,p)$, 
then the operator ideal $C_E$ is a 
  weighted Lorentz ideal, denoted by  $C_{d(w,p)}$.

 We denote the unit vector basis of $\ell_{p,q}$ and $d(w,p)$ by $\{e_k^{\ell_{p,q}}\}_{k=1}^{\infty}$ and $\{e_k^{d(w,p)}\}_{k=1}^{\infty}$, respectively. 
Below, we recall two  well-known results.
\begin{lemma}[see \mbox{\cite[Lemma 2.1]{Dilworth}  or  \cite{CD88})}]\label{Lpqequivalp}
    Let $0<p,q<\infty$ and $\{e_n^{\ell_{p,q}}\}_{n\ge 1}$ be the unit vector basis of $\ell_{p,q}$.
    If $\left \{ y_n:=\sum_{i=p_n+1}^{p_{n+1}} \alpha_i e_i^{\ell_{p,q}}\right\}_{n\ge 1}$ is a semi-normalized   block basic sequence of $\{e_n^{\ell_{p,q} }\}_{n\ge 1}$ such that $\lim_{n\rightarrow \infty} \alpha_n=0$, then there exists a subsequence of $\{y_n\}_{n\ge 1}$ which is equivalent to the unit vector basis of $\ell_q$.
\end{lemma}
\begin{lemma}[\mbox{see \cite[Lemma 1]{ACL} and \cite[Lemma 3.1]{Popa}}]\label{Lorequivalp}
    Let $0<p<\infty$ and $\{e_n^{d(w,p)}\}_{n\ge 1}$ be the unit vector basis of $d(w,p)$.
    If $\left\{ y_n:=\sum_{i=p_n+1}^{p_{n+1}} \alpha_i e_i^{d(w,p)}\right \}_{n\ge 1}$ is a semi-normalized block basic sequence of $\{e_n^{d(w,p)}\}_{n\ge 1}$ such that $\lim_{n\rightarrow \infty} \alpha_n=0$, then there exists a subsequence of $\{y_n\}_{n\ge 1}$ which is equivalent to the unit vector basis of $\ell_p$.
\end{lemma}

\section{A noncommutative version of Sch\"{u}tt's theorem}
\subsection{Proof of Theorem \ref{nvos}}
For an (atomless) noncommutative probability space $(\cM,\tau)$ and  a symmetric function space  $E(0,1)$ such that both  $E(0,1)$  and its K\"othe dual have  the Kruglov property\cite{SA2011}, 
    $E(0,1)$ embeds isomorphically into $L_p(\cN,\nu)$, $1\leq p<2$, for some noncommutative probability space $(\cN, \nu)$ if and only if $E(\cM,\tau)$  embeds isomorphically into $L_p(\hat{\cN}, \hat{\nu})$ for some noncommutative probability space $(\hat{\cN}, \hat{\nu})$\cite[Theorem 4]{HJSZ}.
	Unlike $L_{p,q}$, the space $L_W^p$ may not have the Kruglov property (see e.g. \cite[Theorem 5.1]{AS2005} for a concrete example).
	Thus, \cite[Theorem 4]{HJSZ} is not enough to establish a noncommutative analogue of Sch\"{u}tt's theorem.
	In order to prove Theorem \ref{ifnv},   We need a version of Sch\"utt's theorem for Lorentz spaces on $(0,\infty)$ (see Proposition \ref{ifvos} below, whose proof can be found in   Appendix~\ref{sch-inf}), 
    and   employ \cite[Theorem 1]{HJSZ} instead of \cite[Theorem 4]{HJSZ}.

      \begin{prop}\label{ifvos}
		Let $1\leq p<2$ and $W$ be a nonincreasing, positive, continuous function on $(0, \infty)$ such that $
      \int_{0}^{t}W(s)\,\mathrm{d}s \;<\;\infty
      $ for all $t\in (0,\infty )
     $. 
		Then the following are equivalent:
		\begin{enumerate}
    \item[(i)] $L_{W}^p(0,\infty)$ is isomorphic to a subspace of $L_p(0,1)$. 
    \item[(ii)] $L_{W}^p(0,\infty)$ has cotype $2$ (is $2$-concave). 
    \item[(iii)] There is a $C > 0$ such that for all $x \in (0, \infty)$, we have
		\begin{equation}
			\int_0^x W(t) t^{-p/2} \, dt \leq Cx^{-p/2} \int_0^x W(t) \, dt .
			\label{ineq:ifvos}
		\end{equation}
\end{enumerate}
	\end{prop}

    The following theorem follows from a combination of 
	 Proposition~\ref{ifvos} and  \cite[Theorem 1]{HJSZ}.

         \begin{thm}\label{ifnv}
		Let $1 \leq r\leq p< 2$ and $W$ be a nonincreasing, positive, continuous function on $(0, \infty)$.  
		Let $\cM$ be an atomless semifinite $\sigma$-finite von Neumann algebra equipped with a semifinite infinite faithful normal trace $\tau$.
		If $L_W^p(0,\infty)\cap L_{\infty}(0,\infty)\not=L_r(0,\infty)\cap L_{\infty}(0,\infty)$ in the sense of
sets, then
		the following are equivalent.
        \begin{enumerate}
    \item[(i)] there exists a noncommutative probability space $(\cN,\nu)$
		such that $L_W^p(0,\infty)$ is isomorphic to a
		subspace of $L_r({\cN,\nu})$.
    \item[(ii)] there exists a noncommutative probability space $(\hat{\cN}, \hat{\nu})$ such that $L_W^p(\cM,\tau)$ is isomorphic to a
		subspace of $L_r(\hat{\cN},\hat{\nu})$.
    \item[(iii)] $L_W^p(0,\infty)$ is isomorphic to a subspace of $L_p(0,1)$.
    \item[(iv)] $L_{W}^p(0,\infty)$ has cotype $2$ (is $2$-concave).
    \item[(v)] there is a $C > 0$ such that for all $x \in (0, \infty)$, we have
		\begin{equation*}
			\int_0^x W(t) t^{-p/2} \, dt \leq Cx^{-p/2} \int_0^x W(t) \, dt.
		\end{equation*}
\end{enumerate}
		The implications $(iii)\Leftrightarrow (iv)\Leftrightarrow (v)\Leftrightarrow (i)\Rightarrow(ii)$ hold true when $\cM$ is an arbitrary (not necessarily atomless) semifinite $\sigma$-finite von Neumann algebra equipped with a semifinite (not necessarily infinite) faithful normal trace.
\end{thm}	
        \begin{proof}
		By \cite[Theorem 1]{HJSZ}, statements (i) and (ii) are equivalent. By Proposition~\ref{ifvos}, statements (iii), (iv) and (v) are equivalent. 
        
		(i) $\Rightarrow$ (iv):
        Recall that $L_r(\cN,\nu)$ has cotype $2$\cite{Fack,PisierXu2003}.
        If $L_W^p(0,\infty)\hookrightarrow L_r(\cN,\nu)$, then $L_W^p(0,\infty)$ has cotype $2$.
        
		(iii) $\Rightarrow$ (i): If $L_W^p(0,\infty)\hookrightarrow L_p(0,1)$, then $$L_W^p(0,\infty)\hookrightarrow L_p(0,1)\stackrel{\mbox{\tiny \cite[Prop. 11.1.9]{AK}}}{\hookrightarrow} L_r(0,1) .$$
         This completes the proof. 
		\end{proof}
	As an application, we obtain the following corollary.
		\begin{corollary}\label{main-result1}Let $1\le q\le p<2$ and $ 1\leq r<\infty$.
			For any  semifinite $\sigma$-finite von Neumann algebra $\cM$ equipped with a semifinite infinite faithful normal trace $\tau$, there exists a finite von Neumann algebra $\cN$ equipped with a faithful normal tracial state $\nu$ such that 
			$$L_{p,q}(\cM,\tau)\hookrightarrow L_r(\cN,\nu)$$
			if and only if $r\le q$ and $r\not=p$.  
		\end{corollary}
		\begin{proof}
					Let $W(s):=\frac{q}{p}s^{q/p-1}$, which is nonincreasing, positive and continuous with $\int_0^{\infty} W(s)\,ds=\infty$.
			In particular, $L_W^q(0,\infty)=L_{p,q}(0,\infty)$, which has cotype $2$\cite{Creekmore}.
			
			$(\Rightarrow)$: Assume that $L_{p,q}(\cM,\tau)\hookrightarrow L_r(\cN,\nu)$.
            Since $\ell_{p,q}\hookrightarrow L_{p,q}(\cM,\tau)$ (see e.g. \cite[Propositions 3.0.1 and 3.0.2]{CK2017}), it follows that   $\ell_{p,q} \hookrightarrow L_r(\cN,\nu)$. 
            By \cite{Fack}, $L_r(\cN,\nu)$ has type $\min\{r,2\}$.
            By \cite[Theorem 3.5]{Creekmore}, we have $\min\{r,2\}\leq \min\{p, q\}=q$ which means $1\le r\le q<2$. 
            By \cite[Corollary 1.2]{HRS}, $r\not=p$.

			$(\Leftarrow)$: Assume that  $r\leq q$ and $r\not=p$. Note that  $L_{p,q}(0,\infty )$ has cotype 
            $2$\cite{Creekmore}.
			By (ii) in Theorem~\ref{ifnv}, there exists some noncommutative probability space $(\cN, \nu)$ such that $L_{p,q}(\cM,\tau)$ is isomorphic to a subspace of $L_r(\cN,\nu)$.
            This completes the proof. 
		\end{proof}

    Below, we give the proof of Theorem~\ref{nvos}.
    Note that 
     the function  $W$ in  Theorem~\ref{nvos} is defined only on $(0,1)$, and we may extend  $W$   to a function on  $(0,\infty)$ that satisfies  all assumptions  in Theorem~\ref{ifnv}. 
		
%	\begin{theorem}
%	    \label{nvos}
%		Let $1 \leq r\leq p < 2$ and $ W $ be a nonincreasing, continuous, nonnegative function on $(0, 1)$ with $\int_0^1 W(t)\,dt=1$. Let $\cM$ be an atomless finite von Neumann algebra equipped with a faithful normal trace $\tau$.
%		The following are equivalent.
%\begin{enumerate}
%            \item[(i)] there exists a noncommutative probability space $(\cN, \nu)$ such that $L_W^p(0,1)$ is isomorphic to a subspace of $L_r({\cN,\nu})$.
%		\item[(ii)] there exists a noncommutative probability space $(\hat{\cN},\hat{\nu})$ such that $L_W^p(\cM,\tau)$ is isomorphic to a subspace of $L_r(\hat{\cN},\hat{\nu})$.
%		\item[(iii)] $L_W^p(0,1)$ is isomorphic to a subspace of $L_p(0,1)$.
%		\item[(iv)] $L_{W}^p$ has cotype 2 (is 2-concave). 
%		\item[(v)] there is a $ C > 0$ such that for all $x \in [0, 1]$, we have
%            $$
%			\int_0^x W(t) t^{-p/2} \, dt \leq Cx^{-p/2} \int_0^x W(t) \, dt.
%		$$
%\end{enumerate}
%		The implication $(i)\Rightarrow(ii)$ holds true when $\cM$ is an arbitrary (not necessarily atomless) finite von Neumann algebra equipped with a faithful normal trace.
%	\end{theorem}
\begin{proof}[Proof of Theorem~\ref{nvos}]
Since $(i), (ii)\Rightarrow (iv)$ and $(iii)\Leftrightarrow (iv)\Leftrightarrow (v)$\cite{Schutt}, it suffices to show that $(v)\Rightarrow (i), (ii)$. 
Assume that $W$ satisfies Condition (v).

		Without loss of generality, we assume that $W(1)>0$.
	Indeed,   setting  $W_1(t):=\frac{W(t)+1}{2}$, $t\in (0,1)$, we have  $\int_0^1 W_1(t)dt =1$. 
		By the definition of  a weighted Lorentz space, $L_W^p(0,1)$ and $L_{W_1}^p(0,1)$ coincide, and their norms are equivalent. 
        Moreover,  since  $W$ satisfies Condition $(v)$ of Theorem \ref{nvos}, it follows that $W_1$ possesses the same property.
        That is, for $x \in (0,1)$, there exists a constant $C_1>0$ such that
		\begin{align}\label{finitecase}
			\int_{0}^{x}W_1(t) t^{-p/2} \,dt \leq C_1x^{-p/2} \int_{0}^{x}W_1(t) \,dt.
		\end{align}
        
        		Let $A:=\int_{0}^{1}W_1 (t) t^{-p/2} \,dt$. We have  $0< A\leq C_1$. 
For simplicity, in what follows,  we still denote $W_1$ by $W$.

		Let $\frac{p}{2}-1<\alpha<0$. 
        In what follows, we define 
		$$
		\widetilde{W}(t) := \left\{
		\begin{array}{ll}
			W(t) & \text{if } 0 < t \leq 1, \\
			W(1)t^{\alpha} & \text{if } t > 1,
		\end{array}
		\right.
		$$
        which extends $W(t)$ from $(0,1)$ to $(0,\infty)$.
		Note that $\widetilde{W}(t)$ is a nonincreasing, positive, and continuous function on $(0,\infty)$ satisfying $\int_{0}^{\infty}\widetilde{W}(t)\,dt=\infty$.
        
      Now, we prove that $\widetilde{W}(t)$ satisfies the Condition $(v)$ (with respect to    the unbounded interval $(0,\infty )$) in Theorem~\ref{ifnv}. 
        We denote $I(x):=\int_{0}^{x}\widetilde{W}(t) t^{-p/2} \,dt$ and $S(x):=x^{-p/2}\int_{0}^{x}\widetilde{W}(t) \,dt$.
		For $x \in (1,\infty)$,  we have 
		\begin{align*}
			I(x) &=A+W(1)\int_{1}^{x} t^{\alpha-p/2} \,dt\\
			&=A+\frac{W(1)}{1+\alpha-p/2}x^{1+\alpha-p/2}-\frac{W(1)}{1+\alpha-p/2}  
		\end{align*}
		and 
		\begin{align*}
			S(x) &=x^{-p/2}\left(1+\int_{1}^{x}\widetilde{W}(t) \,dt\right)\\
			  &=x^{-p/2}\left(1+W(1)\int_{1}^{x}t^{\alpha} \,dt\right)\\
			  &=x^{-p/2}+\frac{W(1)}{1+\alpha}x^{1+\alpha-p/2}-\frac{W(1)}{1+\alpha}x^{-p/2}. 
		\end{align*}
		Let 
		$$
		G(x):=\frac{I(x)}{S(x)}, \quad \forall x\ge 1.
		$$
       Note that  $G(x)$ satisfies the following properties:
        \begin{enumerate}
            \item[(1)] $G(x)$ is continuous and positive on $[1,\infty)$,
	    \item[(2)]   $G(1)=A\leq C_1$,
		\item[(3)]  
			\begin{align*}
			\lim_{x \to \infty} G(x)
			=
			\frac{\frac{W(1)}{1+\alpha-p/2}}{\frac{W(1)}{1+\alpha}}
			=
			\frac{1+\alpha}{1+\alpha-p/2}<\infty . 
		\end{align*}  
        \end{enumerate}
		Therefore, $G(x)$ is bounded on $[1,\infty)$.
		There exists a constant $M>0$ such that $0<G(x)\leq M$ for any $x\in [1,\infty)$.
		That is, 
		\begin{align}\label{infinitecase}
		\int_{0}^{x}\widetilde{W}(t) t^{-p/2} \,dt=I(x) \leq MS(x)=Mx^{-p/2} \int_{0}^{x}\widetilde{W}(t) \,dt,\quad x\in [1,\infty).
		\end{align}	
		Let  $\widetilde{C}:=\max\{C_1,M\}$.
        By inequalities (\ref{finitecase}) and (\ref{infinitecase}), we establish the  inequality 	
		\begin{align*}
			\int_{0}^{x}\widetilde{W}(t) t^{-p/2} \,dt \leq \widetilde{C}x^{-p/2} \int_{0}^{x}\widetilde{W}(t) \,dt,\quad x \in (0,\infty),
		\end{align*}
        which satisfies the Condition $(v)$ in Theorem~\ref{ifnv}.
        
        Recalling that $1\le r\le p<2$ and $\frac{p}{2}-1<\alpha<0$, we have $0<\frac{\alpha+1}{p}<\frac{1}{p}\le \frac{1}{r}$.
        Let $\beta$ such that $-\frac{1}{r}<\beta~<-\frac{\alpha+1}{p}~<0$.
        Therefore, we have $0<r\beta+1<\infty$ and $-\infty<p\beta+\alpha+1<0$, which means that 
        \begin{align}\label{pwrne1}
        \int_1^{\infty} t^{r\beta} \,dt =\frac{1}{r\beta+1}\int_1^{\infty} \,d(t^{r\beta+1}) =\infty
    \end{align}
        and 
        \begin{align}\label{pwrne2}
        \int_1^{\infty} t^{p\beta} \widetilde{W}(t)\,dt= W(1)\int_1^{\infty} t^{p\beta+\alpha}\,dt=W(1)\frac{1}{p\beta+\alpha+1}\int_1^{\infty} \,d(t^{p\beta+\alpha+1}) =-\frac{W(1)}{p\beta+\alpha+1}<\infty.
                \end{align}
        Let $$
		f(t) := \left\{
		\begin{array}{ll}
			0 & \text{if } 0 < t \leq 1, \\
			t^{\beta} & \text{if } t > 1.
		\end{array}
		\right.
		$$
        By (\ref{pwrne1}) and (\ref{pwrne2}), we have $f(t)\not\in L_r(0,\infty)\cap L_{\infty}(0,\infty)$ and $f(t)\in L_{\widetilde{W}}^p(0,\infty)\cap L_{\infty}(0,\infty)$, which means that $L_{\widetilde{W}}^p(0,\infty)\cap L_{\infty}(0,\infty) \not= L_r(0,\infty)\cap L_{\infty}(0,\infty)$.
        By Theorem \ref{ifnv},    there exists a noncommutative probability space $(\cN,\nu)$
		such that $L_{\widetilde{W}}^p(0,\infty)\hookrightarrow L_r({\cN,\nu})$ and 
  a noncommutative probability space $(\hat{\cN}, \hat{\nu})$ such that $L_{\widetilde{W}}^p(\cM ,\tau ) \hookrightarrow L_r(\hat{\cN},\hat{\nu})$.
    The proof is competed by noting  that $L_W^p(0,1)\hookrightarrow L_{\widetilde{W}}^p(0,\infty )$  and $L_W^p (\cM,\tau) =  L_{\widetilde{W}}^p(\cM ,\tau )$. 
    \end{proof}
 
 \begin{rem}Figiel, Johnson and Tzafriri  showed that if $X$ is a complemented subspace of $L_W^p(0,1)$ for $1\le p<\infty$, then $X$ either embeds isomorphically into $L_p(0,1)$ (precisely, $X$ is isomorphic to a Hilbert space) or contains a complemented copy of $\ell_p$ \cite[Theorem 5.9]{FJT1975} (see \cite[Theorem 3.1]{Rand03} and \cite{HS26} for its noncommutative extensions).
 
 \end{rem}

	As an application of Theorem~\ref{nvos}, we obtain the following corollary which should be compared with Corollary~\ref{main-result1}.
	
	\begin{corollary}Let $1\le q\le p<2$ and $1\leq r<\infty$.
		For any finite von Neumann algebra $\cM$ on a separable Hilbert space equipped with a faithful normal trace $\tau$, there exists a finite von Neumann algebra $\cN$ equipped with a faithful normal tracial state $\nu$ such that 
		$$L_{p,q}(\cM,\tau)\hookrightarrow L_r(\cN,\nu)$$
		if and only if $r\le q$.  
	\end{corollary}
\subsection{Some remarks concerning $\ell_{p,q}$-sequence spaces isomorphic to subspaces of  noncommutative $L_r$-spaces}

Recall that  the authors in \cite{HSSZ} established the lack of isomorphic embeddings of the sequence space $\ell_{p,q}$ into the corresponding noncommutative Lorentz space $L_{p,q}(\mathcal{M}, \tau)$ over a noncommutative probability space $(\cM,\tau)$  when $p \neq q$.
We study the    isomorphic embeddings of $\ell_{p,q}$ into the noncommutative  space $L_r(\mathcal{N}, \nu)$ over a noncommutative probability space $(\cN,\nu)$, which   complements the  result in \cite{HSSZ}. 

\begin{theorem}\label{lrqintoLr}Let $1\le r<2$ such that $1\le q\le r$ or $1<r<q <\infty$. 
		If $q\not=r$, then 
		$$\ell_{r,q}\not\hookrightarrow L_r(\cN,\nu)$$
		for any finite von Neumann algebra $\cN$ equipped with a faithful normal tracial state $\nu$.
	\end{theorem}
\begin{proof}
    We assume that $q\not=r$  and  there exists a von Neumann algebra $\cN$ such that $\ell_{r,q}\hookrightarrow L_r(\cN,\nu)$.
    By \cite{CF}, $\ell_{r,q}$ contains uniformly the spaces $\ell_r^n$, $n\ge 1$.
    By \cite[Theorem 5.1]{RX}, $\ell_r\hookrightarrow \ell_{r,q}$.
    Therefore, $\ell_r\hookrightarrow L_{r,q}(0,\infty)$.
    By \cite[Theorem 7]{Dilworth}, we have either $r=2$ or $r=q$.
    This is a contradiction.
\end{proof}

\begin{rem}\label{quesanws1}
Let $1\le r<2$ such that $1\le q\le r$ or $1<r<q <\infty$. 
     Note that $\ell_{r,q}$ is the commutative core (see (\ref{commucore}) above) of the ideal (of compact operators) $C_{r,q}$.
     
     If  $C_{r,q}\hookrightarrow L_r(\cN,\nu)$ for some finite von Neumann algebra $(\cN,\nu)$, then we have $\ell_{r,q}\hookrightarrow C_{r,q}\hookrightarrow L_r(\cN,\nu)$.
     By Theorem \ref{lrqintoLr} above, we have $r=q$ (in this case, $\ell_{r,q}=\ell_{r}\hookrightarrow L_r(0,1)$).
           This answers the Question \ref{JSZquestion} in the setting of $C_{r,q}$ when $1\le r<2$ such that $1\le q\le r$ or $1<r<q <\infty$.
           Moreover, the following conditions are equivalent ($ C_{r,q}\hookrightarrow L_r(\cN,\nu) \Rightarrow  (1)\stackrel{Th. \ref{lrqintoLr}}{\Leftrightarrow} (2)\stackrel{Th. \ref{lrqintoLr}}{\Leftrightarrow}  (3)  $)
\begin{enumerate}
    \item[(1)]   $\ell_{r,q}\hookrightarrow L_r(\cN,\nu)$ for some noncommutative probability space $(\cN,\nu)$;
    \item[(2)]  $r=q$;
    \item[(3)]  $\ell_{r,q}\hookrightarrow L_r(0,1)$. 
\end{enumerate}

\end{rem}

The classification of   isomorphic embedding   of $\ell_{p,r}$ into a noncommutative $L_r$-space follows from   their type and cotype.
	
	\begin{prop}\label{proplprintoLr}
Let $1\le r<2$ such that $1\le r\le p<\infty$ or $1<p<r$.
		There exists some finite noncommutative probability space $(\cN,\nu)$ such that
		$$\ell_{p,r}\hookrightarrow L_r(\cN,\nu)$$ 
		if and only if
		$1\le r\le p<2$
	\end{prop}
	\begin{proof}
		$(\Rightarrow)$: 
        Assume that 
         $\ell_{p,r}\hookrightarrow L_r(\cN,\nu)$.
         Recall that (see  \cite[Theorems 3.5 and  3.6]{Creekmore} and  \cite{Fack}):
		\begin{itemize}
			\item if $1< p<r<2$, $\ell_{p,r}$ has type $p$ and $L_r(\cN,\nu)$ has type $r$,
			\item if $p>2$, $\ell_{p,r}$ has cotype $p$ and $L_r(\cN,\nu)$ has cotype $2$,
			\item if $p=2$, $\ell_{p,r}$ has cotype $(2+\epsilon)$ for any $\epsilon>0$ and $L_r(\cN,\nu)$ has cotype $2$.
		\end{itemize}
 For the above three cases,  $\ell_{p,r}\not \hookrightarrow L_r(\cN,\nu)$. Therefore, 
		we obtain  $1\le r\le p<2$.
		
		$(\Leftarrow)$: If $1\le r\le p<2$, then by Corollary \ref{aoifvos}, we have $L_{p,r}(0,\infty)\hookrightarrow L_r(0,1)$.
        Therefore, we have 
        \begin{align}\label{lprintoLr}
        \ell_{p,r}\hookrightarrow L_{p,r}(0,\infty)\hookrightarrow L_r(0,1). 
        \end{align}
        The proof is complete by taking $\cN =L_\infty (0,1)$.
	\end{proof}
    
        \begin{rem}\label{quesanws2}
    %    Note that $\ell_{p,q}$ is the commutative core (see (\ref{commucore}) above) of the ideal (of compact operators) $C_{p,q}$.

     %  Let $1\le p<2$ such that $(p,q)\in \{(1,1) \} \cup (1,2)\times (1,\infty )  $. 
%The following conditions are equivalent:
%\begin{enumerate}
%    \item[(1)] $\ell_{p,q}\hookrightarrow L_p(\cN,\nu)$ for some noncommutative probablity space $(\cN,\nu)$;
%    \item[(2)]  $p=q$;
%    \item[(3)]  $\ell_{p,q}\hookrightarrow L_p(0,1)$. 
%\end{enumerate}
%      The implications $(1),(3)\Rightarrow(2)$ follow from Theorem \ref{lrqintoLr} above. 
%      The implications $(2)\Rightarrow(1),(3)$ are obvious. 
%      This answers Question \ref{JSZquestion} in the setting of $C_{p,q}$ when $(p,q)\in \{(1,1) \} \cup (1,2)\times (1,\infty )  $. 
Let $1\le r<2$ such that $1\le r\le p<\infty$ or $1<p<r$.
Note that $\ell_{p,r}$ is the commutative core (see (\ref{commucore}) above) of the ideal (of compact operators) $C_{p,r}$.

The following conditions are equivalent:
\begin{enumerate}
    \item[(1)] $\ell_{p,r}\hookrightarrow L_r(\cN,\nu)$ for some noncommutative probability space $(\cN,\nu)$;
    \item[(2)]  $1\le r\le p<2$;
    \item[(3)]  $\ell_{p,r}\hookrightarrow L_r(0,1)$. 
\end{enumerate}
The equivalence $(1)\Leftrightarrow(2)$ follows from Proposition \ref{proplprintoLr} above. 
The equivalence  $(2)\Leftrightarrow(3)$ 
follows from the type and cotype of $\ell_{p,r}$ and $L_r(0,1)$ \cite[Theorems 3.5 and  3.6]{Creekmore}  and  \eqref{lprintoLr} above. 

If $C_{p,r}\hookrightarrow L_r(\cN,\nu)$ for some noncommutative probability space $(\cN,\nu)$, then we have $$\ell_{p,r} \hookrightarrow C_{p,r} \hookrightarrow L_r(\cN,\nu).$$
Therefore, we have  $\ell_{p,r}\hookrightarrow L_r(0,1)$.
This  answers  Question \ref{JSZquestion} in the affirmative for
$C_{p,r} $ when $1\le r<2$ such that $1\le r\le p<\infty$ or $1<p<r$. 

%By Proposition \ref{proplprintoLr}, 
 %       $\ell_{p,r}\hookrightarrow L_r(\cN,\nu) $ if and only if $1\le r\le p<2$. 
%On the other hand,   by the type and cotype of $\ell_{p,r}$ and $L_r(0,1)$ \cite[Theorems 3.5 and  3.6]{Creekmore}  and  Corollary~\ref{aoifvos}, $\ell_{p,r}\hookrightarrow L_r(0,1)$ if and only if $1\le r\le p <2$ .  
        \end{rem}

        The case when $2<r<\infty$ is completely different from that of $1\le r<2$. 
	
	\begin{proposition}
		Let $2<r<\infty$.
		There exists some noncommutative probability space $(\cN,\nu)$ such that
		$$\ell_{p,q}\hookrightarrow L_r(\cN,\nu)$$
		if and only if $p=q=2$ or $p=q=r$. 
	\end{proposition}
	\begin{proof}
		$(\Rightarrow)$: %By type and cotype argument above, we obtain that $2<p\le r$ and $2\le q\le r$.
	It follows immediately from the noncommutative Kadec--Pelczynski dichotomy theorem  \cite[Proposition 5.4]{RX}	  
    (see also \cite[Theorem 4.1]{HSSZ}, \cite[Corollary 6]{KP62} and \cite{S96}). 
		%By ultrapower technique~\cite[Chapter 11]{AK}, we obtain that $L_{p,q}(0,1)\hookrightarrow L_r(0,1)$.
		%Since $r>2$, by \cite[Theorem 1.6]{JMST} (see also \cite[Theorem 10.8]{Pisier}), we obtain that $p=q=r$ or $p=q=2$.
		
		$(\Leftarrow)$: If $p=q=2$ or $p=q=r$, then the statement follows from the fact that $\ell_2, \ell_r\hookrightarrow L_r(0,1)$  \cite[Theorem 6.4.19]{AK}.
	\end{proof}

	\begin{rem}
	Let $\cM$ be  a  semifinite infinite von Neumann algebra equipped with a semifinite faithful normal trace $\tau$ and let $2\le r<\infty$.
	By \cite[Theorem 3.1]{S96}, if $(p,q)\not=(2,2)$, then 
	$L_{p,q}(\cM,\tau)$ does not embed into $L_r(\cN,\nu)$  (in particular, $C_{p,q}$)
	for any noncommutative probability space $(\cN,\nu)$.
	\end{rem}

 It was shown in \cite{Schutt}  that $\ell_{p,q}$, $1<p<q $, does not embed into $L_1(0,1)$ via   a deep result due to Kwapien and Sch\"{u}tt\cite{KS89}. 
 Let $X$ be a subspace of $L_p(0,1),\enspace 1\le p< 2$.
  Recall that Rosenthal \cite[Theorem $2'$]{Rosenthal75}  states that 
\begin{enumerate}
 \item  if $p=1$   \begin{enumerate}
        \item and  $X$ is not reflexive, then $X$  contains a complemented copy  of $\ell_p$.
        \item and  $X$ is reflexive,  then there is a $1< p'<2 $ such that $X\hookrightarrow L_{p'}(0,1)$; 
    \end{enumerate}
\item     if $1<p<2$
      \begin{enumerate}
        \item and   the $p$-norm is not equivalent to   $1$-norm  on $X$, then $X$  contains a complemented  copy   of $\ell_p$.
        \item and  the $p$-norm is   equivalent to   $1$-norm  on $X$, then  there is a $p'>p$ such that $X\hookrightarrow L_{p'}(0,1)$.
    \end{enumerate}
\end{enumerate}
%if $p=1$ and $X$ is not reflexive or $p>1$ and the $p$-norm is not equivalent with   $1$-norms  on $X$, 
%then $X$  contains a complemented isomorph of $\ell_p$.
%On the other hand, 
%if $X$ is reflexive in case $p=1$ or if the $p$-norm and $1$-norms on $X$ are equivalent for the  case $p>1$, then there is a $p'>p$ such that $X\hookrightarrow L_{p'}(0,1)$.
 In particular, if $p=1$ and $X$ is reflexive, then there exists some $p'>1$ such that $X\hookrightarrow L_{p'}(0,1)$ and the set of such $p'$ is an open interval \cite[p. 774]{Rosenthal75} (see also \cite{Rosenthal73}). %Let $q$ be the right-end point of this interval.
%{\color{red}(Wrong)In particular, if there exists some $1\le p< 2$ such that $X\hookrightarrow L_p(0,1)$, then there exists some $p<p'$ such that $X\hookrightarrow L_{p'}(0,1)$.
%That is, the set of $1\le p<2$ such that $X\hookrightarrow L_p(0,1)$ forms an open interval \cite[p. 774]{Rosenthal75} (see also \cite{Rosenthal73}).
%Moreover, the interval is nonempty provided $X$ is reflexive.
%Let $q$ be the right-end point of this interval.
%} 
Precisely, Rosenthal introduced the index $I_p(X)$, which equals the infimum of those numbers $K$ so that there is a probability density $\varphi$ on $(0,1)$ (i.e. $\int_0^1 \varphi(t)dt=1$ and $\varphi>0$) with
$$
\left(\int_0^1|x(t)|^p \varphi^{1-p}(t)\,dt\right)^{1/p}\le K\int_0^1|x(t)|\,dt,\enspace \text{for all} \enspace x\in X.
$$
Rosenthal \cite{Rosenthal73,Rosenthal75} showed that for $1<p\le 2$, $I_p(X)<\infty$ if and only if $X$ is isomorphic to a subspace $Y$ of $L_p(\nu)$ (for some measure $\nu$) such that the $p$-norm and $1$-norm are equivalent on $Y$.
Rosenthal proved that the set of $p \in [1, 2)$ such that $I_p(X) < \infty$ forms an open interval \cite[p. 774]{Rosenthal75}. 
Let $q$ be the right-end point of this interval.
Moreover, this interval is nonempty provided that  $X$ is reflexive; if $X \subset L_1$ is reflexive, the interval is of the form $(1, q)$.

Rosenthal posed the following natural question:
\begin{quest}  Let $X$ be a subspace of $L_p(0,1),\enspace 1\le p< 2$.
    Does $X$ contain a subspace isomorphic to $\ell_r$ for some $r$, $1< r\le2$?
\end{quest}
Rosenthal was  able to show that if $X$ contains a subspace that embeds into $L_q(0,1)$ (where $q$ is the right-end point of the interval formed by the set of $p\in [1,2)$ such that $I_p(X)<\infty$), then the above question has an affirmative answer.
%That is, if $X$ is  a subspace of $L_p(0,1),\enspace 1\le p< 2$, then $X$  contains a complemented copy of $\ell_r$ for some $1< r\le2$. 
%there exists some $r$ for $1\le r<2$ such that $\ell_r\hookrightarrow X$, which can answer the following question affirmatively:
Guerre-Delabri\`{e}re and Levy  \cite[Theorem IV.5.5]{GD92}  provided an affirmative answer to this question in full generality  by showing that   for all $\epsilon>0$, $\ell_{p(X)}$ is $(1+\epsilon)$-isomorphic to a subspace of $X$, where $p(X)$ is the Rademacher type of $X$.
Note that 
$p(X)$ is the smallest real number $p$
such that $\ell_p\hookrightarrow X$.
%Therefore, $\ell_{p(X)}\hookrightarrow L_p(0,1)$.
The proof of the theorem is based on the stability of $L_p(0,1)$.
In the noncommutative setting, Nhany showed  that $L_p(\cM)$ is stable, $1\le p<\infty$ if and only if the underlying von Neumann algebra $\cM$ is of type I \cite[La Th\'eor\`{e}me principale]{Nhany}. 
Therefore, it remains unknown whether a noncommutative version of Levy--Guerre-Delabri\`{e}re Theorem holds for noncommutative $L_p$-spaces, $1\le p< 2$. 
Carothers and Dilworth 
 \cite[p. 76]{CD91} provided  another argument for $$\ell_{p,q}\not\hookrightarrow L_1(0,1) $$ 
 when $q>p>1$ by the Guerre-Delabri\`{e}re--Levy  theorem. 
Since both Kwapien--Sch\"{u}tt theorem and  Guerre-Delabri\`{e}re--Levy  theorem are not available in the noncommutative setting, we do not know whether $\ell_{p,q}$ embeds into the predual of a von Neumann algebra.

\section{Isomorphic embedding of $C
_{p_1,q_1}$ into $C_{p_2,q_2}$ }
If $\{e_n\}_{n\ge1}$ is a basis of a $p$-Banach space $X$, a sequence $\{y_n\}_{n\ge 1}$ in $X$ is said to be {\it a block basic sequence of $\{e_n\}_{n\ge1}$ }if there is an increasing sequence of natural numbers $\{p_n\}_{n\ge1}$ such that $y_n=\sum_{i=p_n+1}^{p_{n+1}} \alpha_i e_i$, where  $\alpha_i $'s are  scalars.

\begin{proposition}\label{lemmaporq}
Let $\{y_n\}_{n\ge 1} $ be a semi-normalized block basis of the vector basis $\{e^{\ell_{p,q}}_n \}_{n\ge 1}$ of $\ell_{p,q}$. 
There exists a subsequence  $\{y_{k_n}\}_{n\ge 1} $  of $\{y_n\}_{n\ge 1} $ such that 
$$\norm{\sum_{n=1 }^N y_{k_n}}_{\ell_{p,q}}\approx N^{1/p} $$
or 
$$\norm{\sum_{n=1 }^N y_{k_n}}_{\ell_{p,q}}\approx N^{1/q}.  $$
\end{proposition}
\begin{proof}It suffices to consider the case when $p\ne q$. 
Let $$y_n:=\sum_{i= k_{n}+1 }^{k_{n+1 } } \alpha_{  i} e_{ i}^{\ell_{p,q}} .$$ 
Without loss of generality, we may assume that $\{|\alpha_{  i}| \}_{i= k_{n}+1 }^{k_{n+1 } } $ is nonincreasing. 
Let $\alpha_{(n)} : = (\alpha_{  i})_{i= k_{n}+1 }^{k_{n+1 } }$.

Let $\beta$ be a limit point of $(\alpha_{(n)})_{n\ge 1}$ in the pointwise topology. 
If $\beta=0$, that is, $\norm{ \alpha_{(n)}}_\infty \to 0 $ as $n\to \infty $, then by
Lemma \ref{Lpqequivalp}, passing to a subsequence if necessary, we obtain that $$\norm{\sum_{n=1 }^N y_{k_n}}_{\ell_{p,q}}\approx N^{1/q}.  $$

Now,   assume that $\beta=\lim_{n\to \infty }  \alpha_{(n )}  \ne 0$. 
By the Fatou property of $\ell_{p,q}$, we have $\beta\in \ell_{p,q}$ with $\norm{\beta}_{\ell_{p,q}} \le  \sup_{n\ge 1}\norm{\alpha_{(n)}}_{\ell_{p,q}}$\cite[Lemma 3.6]{NP}. Note that $\norm{ \alpha_{(n)}  - \beta}_\infty  \to 0$ as $n\to \infty $, and  
$$\norm{\sum_{k=1}^n y_k}_{\ell_{p,q}} = \norm{\oplus _{k=1}^n \alpha_{(k)}}_{\ell_{p,q}} = \norm{\oplus _{k=1}^n \left( \beta - (\beta  - \alpha_{(k)}) \right) }_{\ell_{p,q}}. $$
We have $$\left|\norm{\oplus _{k=1}^n  \beta   }_{\ell_{p,q}} -  \norm{\oplus _{k=1}^n \left(  \beta  - \alpha_{(k)}  \right) }_{\ell_{p,q}}\right| \lesssim   \norm{\sum_{k=1}^n y_k}_{\ell_{p,q}}  \lesssim \norm{\oplus _{k=1}^n  \beta   }_{\ell_{p,q}} +  \norm{\oplus _{k=1}^n \left(  \beta  - \alpha_{(k)}  \right) }_{\ell_{p,q}} . $$
If $ \norm{ \left(  \beta  - \alpha_{(k)}  \right) }_{\ell_{p,q}}\to 0$ as $k\to \infty$, then  
passing to a subsequence if necessary, we have  $$\norm{\sum_{n=1 }^N y_{k_n}}_{\ell_{p,q}}\approx N^{1/p}.  $$
If $ \norm{ \left(  \beta  - \alpha_{(k)}  \right) }_{\ell_{p,q}}\not \to 0$ as $k\to \infty$, then 
passing to a subsequence if necessary, we have  
$$\left|N^{1/p}\norm{  \beta   }_{\ell_{p,q}} -  N^{1/q}\right|  \lesssim \norm{\sum_{k=1}^N y_k}_{\ell_{p,q}}  \lesssim N^{1/p}\norm{  \beta   }_{\ell_{p,q}} +  N^{1/q} . $$
If $p<q$, then  
we have  
$$ \norm{\sum_{k=1}^N y_k}_{\ell_{p,q}} \approx  N^{1/p}  . $$
If $p>q$, then we have 
$$ \norm{\sum_{k=1}^N y_k}_{\ell_{p,q}} \approx  N^{1/q}    . $$
This completes the proof. 
\end{proof}

The following result is well-known in the setting of Banach spaces, see e.g. \cite[Theorem 4.e.5]{LT1}. 
Below, we establish this result in the setting of quasi-Banach spaces $\ell_{p,q}$. 
\begin{lemma}\label{classical}
Let $0<p_1,p_2,q_1,q_2 <\infty $. 
If $\ell_{p_1,q_1}\hookrightarrow \ell_{p_2,q_2}$, then 
\begin{enumerate}
\item either $p_1=q_1=q_2$,
\item or $p_1\ne q_1$,  $p_1=p_2$ and $q_1=q_2$.
\end{enumerate}
\end{lemma}
\begin{proof}
By \cite[Lemma 2.1]{KR}, the standard basis of $\ell_{p_1,q_1}$ is equivalent to a block basis $(y_n)_{n \ge 1 }$ in $\ell_{p_2,q_2}$.  
Let $$y_n:=\sum_{i= k_{n}+1 }^{k_{n+1 } } \alpha_{  i} e_{ i}^{\ell_{p_2,q_2}} .$$ Since $\ell_{p_2,q_2 }$ is symmetric, we may assume that $\alpha_{(n)}: =(\alpha_{ i})_{i= k_{n}+1 }^{k_{n+1 } }$ is decreasing. 
%By Lemma \ref{Lpqequivalp}, 
%we obtain that $\left( \frac{1}{2^{n/p_1} } \sum_{i=2^n +1}^{2^{n+1}} e_i ^{\ell_{p_1,q_1}} \right)_{n\ge 1}$ contains a subsequence which is equivalent to the vector basis of $\ell_{q_1}$, and, passing to a further subsequence, 
%$\left( \frac{1}{2^{n/p_1} }T( \sum_{i=2^n +1}^{2^{n+1}} e_i ^{\ell_{p_1,q_1}}) \right)_{n\ge 1}$ contains a subsequence which is equivalent to the vector basis of $\ell_{q_2}$. 
%Since $T$ is an isomorphism, 
%we  have  $q_1=q_2$. 

Let  $\overline{\text{span}\{y_n\}}$ be  the  closure of $\text{span}\{y_n\}$ in  $\ell_{p_2,q_2}$ and 
let $T: \ell_{p_1,q_1} \rightarrow \overline{\text{span}\{y_n\}}$ be an isomorphism such that $T(e_n^{\ell_{p_1,q_1}}) =y_n$ for every $n\ge 1$.
By Lemma~\ref{Lpqequivalp}, 
we obtain that $\left( \frac{\sum_{i=2^n +1}^{2^{n+1}} e_i ^{\ell_{p_1,q_1}}}{\norm{\sum_{i=2^n +1}^{2^{n+1}} e_i ^{\ell_{p_1,q_1}}}_{\ell_{p_1,q_1}} } \right)_{n\ge 1}$ contains a subsequence which is equivalent to the vector basis of $\ell_{q_1}$.
Since $T$ is an isomorphism, we have $$\norm{\sum_{i=2^n +1}^{2^{n+1}} e_i ^{\ell_{p_1,q_1}}}_{\ell_{p_1,q_1}} \approx \norm{T( \sum_{i=2^n +1}^{2^{n+1}} e_i ^{\ell_{p_1,q_1}})}_{\ell_{p_2,q_2}}.$$
Passing to a further subsequence, by Lemma \ref{Lpqequivalp} again, 
$\left( \frac{T( \sum_{i=2^n +1}^{2^{n+1}} e_i ^{\ell_{p_1,q_1}})}{\norm{T( \sum_{i=2^n +1}^{2^{n+1}} e_i ^{\ell_{p_1,q_1}})}_{\ell_{p_2,q_2}}} \right)_{n\ge 1}$ (which is equivalent to $
\left( \frac{  \sum_{i=2^n +1}^{2^{n+1}} y_i }{\norm{  \sum_{i=2^n +1}^{2^{n+1}} y_i }_{\ell_{p_2,q_2}}} \right)_{n\ge 1}$)
 contains a subsequence which is equivalent to the vector basis of $\ell_{q_2}$. 
Then, 
we  have  $q_1=q_2$ (see \cite[Proposition 2.9]{KPR} above). 

Since $\{e_n^{\ell_{p_1,q_1}}\}_{n\ge 1} \sim \{y_n\}_{n\ge 1}$, we have $$\norm{\sum_{i=1}^n e_i^{\ell_{p_1,q_1}}}_{\ell_{p_1,q_1}} \approx \norm{\sum_{i=1}^n y_i}_{\ell_{p_2,q_2}}.$$
By Proposition~\ref{lemmaporq}, we obtain that $p_1=p_2$  or $q_2$. 
In the first case, we have $p_1=p_2$, $q_1=q_2$,  and in the second case, we have $p_1=q_2 =q_1$. This completes the proof. 
\end{proof}

The direct sum   $X\oplus  Y$ of two Banach spaces $X$ and   $Y$, consists of all elements  $(x,y)$ equipped  with the norm 
$$
\norm{(x,y)}_{X\oplus Y}:= \norm{x} _X+\norm{y} _Y ,~ x\in X, y\in Y.
$$

%By

 It is a  corollary of    \cite[Theorem 5.1]{HSZ25} (see also \cite{HSY}, and see  \cite[Corollary 2.8]{Arazy81} for the case when $p_1,q_1,p_2,q_2>1$) that 
 if $\ell_{p_1,q_1}\hookrightarrow C_{p_2,q_2}$ for $0<p_1,p_2,q_1,q_2 <\infty $, then $\ell_{p_1,q_1}\hookrightarrow \ell_{p_2,q_2}\oplus \ell_2$.
Having this result at hand,   we characterize these $\ell_{p_1,q_1}$-spaces that embed isomorphically into  $C_{p_2,q_2}$.
 % and \cite[Corollary 2.8]{Arazy81} (consider $\ell_{\max\{p_2,2\},\max\{q_2,2\}}$)
%, we obtain that 
\begin{proposition}\label{directsuml2}
Let $0<p_1,p_2,q_1,q_2 <\infty $.
If $\ell_{p_1,q_1}\hookrightarrow C_{p_2,q_2}$ (in particular,  $\ell_{p_1,q_1}\hookrightarrow \ell_{p_2,q_2}\oplus \ell_2$), then  
\begin{enumerate}
\item $p_1=q_1=2$,
\item or $p_1=q_1= q
_2$, 
\item or $p_1=p_2$, $q_1=q_2$. 
\end{enumerate}

\end{proposition} 
 \begin{proof}

 By    \cite[Theorem 5.1]{HSZ25}, we may assume that $T$ is an isomorphism from $\ell_{p_1,q_1}$ into $ \ell_{p_2,q_2}\oplus \ell_2$ such that  $T(e_k^{\ell_{p_1,q_1}}) \sim a_k \oplus b_k \in \ell_{p_2,q_2}\oplus \ell_2$.  By \cite[Lemma 2.1]{KR07}, we may assume that the images $T(e_k^{\ell_{p_1,q_1}} )$ are disjoint. 
%If $b_k\to 0$ in $\ell_2$ as $k\to \infty$, then the assertion follows from Lemma \ref{classical}. 

(A): If $b_k\to 0$ in $\ell_2$ as $k\to \infty$, then the assertion follows from Lemma \ref{classical}.
  
(B): If $b_k\not\rightarrow 0$ in $\ell_2$ as $k\to \infty $ and $a_k \to 0$ in $\ell_{p_2,q_2 }$ as $k\to \infty $, then we have 
 $$\left\{ T(e_k^{\ell_{p_1,q_1}} )\right\}_{k\ge 1}\sim  \left\{  e_{k}^{\ell_2}\right\}_{k\ge 1}.$$
 In this case, we have  $p_1=q_1 =2$. 
 
(C): If $b_k\not\rightarrow 0$ in $\ell_2$ as $k\to \infty $ and $a_k \not\to 0$ in $\ell_{p_2,q_2 }$ as $k\to \infty $, 
%Now, we assume that $b_k\not\rightarrow 0$ in $\ell_2$ as $k\to \infty $.
%If  $a_k \to 0$ in $\ell_{p_2,q_2 }$ as $k\to \infty $, then we have 
 %$$\left\{ T(e_k^{\ell_{p_1,q_1}} )\right\}_{k\ge 1}\sim  \left\{  e_{k}^{\ell_2}\right\}_{k\ge 1}.$$
 %In this case, we have  $p_1=q_1 =2$. 
%If $a_k \not\to 0$ in $\ell_{p_2,q_2}$, 
then  passing to a subsequence if necessary, 
 we have 
\begin{eqnarray*}
\left\{ e_k^{\ell_{q_1} }\right\}_{k\ge 1}
&\mathrel{\overset{\text{\tiny Lem. \ref{Lpqequivalp}}}{\sim}}
&\left\{ \frac{1}{(2^k)^{1/p_1}} \sum_{n=2^k +1}^{2^{k+1 }}  T(e_n) \right\}_{k\ge 1}\\
&\sim 
&\left\{    \frac{1}{(2^k ) ^{1/p_1}}\sum_{n=2^k +1}^{2^{k+1 }}  a_n \oplus  \frac{1}{(2^k ) ^{1/p_1}} \sum_{n=2^k +1}^{2^{k+1 }}   b_n \right\}_{k\ge 1}
%&\mathrel{\overset{\text{\tiny Lem. \ref{Lpqequivalp}}}{\sim}}
%&\left\{    e_k^{\ell_{q_2}} \oplus  \frac{1}{(2^k ) ^{1/p_1}} \sum_{n=2^k +1}^{2^{k+1 }}   b_n \right\}_{k\ge 1}
%&\mathrel{\overset{\text{\tiny\cite[Lemma 2.1]{Dilworth}}}{\sim}}
%&\left\{e_k^{\ell_{q_2}}+e_k^{\ell_2}\right\}_{k\ge 1}
.
\end{eqnarray*}
Note that 
$   \left\{   \frac{1}{(2^k ) ^{1/p_1}}\sum_{n=2^k +1}^{2^{k+1 }}  a_n    \right\}_{k\ge 1}$ vanishes or 
is equivalent to $\{e_k^{\ell_{q_2}}\}_{k\ge 1}$ (see Lemma~\ref{Lorequivalp}), 
 $\left\{ \frac{1}{(2^k ) ^{1/p_1}} \sum_{n=2^k +1}^{2^{k+1 }}   b_n \right\}_{k\ge 1}$      vanishes or is equivalent to  $\{e_k^{\ell_{2}}\}_{k\ge 1}$.

\begin{enumerate}
    \item[(a)] If $q_2\ge 2$, then we have 
    $$
     \left\{ e_k^{\ell_{q_1} }\right\}_{k\ge 1} \sim  \left\{ e_k^{\ell_{2} }\right\}_{k\ge 1} .
    $$
    \item[(b)] If $q_2<2$, then we have 
      $$
     \left\{ e_k^{\ell_{q_1} }\right\}_{k\ge 1} \sim  \left\{ e_k^{\ell_{q_2} }\right\}_{k\ge 1} .
    $$
\end{enumerate}
This shows that $q_1= q_2$ or $2$. 
 
On the other hand, we have 
$$ n^{1/p_1} = \norm{\sum_{k=1}^n e_k^{\ell_{p_1,q_1} } }_{\ell_{p_1,q_1}} \approx  \norm{\sum_{k=1}^n T\left(e_k^{\ell_{p_1,q_1 } }\right) }_{\ell_{p_2,q_2}\oplus \ell_2 }\approx  
\norm{\sum_{k=1}^n a_k \oplus b_k}_{\ell_{p_2,q_2}\oplus \ell_2} . $$
By Proposition \ref{lemmaporq},  we obtain that
  $$p_1 = p_2, q_2 \mbox{ or } 2 . $$
Therefore, it suffices to consider  the cases when 
$$\ell_{p_1,q_1}=\ell_{2,2},\ell_{q_2,q_2}, \ell_{p_2,q_2} , \ell_{2,q_2},\ell_{p_2, 2}, \ell_{q_2, 2} ,   $$ and show  that $$\ell_{2,q_2},\ell_{p_2, 2}, \ell_{q_2, 2}\not\hookrightarrow \ell_{p_2,q_2}\oplus \ell_2$$  if $q_2\ne 2$ and $ p_2\ne 2$.

Without loss of generality, we may assume that $a_k=\oplus _{i=0}^{k-1} 0 \oplus \mu(a_k )\oplus _{i=k+1}^{\infty } 0 $, $k=1,2,\cdots$. 
Let $a$ be a pointwise limit of $\mu(a_k)$'s. By \cite[Lemma 3.6]{NP}, $a\in \ell_{p_2,q_2}$. 
In particular,    $\mu(a_k)-a \to 0$ in measure (possibly goes to $0$ in norm). 
Passing to a subsequence if necessary, we write 
 $$T(e_k^{\ell_{p_1,q_1}}) \sim   \Big(  \oplus_{i=0}^{k-1} 0 \oplus   (  a+ (\mu(a_k) -a )) \oplus_{i=k+1 }^{\infty} 0 \Big)   \bigoplus\Big( \oplus_{i=0}^{k-1} 0 \oplus   b_k \oplus_{i=k+1 }^{\infty} 0 \Big)  \in \ell_{p_2,q_2}\oplus \ell_2 $$
 where $k=1,2,\cdots$.
\begin{enumerate}
\item  If  
 $\ell_{2,q_2} \hookrightarrow \ell_{p_2,q_2}\oplus \ell_2$, then, passing to a subsequence if necessary,  we have 
\begin{align*}
\{ e^{\ell_{q_2}}_k \}_{k\ge 1 }
&\mathrel{\overset{\text{\tiny Lem. \ref{Lpqequivalp}}}{\sim}}
\left\{ \frac{1}{(2^k )^{1/2}}\sum_{n=2^k+1}^{2^{k+1 } }
(a_n  \oplus b_n )  \right\} _{k\ge 1 }.
\end{align*}
The right-hand-side (passing to a further subsequence if necessary) is equivalent to $\{
 e^{\ell_{q_2}} \oplus   e^{\ell_{2}}\}_{k\ge 1}$ or $\{
   e^{\ell_{2}}\}_{k\ge 1}$. 
   If $q_2\ge 2$, then we have $\{ e^{\ell_{q_2}}_k \}_{k\ge 1 } \sim \{ e^{\ell_{2 }}_k \}_{k\ge 1 } 
   $, i.e.,   $q_2=2$, which is a contradiction to our assumption. 
 Hence, we have 
 $q_2<  2$.
 This together $$n^{1/2}=  \norm{\sum_{k=1}^n  e_k  ^{\ell_{2,q_2}} }_{\ell_{2,q_2 }} \approx  \norm{\sum_{k=1}^{n} \left(  \oplus_{i=0}^{k-1} 0 \oplus   (  a+ (\mu(a_k)-a )) \oplus_{i=k+1 }^{\infty} 0 \right) }_{\ell_{p_2,q_2}}  +    n^{1/2}$$ implies that 
\begin{enumerate}
\item  if  $p_2=q_2$, then we have $n^{1/2} =n^{1/p_2}$, i.e., $2=p_2=q_2$, which is a contradiction; 
\item if $p_2\ne q_2$, then, by $\mu(a_k)-a\to 0$ pointwise,  we have $n^{1/2} \gtrsim    n^{1/ p_2 }\norm{a}_{\ell_{p_2,q_2}}  $.  Hence, we have either $a=0$ or $p_2 >2 $. 
\end{enumerate} 
 
 If 
$\norm{ \mu(a_k)-a}_{\ell_{p_2,q_2 }} \not\to 0$ as $k\to \infty $,  
  then by Lemma \ref{Lpqequivalp}, we have $$\left\{\left(  \oplus_{i=0}^{k-1} 0 \oplus   (\mu(a_k)-a )\oplus_{i=k+1 }^{\infty} 0 \right) \right\}_{k=1}^\infty \sim \{e_k^{\ell_{q_2}}\}_{k=1}^\infty  .$$
   Moreover, passing to a subsequence if necessary, we may assume that $a $ is disjoint with $\mu(a_k)-a$ for $k=1,2,\cdots$ \footnote{Indeed, let $a=\sum_{k=1}^\infty \alpha _k e_k^{\ell_{p_2,q_2}}$. 
   There exists a sufficiently large number $n_1$ such that $\norm{
   a- \sum_{k=1}^{n_1} \alpha _k e_k^{\ell_{p_2,q_2}}}_{\ell_{p_2,q_2}}<\frac12 $. 
   Since $\mu(a_k)-a\to 0$ as $k\to \infty $ in measure, it follows that $ (\mu(a_k)-a ) \sum_{k=1}^{n_1}  e_k^{\ell_{p_2,q_2}}\to0 $ in quasi-norm as $k\to \infty$.
   Hence, there exists a sufficiently large number $k_1$,  such that $\norm{
 (\mu(a_{k_1})-a ) \sum_{k=1 }^{n_1 }  e_k^{\ell_{p_2,q_2} 
   }}_{\ell_{p_2,q_2} }<\frac14$. 
   In particular,   $ (\mu(a_{k_1})-a ) \sum_{k=1 }^{n_1 }  e_k^{\ell_{p_2,q_2}} $  and $\sum_{k=n_1+1}^{\infty } \alpha _k e_k^{\ell_{p_2,q_2}}$ are mutually disjoint. 
   Arguing inductively, we may construct two mutually disjoint sequences which are equivalent to $ \{ \oplus_{i=1}^{k-1} 0 \oplus     a \oplus_{i=k+1 }^{\infty} 0\}_{k\ge 1} $ and $ \{ \oplus_{i=1}^{k-1} 0 \oplus    ( \mu(a_k)-a) \oplus_{i=k+1 }^{\infty} 0\}_{k\ge 1} $ respectively. 
  }. Hence, 
   for any finite sequence $\alpha=(\alpha_k )$, we have 
   $$\norm{\sum_{k=1}^\infty \alpha_k
     e_{k}^{\ell_{2,q_2}}}_{\ell_{2,q_2}}\gtrsim  \norm{\sum_{k=1}^\infty \alpha_k  e  _{k}^{\ell_{ q_2}}  } _{\ell_{q_2 }}  ,$$ which implies that $q_2\ge 2$ (see \cite[p. 217]{BS}). 
This is a contradiction with the fact that $q_2<2$. 
Hence, we may assume that 
$\norm{ \mu(a_k)-a}_{\ell_{p_2,q_2 }} \to 0$ as $k\to \infty $.  
If $a=0$, then we have 
$\{ e_{k}^{\ell_{2,q_2}} \}_{k\ge 1} \sim \{  e_{k}^{\ell_{ 2 }} \}_{k\ge 1}$, i.e., $q_2=2$;
if $a\ne 0$, then by  $p_2>2$ and $$\norm{ \frac{1}{(2^k )^{1/2}}\oplus _{n=2^k+1}^{2^{k+1 } }
a   }_{\ell_{p_2,q_2}} =\frac{2^{k/{p_2}}}{2 ^{k/2}}   \norm{a}_{\ell_{p_2,q_2}}\to 0 \enspace \text{as} \enspace k\to \infty, $$ we have  
\begin{align*}
\{ e^{\ell_{q_2}}_k \}_{k\ge 1 } 
\mathrel{\overset{\text{\tiny Lem. \ref{Lpqequivalp}}}{\sim}}
\left\{ \frac{1}{(2^k )^{1/2}} \left( \oplus _{n=2^k+1}^{2^{k+1 } }
(a  \oplus e^{\ell_{ 2}}_k  ) \right) \right\} _{k\ge 1 }\sim \{ e^{\ell_{ 2}}_k \}_{k\ge 1 } .
\end{align*}
Hence, $q_2=2$. This shows that $\ell_{2,q_2}\not\hookrightarrow \ell_{p_2,q_2}\oplus \ell_2$ when $p_2, q_2\ne 2$. 
\item  Assume that $ \ell_{p_2, 2}\hookrightarrow \ell_{p_2,q_2}\oplus \ell_2$.

By   assumption (C),  $b_k  \not\to  0$ in $\ell_2$ as $k \to \infty $.
Then we have
$$ n^{1/p_2} = \norm{\sum_{k=1}^n e_k^{\ell_{p_2,2} } }_{\ell_{p_2,2}} \approx  \norm{\sum_{k=1}^n T\left(e_k^{\ell_{p_2,2 } }\right) }_{\ell_{p_2,q_2}\oplus \ell_2 }\approx  
\norm{\sum_{k=1}^n a_k \oplus b_k}_{\ell_{p_2,q_2}\oplus \ell_2} . $$
Since $p_2,q_2\ne 2$, it follows from Proposition \ref{lemmaporq} that $$n^{1/p_2}\approx n^{1/p_2}+n^{1/2}$$ or 
$$n^{1/p_2}\approx n^{1/q_2} +n^{1/2} , $$
both of which  imply that    $p_2<   2$. % or $p_2= q_2<2$.
Therefore,  
we have $ \frac{1}{n^{1/p_2}}\sum_{k=i}^{i+n} b_k  \to 0$ in  the $\ell_2$-norm for any $i$ as $n\to \infty  $. 
%In the second case, passing to a subsequence if necessary, we have
%$p_2=q_2$, i.e., we have $\ell_{p_2,2}\hookrightarrow \ell_{p_2}\oplus \ell_2$, which is impossible. Indeed, we have 
%\begin{eqnarray*}
%\{e^{\ell_2}_k\}_{k\ge 1}
%&\mathrel{\overset{\text{\tiny\cite[Lemma 2.1]{Dilworth}}}{\sim}}
%&\left \{ \frac{1}{(2^k)^{1/ p_2 }}  \sum_{n=2^k+1}^{2^{k+1 }} e^{\ell_{p_2,2}}_n  \right  \}_{k\ge 1}\\
%&\sim
%&\left\{    \frac{1}{(2^k ) ^{1/p_2}}\sum_{n=2^k +1}^{2^{k+1 }}  a_n \oplus  \frac{1}{(2^k ) ^{1/p_2}} \sum_{n=2^k +1}^{2^{k+1 }}   b_n \right\}_{k\ge 1}\\
%&\sim
%&\{e^{\ell_{p_2 }}_k\}_{k\ge 1} .  
%\end{eqnarray*}
%Hence, $p_2=2$, which is a contradiction
%For the first case,  
Passing to a subsequence if necessary,
we have 
\begin{eqnarray*}
\{e^{\ell_2}_k\}_{k\ge 1}
&\mathrel{\overset{\text{\tiny Lem. \ref{Lpqequivalp}}}{\sim}}
&\left \{ \frac{1}{(2^k)^{1/ p_2 }}  \sum_{n=2^k+1}^{2^{k+1 }} e^{\ell_{p_2,2}}_n  \right  \}_{k\ge 1}\\
&\sim
&\left\{    \frac{1}{(2^k ) ^{1/p_2}}\sum_{n=2^k +1}^{2^{k+1 }}  a_n \oplus  \frac{1}{(2^k ) ^{1/p_2}} \sum_{n=2^k +1}^{2^{k+1 }}   b_n \right\}_{k\ge 1}\\
&\sim
&\{e^{\ell_{q_2 }}_k\}_{k\ge 1} .  
\end{eqnarray*}
That is, 
%$\ell_{p_2,2}$ contains $\ell_{q_2}$, which is possible only when 
$q_2=2$, which is a contradiction. %Hence, $q_2>2$. 
%However, 
% $p_2<2$ implies that 
% $$e^{\ell_{2}} \sim \frac{1}{n^{1/p_2}} \sum e^{\ell_{p_2,2}} \sim   e^{\ell_{q_2}} . $$
%Contradiction.
\item  Assume that  $  \ell_{q_2, 2}\hookrightarrow \ell_{p_2,q_2}\oplus \ell_2$. 

If $q_2<2$, then passing to a subsequence if necessary,  we have 
\begin{eqnarray*}
\left\{e^{\ell_2}_k \right\}_{k\ge 1}
&\mathrel{\overset{\text{\tiny Lem. \ref{Lpqequivalp}}}{\sim}}
&\left \{ \frac{1}{(2^k)^{1/ q_2 }}  \sum_{n=2^k+1}^{2^{k+1 }} e^{\ell_{q_2,2}}_n  \right  \}_{k\ge 1}\\
&\sim
&\left\{    \frac{1}{(2^k ) ^{1/q_2}}\sum_{n=2^k +1}^{2^{k+1 }}  a_n \oplus  \frac{1}{(2^k ) ^{1/q_2}} \sum_{n=2^k +1}^{2^{k+1 }}   b_n \right\}_{k\ge 1}.  
\end{eqnarray*}
The right-hand-side either vanishes, or else is equivalent to  $\{ e_k ^{\ell_{q_2 } }  \}_{k\ge 1}$  (that is, $q_2=2$), both of which   are impossible. 

If $q_2>2$ and $b_k\not \to  0$ in $\ell_2 $ as $k\to \infty  $, then we have 
$$ n^{1/q_2} = \norm{\sum_{k=1}^n e_k^{\ell_{q_2,2} } }_{\ell_{q_2,2}} \approx  \norm{\sum_{k=1}^n T\left(e_k^{\ell_{q_2,2 } }\right) }_{\ell_{p_2,q_2}\oplus \ell_2 }\approx  
\norm{\sum_{k=1}^n a_k \oplus b_k}_{\ell_{p_2,q_2}\oplus \ell_2} ,$$
which means $n^{1/q_2} \gtrsim n^{1/2}$.
This is a contradiction.
Therefore, $b_k\to 0$ in $\ell_2$ as $k\to \infty  $.
Hence, $  \ell_{q_2, 2} \hookrightarrow \ell_{p_2,q_2} $ and this implies that $q_2=2$ (see Lemma~\ref{classical}), which contradicts our assumption.
\end{enumerate}
This  completes the proof. 
 \end{proof}

 Below, we present  the main result of this section.

 \begin{theorem}Let $0<p_1,p_2,q_1,q_2 <\infty $. Then, $C_{p_1,q_1}\hookrightarrow C_{p_2,q_2}$ if and only if 
 \begin{enumerate}
 \item either $p_1=q_1=2$,
 \item or $p_1=p_2$ and $q_1=q_2$. 
 \end{enumerate}
 
 \end{theorem}
 \begin{proof}
If $C_{p_1,q_1}\hookrightarrow C_{p_2,q_2}$, then we have $$\ell_{p_1,q_1}\hookrightarrow C_{p_1,q_1}\hookrightarrow C_{p_2,q_2} .$$
By Proposition~\ref{directsuml2} above, we have either (i) $p_1=q_1=2$ or (ii) $p_1=p_2$ and $q_1=q_2$ or (iii) $p_1=q_1=q_2$.
It suffices to show that the third case is impossible, that is, $C_{q} \not\hookrightarrow C_{p,q}$, $p\ne q\ne 2$.

Assume by the contradiction that  $C_{q}  \hookrightarrow C_{p,q}$, $p\ne q\ne 2$.   By Lemma \ref{classical}, $\ell_2\not\hookrightarrow \ell_{p,q},\ell_q$ when $p\ne q\not=2$. 
Moreover, $c_0\not\hookrightarrow \ell_{p,q},\ell_q$ (see e.g., 
\cite[Proposition A.3]{HSY}; see \cite[Theorem 5.9.6]{DPS2023} or \cite[Theorems 2.4.12 and 5.1.14]{Meyer} if $\ell_{p,q}$ is a Banach space). 
By  Lemma \ref{arazylemmaforisomorphism},  $\{ e_{ii}^{C_q} \}_{i\ge 1}$ in $C_{q}$ is equivalent to $\{b\otimes  e_{ii}^{C_{p,q}} \}_{i\ge 1}$ in $C_{p,q}$. 
Therefore, for any $N\ge 1$, we have  
$$N^{1/q} =\norm{\sum_{i=1}^N e_{ii}^{C_q}}_{C_q} \approx \norm{\sum_{i=1}^N b \otimes  e_{ii}^{C_{p,q}}}_{C_{p,q}} \approx _{b,p,q }N^{1/p } ,$$ 
 which implies that $p=q$.
This completes the proof.  
\end{proof}

\section{Embeddings between weighted Lorentz ideals}
In this section, we study the isomorphic embeddings between operator ideals generated by weighted Lorentz sequence spaces $d(w,p)$, $0<p<\infty$. 
In particular, we show that $C_p\not \hookrightarrow C_{d(w,p)}$ when  $0<p<\infty$ and $w\subset c_0$, which is in contrast with the classical result that $\ell_p\hookrightarrow  d(w,p) $\cite{LT1,Dilworth,ACL, Popa}.

%\subsection{Isomorphic embedding of $C
%_{d(w_1,p_1)}$ into $C_{d(w_2,p_2)}$}

\begin{lemma}\label{dcembedding}
         Let $0<p_1, p_2<\infty$.
    Let $d(w_1,p_1)$ and $d(w_2,p_2)$ be the weighted Lorentz sequence space where $w_1=\{w_n^1\}_{n\geq 1}$ and $w_2=\{w_n^2\}_{n\ge 1}$ are nonincreasing sequences.
    Let $d(w_1,p_1)\hookrightarrow C_{d(w_2,p_2)}$.
	Then, $d(w_1,p_1)\hookrightarrow d(w_2,p_2)\oplus \ell_2$, and we have
		\begin{enumerate}
			\item either $d(w_1,p_1)\approx\ell_2$,
			\item or $p_1=p_2$,             
		\end{enumerate}
\end{lemma}

	\begin{proof}

		By   \cite[Theorem 5.1]{HSZ25} (see also \cite{HSY}), we may assume that $T$ is an isomorphism from $d(w_1,p_1)$ into $ d(w_2,p_2)\oplus \ell_2$.
		Let $$ \left\{ T(e_k^{d(w_1,p_1)}) \right\}_{k\ge 1} \sim  \{ a_k\oplus
		b_k \}_{k\ge 1}\subset  d(w_2,p_2)\oplus \ell_2 .$$
		Without loss of generality, we may assume that the images $T(e_k^{d(w_1,p_1)})$ are disjoint~\cite[Lemma 2.1]{KR07}.
		
		\textbf{(1)}: If $b_k \rightarrow 0$ in $\ell_2$ as $k\rightarrow \infty$,  then passing to a subsequence if necessary,
        we have 
        $ \left\{ T(e_k^{d(w_1,p_1)}) \right\}_{k\ge 1} \sim  \{ a_k  \}_{k\ge 1}$, that is, 
        $d(w_1,p_1)\hookrightarrow d(w_2,p_2)$. 
        By Remark~\ref{equivapp}, we have $p_1=p_2$.

		\textbf{(2)}: If $b_k \not\rightarrow 0$ in $\ell_2$ as $k\rightarrow \infty$ and $a_k\rightarrow 0$ in $d(w_2,p_2)$ as $k\rightarrow \infty$, then passing to a subsequence if necessary, we have
		$$
		\left\{T(e_k^{d(w_1,p_1)})\right\}_{k\ge 1}\sim \left\{ e_k^{\ell_2}\right\}_{k\ge 1}.
		$$
		In this case, $d(w_1,p_1)\approx\ell_2$.
		
		\textbf{(3)}: If $b_k \not\rightarrow 0$ in $\ell_2$ as $k\rightarrow \infty$ and $a_k\not \rightarrow 0$ in $d(w_2,p_2)$ as $k\rightarrow \infty$,
        then passing to a subsequence if necessary, we have
		\begin{eqnarray*}
&\quad &		\left\{ e_k^{\ell_{p_1}}\right\}_{k\ge 1}\\
		&\mathrel{\overset{\text{\tiny Lemma~\ref{Lorequivalp}}}{\sim}}
		&\left\{ \frac{\sum_{n=2^k +1}^{2^{k+1}} T(e_n ^{d(w_1,p_1)})}{\norm{\sum_{n=2^k +1}^{2^{k+1}} T(e_n ^{d(w_1,p_1)})}_{d(w_2,p_2)}}  \right\}_{k\ge 1}\\
        &\sim
        &\left\{    \frac{\sum_{n=2^k +1}^{2^{k+1 }}  a_n }{\norm{\sum_{n=2^k +1}^{2^{k+1}} T(e_n ^{d(w_1,p_1)})}_{d(w_2,p_2)}} \oplus  \frac{\sum_{n=2^k +1}^{2^{k+1 }}  b_n }{\norm{\sum_{n=2^k +1}^{2^{k+1}} T(e_n ^{d(w_1,p_1)})}_{d(w_2,p_2)}}  \right\}_{k\ge 1}.
		\end{eqnarray*}
Note that 
$   \left\{  \frac{\sum_{n=2^k +1}^{2^{k+1 }}  a_n }{\norm{\sum_{n=2^k +1}^{2^{k+1}} T(e_n ^{d(w_1,p_1)})}_{d(w_2,p_2)}}     \right\}_{k\ge 1}$ vanishes or 
is equivalent to $\{e_k^{\ell_{p_2}}\}_{k\ge 1}$ (see Lemma~\ref{Lorequivalp}), 
 $\left\{ \frac{\sum_{n=2^k +1}^{2^{k+1 }}  b_n }{\norm{\sum_{n=2^k +1}^{2^{k+1}} T(e_n ^{d(w_1,p_1)})}_{d(w_2,p_2)}} \right\}_{k\ge 1}$      vanishes or is equivalent to  $\{e_k^{\ell_{2}}\}_{k\ge 1}$.
This shows that  $p_1=p_2$ or $2$.
		On the other hand,  we have 
        \begin{align*}
            1=\frac{\norm{\sum_{i=1}^n e_i^{d(w_1,p_1)}}_{d(w_1,p_1)}}{\left(\sum_{i=1}^n w_i^1\right)^{1/p_1}}
            \approx
            \frac{\norm{\sum_{i=1}^n a_i}_{d(w_2,p_2)}}{\left(\sum_{i=1}^n w_i^1\right)^{1/p_1}}
            \oplus
            \frac{\norm{\sum_{i=1}^n e_i^{\ell_2}}_{\ell_2}}{\left(\sum_{i=1}^n w_i^1\right)^{1/p_1}}
        \end{align*}
	If 
    \begin{itemize}
        \item $p_1>2$,
        \item or $p_1=2$ and $\inf_{1\le n<\infty} w_n^1=0$,
    \end{itemize}
then
    \begin{align*}
    1
    &\approx
    \lim_{n\rightarrow\infty}\frac{\norm{\sum_{i=1}^n a_i}_{d(w_2,p_2)}}{\left(\sum_{i=1}^n w_i^1\right)^{1/p_1}}
            \oplus
            \lim_{n\rightarrow\infty}\frac{\norm{\sum_{i=1}^n e_i^{\ell_2}}_{\ell_2}}{\left(\sum_{i=1}^n w_i^1\right)^{1/p_1}}\\
            &\ge
            \lim_{n\rightarrow\infty}\frac{\norm{\sum_{i=1}^n e_i^{\ell_2}}_{\ell_2}}{\left(\sum_{i=1}^n w_i^1\right)^{1/p_1}}
            =\lim_{n\rightarrow\infty}\frac{n^{1/2}}{\left(\sum_{i=1}^n w_i^1\right)^{1/p_1}}=\infty.
    \end{align*}
    This implies that $0<p_1\le2$ and when $p_1=2$, $d(w_1,p_1)$ has to be isomorphic to $\ell_2$.
    This means either $d(w_1,p_1)\approx\ell_2$ or $0<p_1=p_2<2$.
	\end{proof}

  It is well known that $\ell_p\hookrightarrow  d(w,p) $ for $0<p<\infty$ (see e.g. \cite{LT1,Dilworth,ACL, Popa}).
  The following proposition shows that this result can not extend to the noncommutative setting, i.e., if $p\not=2$ and $d(w,p)\not\approx \ell_p$ then $C_p \not\hookrightarrow C_{d(w,p)}$.

\begin{theorem}\label{ccequiva}
    Let $0<p,q<\infty$ and $d(w,p)$ be the weighted Lorentz sequence space where $w=\{w_n\}_{n\geq 1}$ is nonincreasing sequence.
        Then, $C_{q}\hookrightarrow C_{d(w,p)}$ if and only if  
             \begin{enumerate}
 \item either $q=2$,
 \item or $p=q\not=2$ and $d(w,p)\approx\ell_p$. 
 \end{enumerate}
\end{theorem}
\begin{proof}
It suffices to prove the necessity. 
It follows from $C_{q}\hookrightarrow C_{d(w,p)}$ that $\ell_q\hookrightarrow C_{d(w,p)}$.
By Lemma~\ref{dcembedding}, we have either $\ell_q\approx \ell_2$ or $p=q$.
 Assume by contradiction that   
 $p=q \ne 2$ (i.e., 
 $C_p\hookrightarrow C_{d(w,p)}$)  and $d(w,p)\not\approx\ell_p$ (in particular,  $\lim_{n \rightarrow\infty} w_{n}= 0$).
     It is well-known that
     $c_0\not\hookrightarrow d(w,p),\ell_p$\cite[Proposition A.3]{HSY}.
    By Lemma~\ref{classical}, $\ell_2\not\hookrightarrow\ell_p$ when $p\not=2$.
   Note  $\ell_2\not\hookrightarrow d(w,p)$ when $p\not=2$.
   Indeed, 
    if $\ell_2\hookrightarrow d(w,p)$, then by \cite[Theorem 1]{ACL} and \cite[Theorem 4.2]{Popa}, we have $\ell_p \hookrightarrow\ell_2$, which is impossible when  $p\not=2$ (see 
   Lemma~\ref{classical}).

    By Lemma~\ref{arazylemmaforisomorphism}, there exists $b\in C_{d(w,p)}$ such that $\{e_{ii}^{C_p}\}_{i\ge 1}$ in $C_p$ is equivalent to $\{b\otimes e_{ii}^{C_{d(w,p)}}\}_{i\ge 1}$ in $C_{d(w,p)}$.
    Therefore, for any $N\ge 1$, we have 
    \begin{align*}
        N^{1/p}=\norm{\sum_{i=1}^N e_{ii}^{C_p}}_{C_p}
        \approx 
        \norm{\sum_{i=1}^N b\otimes e_{ii}^{C_{d(w,p)}}}_{C_{d(w,p)}}
        %=
        %\left(\sum_{n=1}^{\infty} b^*(n)(s_{Ni}-s_{N(i-1)}) \right)^{1/p}
        \approx
        \left(\int_0^{\infty} \left(B(\frac{t}{N})\right)^p W(t)\,dt\right)^{1/p}
        %\stackrel{m=t/N}{\approx}
        %\left(N\int_0^{\infty} b^*(m)^p %w(Nm)\,dm\right)^{1/p}.
    \end{align*}
    where for $t\in (n-1,n]$, $n=1,2,\cdots$,  
    $$
    B(t)=\mu(n-1;b), \enspace W(t)=w_n.
    $$
     Let $m=t/N$. 
     Then, we have
     $$
     1\approx \int_0^{\infty} B^*(m)^p W(Nm)\,dm.
     $$
     Since $B$ is bounded and $\lim_{N \rightarrow\infty} w_{Nm}= 0$, it follows from the Dominated Convergence Theorem that  $\int_0^{\infty} B(m)^p W(Nm)\,dm\rightarrow 0$ as $N\rightarrow \infty$, which is a contradiction.
\end{proof}

Below, we establish  equivalent conditions for $C_{d(w_1,p_1)}\hookrightarrow C_{d(w_2,p_2)}$ when $d(w_2,p_2)$ 
has only two nonequivalent symmetric basic sequences.

\begin{theorem}\label{weighted}
     Let $0<p_1, p_2<\infty$.
    Let $d(w_1,p_1)$ and $d(w_2,p_2)$ be two weighted Lorentz sequence spaces,  where $w_1=\{w_n^1\}_{n\geq 1}$ and $w_2=\{w_n^2\}_{n\ge 1}$ are nonincreasing sequences.
    Let $s_n:=\sum_{i=1}^n w_i^2$.
    If $\sup_{1\le n
    ,k<+\infty} s_{nk}/s_ns_k<+\infty$, then $C_{d(w_1,p_1)}\hookrightarrow C_{d(w_2,p_2)}$ if and only if 
     \begin{enumerate}
 \item either $d(w_1,p_1)\approx\ell_2$,
 \item or $d(w_1,p_1)\approx d(w_2,p_2)$ (equivalently, their symmetric bases are equivalent). 
 \end{enumerate}
\end{theorem}
\begin{proof}
$(\Leftarrow)$. This follows from the uniqueness of symmetric basis of $d(w,p)$ (see \cite[Theorem 4]{ACL} and \cite[Proposition 3]{NO}). 

$(\Rightarrow)$.
It follows from $C_{d(w_1,p_1)}\hookrightarrow C_{d(w_2,p_2)}$ that $d(w_1,p_1)\hookrightarrow C_{d(w_2,p_2)}$.
    By Lemma~\ref{dcembedding}, we have either $d(w_1,p_1)\approx\ell_2$ or $p_1=p_2$.
  It suffices to consider the case when $d(w_1,p_1)\not\approx\ell_2$ and  $p_1=p_2$. 

    \textbf{(1)}: If $d(w_1,p_1)\not\approx\ell_2$ and $p_1=p_2=2$, then $$d(w_1,2)\hookrightarrow d(w_2,2)\oplus \ell_2 \approx  d(w_2,2) .$$
    If $d(w_2,2)\approx \ell_2$, then $d(w_1,2) \hookrightarrow \ell_2$, which means $d(w_1,p_1) \approx\ell_2$.
    It is a contradiction.
    If $d(w_2,2)\not\approx \ell_2$, then by Corollary~\ref{ddequiva}, we have $d(w_1,p)\approx \ell_2$ or  $d(w_1,2)\approx d(w_2,2)$.
    Hence, $d(w_1,2)\approx d(w_2,2)$.
    
    \textbf{(2)}: If $d(w_1,p_1)\not\approx\ell_2$ and $p_1=p_2\not =2$, then $c_0,\ell_2\not\hookrightarrow d(w_1,p_1), d(w_2,p_2)$ (see the proof of Proposition~\ref{ccequiva}). 
    By Lemma \ref{arazylemmaforisomorphism}, there exists $b\in C_{d(w_2,p_2)}$ such that $\{ e_{ii}^{C_{d(w_1,p_1)}} \}_{i\ge 1}$ in $C_{d(w_1,p_1)}$ is equivalent to $\{b\otimes  e_{ii}^{C_{d(w_2,p_2)}} \}_{i\ge 1}$ in $C_{d(w_2,p_2)}$, which implies that $$d(w_1,p_1)\hookrightarrow d(w_2,p_2).$$
    Let $p:=p_1=p_2$.
    If  $d(w_2,p)\approx\ell_{p}$, then $d(w_1,p)\hookrightarrow \ell_p$, which means $$d(w_1,p)\stackrel{\tiny \mbox{\cite[Theorem 2.1]{Popa}}}{\approx} \ell_p\approx d(w_2,p).$$
    If  $d(w_2,p)\not\approx\ell_{p}$, then by Corollary~\ref{ddequiva}, we have either $$
    d(w_1,p)\approx \ell_p \mbox{ or } d(w_1,p)\approx d(w_2,p).
    $$
    By Theorem~\ref{ccequiva}, the first case is impossible and we have $d(w_1,p)\approx d(w_2,p)$.
\end{proof}

We do not know whether the condition   $\sup_{1\le n
    ,k<+\infty} s_{nk}/s_ns_k<+\infty$ in Theorem \ref{weighted} is necessary.

\appendix

\section{A decomposition theorem for an isomorphism between  quasi-Banach operator ideals}\label{Arazy-qbs}

We now present a proof for the quasi-Banach version of  a result due to    Arazy \cite[Corollary 5.3]{Arazy81}. 
The   desired conclusion can be obtained by a slight  modification of  Arazy's original proof.
For completeness, we present a full proof below by applying results obtained in   \cite{HSY}.

%\begin{definition}
%Let $X$ be a Banach space. A sequence $\{X_k\}_{k=1}^\infty$ of nonzero closed subspaces of $X$ is called a \emph{Schauder decomposition} of $X$ if every $x\in X$ has a unique representation of the form $x=\sum_{k=1}^\infty x_k$ with $x_k\in X_k$ for every $k$. If $\dim(X_k)=1$ for every $k$, and $x_k\in X_k$ such that $X_k={\rm span}\{x_k\}$ for every $k$, then the sequence $\{x_k\}_{k=1}^\infty$ is called a \emph{Schauder basis} of $X$.
%\end{definition}

Let $\{K_i\}_i$ and $\{L_i\}_i$ be two finite or infinite sequences of mutually orthogonal closed subspaces of a separable infinite-dimensional Hilbert space $H$.
Denote by $p_X$ the self-adjoint projection onto a closed subspace  $X$ of $H$. 

\begin{definition}
A sequence 
$\{x_{i,j}\}_{(i,j)\in G}$   in $\mathcal B(H)$ is called a \emph{block sequence associated with $\{K_i\}_i$ and $\{L_j\}_j$} if $x_{i,j}=p_{_{K_i}}x_{i,j}p_{_{L_j}}$ for all $(i,j)\in G$, and is denoted by
\[
\{x_{i,j}\}_{(i,j)\in G}\stackrel{\mbox{{\rm\tiny b.s}}}{\boxplus}\{K_i\}_i\otimes\{L_j\}_j.
\]
\end{definition}
In particular, when $\{K_i\}_i$ \big(resp. $\{L_j\}_j$\big) contains only one element $K'$ \big(resp. $L'$\big), we may write $\{x_{i,j}\}_{(i,j)\in G}$ as $\{x_j\}_j$  \big(resp. $\{x_i\}_i$\big) and
\[
\{x_j\}_j\stackrel{\mbox{\tiny b.s}}{\boxplus}K'\otimes\{L_j\}_j\;\;\;\;\big({\rm resp.}\;\{x_i\}_i\stackrel{\mbox{\tiny b.s}}{\boxplus}\{K_i\}_i\otimes L'\big).
\]

\begin{definition}
Let $\{x_{i,j}\}_{(i,j)\in G}$ be a block sequence associated with $\{K_i\}_i$ and $\{L_j\}_j$ in $\mathcal B(H)$, and $\big\{u_i:K_i\to H'\big\}_i$ and $\big\{v_j:L_j\to H'\big\}_j$ be the two sequences of isometries. The sequence $\{x_{i,j}\}_{(i,j)\in G}$ is said to be  \emph{generated by the  system $\{u_i:K_i\to H'\}_i\otimes\{v_j:L_j\to H'\}_j$ of isometries  and an operator} if there exist an operator $a\in\mathcal B(H')$ such that
\[
u_ix_{i,j}{v_j}^\ast=a
\]
for every $(i,j)\in G$, and is denoted by
\[
\{x_{i,j}\}_{(i,j)\in G}\stackrel{\mbox{{\rm\tiny b.s}}}{\boxplus}\{u_i:K_i\to H'\}_i\otimes\{v_j:L_j\to H'\}_j\curvearrowleft a.
\] 
If one only   emphasizes the existence of such an isometries system, then the sequence $\{x_{i,j}\}_{(i,j)\in G}$ is said to be  \emph{generated by an operator}, and is denoted by
\[
\{x_{i,j}\}_{(i,j)\in G}\stackrel{\mbox{{\rm\tiny b.s}}}{\boxplus}\{K_i\}_i\otimes\{L_j\}_j\curvearrowleft a
\]
or
\[
\{x_{i,j}\}_{(i,j)\in G}\stackrel{\mbox{{\rm\tiny b.s}}}{\boxplus}\{K_i\}_i\otimes\{L_j\}_j\curvearrowleft\;
\]
if the  choice of $a$ is irrelevant. 
\end{definition}

\begin{definition}\label{definieij}
   Let $F$ be a quasi-symmetric sequence spaces, and $K$ and $L$ be two closed subspaces of $H$, we denote by \emph{$T_F(L,K)$} the closed subspace $\big[(\cdot|\eta_l)\xi_k\big]_{1\leq l\leq k<\infty}$ of $C_F$ for some orthonormal basis $\{\xi_k\}_{k=1}^\infty$ of $K$ and orthonormal basis $\{\eta_k\}_{k=1}^\infty$ of $L$. 
   %Here, the specific choice of $\{\xi_k\}_{k=1}^\infty$ and $\{\eta_k\}_{k=1}^\infty$ is immaterial, and
   We denote the operator $(\cdot|\eta_l)\xi_k$ by \emph{$e^{K,L}_{k,l}$} for every $k,l\geq1$.
\end{definition}

Recall that $T_F$ is the lower triangular part of $C_F$. 
By \cite[Remark 4.15]{HSY}, we have the following theorem:

\begin{theorem}\label{tftoce}
Suppose that $E$ and $F$ are quasi-symmetric sequence spaces with $c_0,\ell_2\not\hookrightarrow E$. If $T\in\mathcal B\big(T_F,C_E\big)$, then for any given positive number $\varepsilon$, there exist two closed subspaces $K$ and $L$ of $H$ with $T_F(L,K)\subset T_F$, and two operators $\tilde{R}_1,\tilde{R}_2\in\mathcal B\big(C_F(L,K),C_E\big)$ having the forms
\[
\tilde{R}_1(e^{K,L}_{k,l})=\left\{
\begin{array}{rcl}
x_{_{(k,l)}}+y_{_{(l,k)}}, &  & {1\leq l\leq k<\infty,}\\
0\;\;\;\;\;\;\;\;\;\;, &  & \text{otherwise,}\\
\end{array}
\right.
\]
\[
\tilde{R}_2(e^{K,L}_{k,l})=x_{2k+1,2l}+y_{2l,2k+1}\;\;\;\;\;\;\;\;{\rm for\;every}\;k,l\geq1,
\]
such that $\big\Vert T|_{T_F(L,K)}-\big(\tilde{R}_1+\tilde{R}_2\big)|_{T_F(L,K)}\big\Vert<\varepsilon$, where
\begin{itemize}
\item [(1)]
There are two sequences $\{F_{(k,l)}\}_{1\leq l\leq k<\infty}$ and $\{G_{(l,k)}\}_{1\leq l\leq k<\infty}$   of mutually orthogonal closed subspaces of $H$ such that
\[
\{x_{_{(k,l)}}\}_{1\leq l\leq k<\infty}\stackrel{\mbox{{\rm\tiny b.s}}}{\boxplus}\{F_{(k,l)}\}_{1\leq l\leq k<\infty}\otimes H\curvearrowleft\;
\]
and
\[
\{y_{_{(l,k)}}\}_{1\leq l\leq k<\infty}\stackrel{\mbox{{\rm\tiny b.s}}}{\boxplus} H\otimes\{G_{(l,k)}\}_{1\leq l\leq k<\infty}\curvearrowleft\;,
\]
\item [(2)]
There are two sequences $\{K_i\}_{i=1}^\infty$ and $\{L_i\}_{i=1}^\infty$ of   mutually orthogonal closed subspaces of $H$ such that
\begin{equation}\label{x2k+12l}
\{x_{2k+1,2l}\}_{k,l=1}^\infty\stackrel{\mbox{{\rm\tiny b.s}}}{\boxplus}\{K_{2k+1}\}_{k=1}^\infty\otimes\{L_{2l}\}_{l=1}^\infty\curvearrowleft\;
\end{equation}
and
\begin{equation}\label{y2l2k+1}
\{y_{2l,2k+1}\}_{k,l=1}^\infty\stackrel{\mbox{{\rm\tiny b.s}}}{\boxplus}\{K_{2l}\}_{l=1}^\infty\otimes\{L_{2k+1}\}_{k=1}^\infty\curvearrowleft\;.
\end{equation}
\end{itemize}
\end{theorem}

The following two figures are helpful for understanding the forms of the sequences $\{x_{_{(k,l)}}\}_{1\leq l\leq k<\infty}$, $\{y_{_{(l,k)}}\}_{1\leq l\leq k<\infty}$, $\{x_{2k+1,2l}\}_{k,l=1}^\infty$ and $\{y_{2l,2k+1}\}_{k,l=1}^\infty$:
\vspace{0.4cm}

%$\{x_{_{(k,l)}}\}_{1\leq l\leq k<\infty}$ is the horizontal bar in the figure below, and $\{y_{_{(l,k)}}\}_{1\leq l\leq k<\infty}$ is the vertical bar in the figure below.
%\noindent
%\\
% Pattern Info

\tikzset{
pattern size/.store in=\mcSize, 
pattern size = 5pt,
pattern thickness/.store in=\mcThickness, 
pattern thickness = 0.3pt,
pattern radius/.store in=\mcRadius, 
pattern radius = 1pt}
\makeatletter
\pgfutil@ifundefined{pgf@pattern@name@_r2sae6z0w}{
\pgfdeclarepatternformonly[\mcThickness,\mcSize]{_r2sae6z0w}
{\pgfqpoint{0pt}{0pt}}
{\pgfpoint{\mcSize+\mcThickness}{\mcSize+\mcThickness}}
{\pgfpoint{\mcSize}{\mcSize}}
{
\pgfsetcolor{\tikz@pattern@color}
\pgfsetlinewidth{\mcThickness}
\pgfpathmoveto{\pgfqpoint{0pt}{0pt}}
\pgfpathlineto{\pgfpoint{\mcSize+\mcThickness}{\mcSize+\mcThickness}}
\pgfusepath{stroke}
}}
\makeatother

% Pattern Info
 
\tikzset{
pattern size/.store in=\mcSize, 
pattern size = 5pt,
pattern thickness/.store in=\mcThickness, 
pattern thickness = 0.3pt,
pattern radius/.store in=\mcRadius, 
pattern radius = 1pt}
\makeatletter
\pgfutil@ifundefined{pgf@pattern@name@_c3llu2tdm}{
\pgfdeclarepatternformonly[\mcThickness,\mcSize]{_c3llu2tdm}
{\pgfqpoint{0pt}{0pt}}
{\pgfpoint{\mcSize+\mcThickness}{\mcSize+\mcThickness}}
{\pgfpoint{\mcSize}{\mcSize}}
{
\pgfsetcolor{\tikz@pattern@color}
\pgfsetlinewidth{\mcThickness}
\pgfpathmoveto{\pgfqpoint{0pt}{0pt}}
\pgfpathlineto{\pgfpoint{\mcSize+\mcThickness}{\mcSize+\mcThickness}}
\pgfusepath{stroke}
}}
\makeatother

% Pattern Info
 
\tikzset{
pattern size/.store in=\mcSize, 
pattern size = 5pt,
pattern thickness/.store in=\mcThickness, 
pattern thickness = 0.3pt,
pattern radius/.store in=\mcRadius, 
pattern radius = 1pt}
\makeatletter
\pgfutil@ifundefined{pgf@pattern@name@_xgzcfukxn}{
\pgfdeclarepatternformonly[\mcThickness,\mcSize]{_xgzcfukxn}
{\pgfqpoint{0pt}{0pt}}
{\pgfpoint{\mcSize+\mcThickness}{\mcSize+\mcThickness}}
{\pgfpoint{\mcSize}{\mcSize}}
{
\pgfsetcolor{\tikz@pattern@color}
\pgfsetlinewidth{\mcThickness}
\pgfpathmoveto{\pgfqpoint{0pt}{0pt}}
\pgfpathlineto{\pgfpoint{\mcSize+\mcThickness}{\mcSize+\mcThickness}}
\pgfusepath{stroke}
}}
\makeatother

% Pattern Info
 
\tikzset{
pattern size/.store in=\mcSize, 
pattern size = 5pt,
pattern thickness/.store in=\mcThickness, 
pattern thickness = 0.3pt,
pattern radius/.store in=\mcRadius, 
pattern radius = 1pt}
\makeatletter
\pgfutil@ifundefined{pgf@pattern@name@_pegnvansc}{
\pgfdeclarepatternformonly[\mcThickness,\mcSize]{_pegnvansc}
{\pgfqpoint{0pt}{0pt}}
{\pgfpoint{\mcSize+\mcThickness}{\mcSize+\mcThickness}}
{\pgfpoint{\mcSize}{\mcSize}}
{
\pgfsetcolor{\tikz@pattern@color}
\pgfsetlinewidth{\mcThickness}
\pgfpathmoveto{\pgfqpoint{0pt}{0pt}}
\pgfpathlineto{\pgfpoint{\mcSize+\mcThickness}{\mcSize+\mcThickness}}
\pgfusepath{stroke}
}}
\makeatother

% Pattern Info
 
\tikzset{
pattern size/.store in=\mcSize, 
pattern size = 5pt,
pattern thickness/.store in=\mcThickness, 
pattern thickness = 0.3pt,
pattern radius/.store in=\mcRadius, 
pattern radius = 1pt}
\makeatletter
\pgfutil@ifundefined{pgf@pattern@name@_6xgdbg32s}{
\pgfdeclarepatternformonly[\mcThickness,\mcSize]{_6xgdbg32s}
{\pgfqpoint{0pt}{0pt}}
{\pgfpoint{\mcSize+\mcThickness}{\mcSize+\mcThickness}}
{\pgfpoint{\mcSize}{\mcSize}}
{
\pgfsetcolor{\tikz@pattern@color}
\pgfsetlinewidth{\mcThickness}
\pgfpathmoveto{\pgfqpoint{0pt}{0pt}}
\pgfpathlineto{\pgfpoint{\mcSize+\mcThickness}{\mcSize+\mcThickness}}
\pgfusepath{stroke}
}}
\makeatother

% Pattern Info
 
\tikzset{
pattern size/.store in=\mcSize, 
pattern size = 5pt,
pattern thickness/.store in=\mcThickness, 
pattern thickness = 0.3pt,
pattern radius/.store in=\mcRadius, 
pattern radius = 1pt}
\makeatletter
\pgfutil@ifundefined{pgf@pattern@name@_ar3mleoup}{
\pgfdeclarepatternformonly[\mcThickness,\mcSize]{_ar3mleoup}
{\pgfqpoint{0pt}{0pt}}
{\pgfpoint{\mcSize+\mcThickness}{\mcSize+\mcThickness}}
{\pgfpoint{\mcSize}{\mcSize}}
{
\pgfsetcolor{\tikz@pattern@color}
\pgfsetlinewidth{\mcThickness}
\pgfpathmoveto{\pgfqpoint{0pt}{0pt}}
\pgfpathlineto{\pgfpoint{\mcSize+\mcThickness}{\mcSize+\mcThickness}}
\pgfusepath{stroke}
}}
\makeatother

% Pattern Info
 
\tikzset{
pattern size/.store in=\mcSize, 
pattern size = 5pt,
pattern thickness/.store in=\mcThickness, 
pattern thickness = 0.3pt,
pattern radius/.store in=\mcRadius, 
pattern radius = 1pt}
\makeatletter
\pgfutil@ifundefined{pgf@pattern@name@_7ntc0297c}{
\pgfdeclarepatternformonly[\mcThickness,\mcSize]{_7ntc0297c}
{\pgfqpoint{0pt}{-\mcThickness}}
{\pgfpoint{\mcSize}{\mcSize}}
{\pgfpoint{\mcSize}{\mcSize}}
{
\pgfsetcolor{\tikz@pattern@color}
\pgfsetlinewidth{\mcThickness}
\pgfpathmoveto{\pgfqpoint{0pt}{\mcSize}}
\pgfpathlineto{\pgfpoint{\mcSize+\mcThickness}{-\mcThickness}}
\pgfusepath{stroke}
}}
\makeatother

% Pattern Info
 
\tikzset{
pattern size/.store in=\mcSize, 
pattern size = 5pt,
pattern thickness/.store in=\mcThickness, 
pattern thickness = 0.3pt,
pattern radius/.store in=\mcRadius, 
pattern radius = 1pt}
\makeatletter
\pgfutil@ifundefined{pgf@pattern@name@_xgxpvldlz}{
\pgfdeclarepatternformonly[\mcThickness,\mcSize]{_xgxpvldlz}
{\pgfqpoint{0pt}{-\mcThickness}}
{\pgfpoint{\mcSize}{\mcSize}}
{\pgfpoint{\mcSize}{\mcSize}}
{
\pgfsetcolor{\tikz@pattern@color}
\pgfsetlinewidth{\mcThickness}
\pgfpathmoveto{\pgfqpoint{0pt}{\mcSize}}
\pgfpathlineto{\pgfpoint{\mcSize+\mcThickness}{-\mcThickness}}
\pgfusepath{stroke}
}}
\makeatother

% Pattern Info
 
\tikzset{
pattern size/.store in=\mcSize, 
pattern size = 5pt,
pattern thickness/.store in=\mcThickness, 
pattern thickness = 0.3pt,
pattern radius/.store in=\mcRadius, 
pattern radius = 1pt}
\makeatletter
\pgfutil@ifundefined{pgf@pattern@name@_0oj3xqw0w}{
\pgfdeclarepatternformonly[\mcThickness,\mcSize]{_0oj3xqw0w}
{\pgfqpoint{0pt}{-\mcThickness}}
{\pgfpoint{\mcSize}{\mcSize}}
{\pgfpoint{\mcSize}{\mcSize}}
{
\pgfsetcolor{\tikz@pattern@color}
\pgfsetlinewidth{\mcThickness}
\pgfpathmoveto{\pgfqpoint{0pt}{\mcSize}}
\pgfpathlineto{\pgfpoint{\mcSize+\mcThickness}{-\mcThickness}}
\pgfusepath{stroke}
}}
\makeatother

% Pattern Info
 
\tikzset{
pattern size/.store in=\mcSize, 
pattern size = 5pt,
pattern thickness/.store in=\mcThickness, 
pattern thickness = 0.3pt,
pattern radius/.store in=\mcRadius, 
pattern radius = 1pt}
\makeatletter
\pgfutil@ifundefined{pgf@pattern@name@_nr9a77uzn}{
\pgfdeclarepatternformonly[\mcThickness,\mcSize]{_nr9a77uzn}
{\pgfqpoint{0pt}{-\mcThickness}}
{\pgfpoint{\mcSize}{\mcSize}}
{\pgfpoint{\mcSize}{\mcSize}}
{
\pgfsetcolor{\tikz@pattern@color}
\pgfsetlinewidth{\mcThickness}
\pgfpathmoveto{\pgfqpoint{0pt}{\mcSize}}
\pgfpathlineto{\pgfpoint{\mcSize+\mcThickness}{-\mcThickness}}
\pgfusepath{stroke}
}}
\makeatother

% Pattern Info
 
\tikzset{
pattern size/.store in=\mcSize, 
pattern size = 5pt,
pattern thickness/.store in=\mcThickness, 
pattern thickness = 0.3pt,
pattern radius/.store in=\mcRadius, 
pattern radius = 1pt}
\makeatletter
\pgfutil@ifundefined{pgf@pattern@name@_24e9ikrdj}{
\pgfdeclarepatternformonly[\mcThickness,\mcSize]{_24e9ikrdj}
{\pgfqpoint{0pt}{-\mcThickness}}
{\pgfpoint{\mcSize}{\mcSize}}
{\pgfpoint{\mcSize}{\mcSize}}
{
\pgfsetcolor{\tikz@pattern@color}
\pgfsetlinewidth{\mcThickness}
\pgfpathmoveto{\pgfqpoint{0pt}{\mcSize}}
\pgfpathlineto{\pgfpoint{\mcSize+\mcThickness}{-\mcThickness}}
\pgfusepath{stroke}
}}
\makeatother

% Pattern Info
 
\tikzset{
pattern size/.store in=\mcSize, 
pattern size = 5pt,
pattern thickness/.store in=\mcThickness, 
pattern thickness = 0.3pt,
pattern radius/.store in=\mcRadius, 
pattern radius = 1pt}
\makeatletter
\pgfutil@ifundefined{pgf@pattern@name@_5x8yqsjy0}{
\pgfdeclarepatternformonly[\mcThickness,\mcSize]{_5x8yqsjy0}
{\pgfqpoint{0pt}{-\mcThickness}}
{\pgfpoint{\mcSize}{\mcSize}}
{\pgfpoint{\mcSize}{\mcSize}}
{
\pgfsetcolor{\tikz@pattern@color}
\pgfsetlinewidth{\mcThickness}
\pgfpathmoveto{\pgfqpoint{0pt}{\mcSize}}
\pgfpathlineto{\pgfpoint{\mcSize+\mcThickness}{-\mcThickness}}
\pgfusepath{stroke}
}}
\makeatother
\tikzset{every picture/.style={line width=0.75pt}} %set default line width to 0.75pt        

\scalebox{0.8}{ 
\begin{tikzpicture}[x=0.75pt,y=0.75pt,yscale=-1,xscale=1]
%uncomment if require: \path (0,764); %set diagram left start at 0, and has height of 764

%Straight Lines [id:da7723356400252037] 
\draw [line width=1.5]    (220.44,130.24) -- (569.67,130.24) ;
%Straight Lines [id:da2966426746792654] 
\draw [line width=1.5]    (220.44,130.24) -- (220.67,481.25) ;
%Shape: Brace [id:dp08810412131114276] 
\draw   (220.11,150.86) .. controls (216.12,150.85) and (214.12,152.84) .. (214.11,156.83) -- (214.11,156.83) .. controls (214.1,162.52) and (212.1,165.37) .. (208.11,165.36) .. controls (212.1,165.37) and (214.08,168.22) .. (214.07,173.92)(214.08,171.35) -- (214.07,173.92) .. controls (214.06,177.91) and (216.05,179.9) .. (220.04,179.91) ;
%Shape: Brace [id:dp6500567071345742] 
\draw   (220.41,200.61) .. controls (216.38,200.56) and (214.35,202.55) .. (214.3,206.57) -- (214.3,206.57) .. controls (214.23,212.32) and (212.18,215.16) .. (208.16,215.11) .. controls (212.18,215.16) and (214.15,218.06) .. (214.08,223.8)(214.12,221.22) -- (214.08,223.8) .. controls (214.03,227.82) and (216.02,229.86) .. (220.04,229.91) ;
%Shape: Brace [id:dp22371577658633934] 
\draw   (220.41,250.61) .. controls (216.38,250.6) and (214.37,252.61) .. (214.36,256.63) -- (214.36,256.63) .. controls (214.35,262.38) and (212.33,265.25) .. (208.3,265.24) .. controls (212.33,265.25) and (214.33,268.12) .. (214.32,273.87)(214.32,271.28) -- (214.32,273.87) .. controls (214.31,277.89) and (216.31,279.91) .. (220.33,279.92) ;
%Shape: Brace [id:dp265801352721051] 
\draw   (220.41,300.61) .. controls (216.38,300.61) and (214.37,302.62) .. (214.37,306.65) -- (214.37,306.65) .. controls (214.37,312.4) and (212.36,315.28) .. (208.33,315.28) .. controls (212.36,315.28) and (214.37,318.16) .. (214.37,323.91)(214.37,321.32) -- (214.37,323.91) .. controls (214.37,327.94) and (216.38,329.95) .. (220.41,329.95) ;
%Shape: Brace [id:dp025600710386486702] 
\draw   (220.78,350.86) .. controls (216.7,350.8) and (214.63,352.81) .. (214.57,356.89) -- (214.57,356.89) .. controls (214.48,362.72) and (212.4,365.6) .. (208.32,365.54) .. controls (212.4,365.6) and (214.4,368.54) .. (214.31,374.37)(214.35,371.75) -- (214.31,374.37) .. controls (214.25,378.45) and (216.26,380.52) .. (220.33,380.58) ;
%Shape: Brace [id:dp09856799169973551] 
\draw   (220.11,400.86) .. controls (216.12,400.85) and (214.12,402.84) .. (214.11,406.83) -- (214.11,406.83) .. controls (214.1,412.52) and (212.1,415.37) .. (208.11,415.36) .. controls (212.1,415.37) and (214.08,418.22) .. (214.07,423.92)(214.08,421.35) -- (214.07,423.92) .. controls (214.06,427.91) and (216.05,429.9) .. (220.04,429.91) ;
%Shape: Brace [id:dp39901725568222335] 
\draw   (270.02,130.4) .. controls (270.11,126.39) and (268.16,124.34) .. (264.15,124.24) -- (264.15,124.24) .. controls (258.42,124.11) and (255.6,122.03) .. (255.69,118.02) .. controls (255.6,122.03) and (252.68,123.97) .. (246.95,123.83)(249.53,123.89) -- (246.95,123.83) .. controls (242.94,123.73) and (240.89,125.69) .. (240.79,129.7) ;
%Shape: Brace [id:dp11182617106058079] 
\draw   (320.02,130.4) .. controls (320.04,126.44) and (318.07,124.45) .. (314.1,124.42) -- (314.1,124.42) .. controls (308.43,124.39) and (305.61,122.4) .. (305.64,118.43) .. controls (305.61,122.4) and (302.77,124.36) .. (297.11,124.33)(299.66,124.35) -- (297.11,124.33) .. controls (293.14,124.31) and (291.15,126.28) .. (291.13,130.25) ;
%Shape: Brace [id:dp3397882726156697] 
\draw   (370.02,130.4) .. controls (370.09,126.44) and (368.14,124.43) .. (364.17,124.36) -- (364.17,124.36) .. controls (358.5,124.26) and (355.7,122.23) .. (355.77,118.26) .. controls (355.7,122.23) and (352.84,124.16) .. (347.17,124.07)(349.72,124.11) -- (347.17,124.07) .. controls (343.21,124) and (341.2,125.95) .. (341.13,129.92) ;
%Shape: Brace [id:dp17093636261788347] 
\draw   (420.02,130.4) .. controls (420.09,126.44) and (418.14,124.43) .. (414.17,124.36) -- (414.17,124.36) .. controls (408.5,124.26) and (405.7,122.23) .. (405.77,118.26) .. controls (405.7,122.23) and (402.84,124.16) .. (397.17,124.07)(399.72,124.11) -- (397.17,124.07) .. controls (393.21,124) and (391.2,125.95) .. (391.13,129.92) ;
%Shape: Brace [id:dp5817894176499134] 
\draw   (469.8,130.4) .. controls (469.82,126.42) and (467.84,124.42) .. (463.86,124.4) -- (463.86,124.4) .. controls (458.17,124.37) and (455.34,122.37) .. (455.36,118.38) .. controls (455.34,122.37) and (452.49,124.34) .. (446.8,124.31)(449.36,124.32) -- (446.8,124.31) .. controls (442.81,124.29) and (440.81,126.27) .. (440.79,130.25) ;
%Shape: Brace [id:dp034683602023634075] 
\draw   (520.18,130.41) .. controls (520.25,126.33) and (518.24,124.26) .. (514.16,124.19) -- (514.16,124.19) .. controls (508.33,124.1) and (505.45,122.01) .. (505.52,117.93) .. controls (505.45,122.01) and (502.51,124) .. (496.68,123.9)(499.3,123.95) -- (496.68,123.9) .. controls (492.6,123.83) and (490.53,125.84) .. (490.46,129.92) ;
%Straight Lines [id:da20365509062432663] 
\draw [pattern=_r2sae6z0w,pattern size=6pt,pattern thickness=0.75pt,pattern radius=0pt, pattern color={rgb, 255:red, 0; green, 0; blue, 0}]   (570.44,180.24) -- (220.67,180.25) -- (220.67,150.92) -- (570.44,150.91) ;
%Straight Lines [id:da5963736382159682] 
\draw [pattern=_c3llu2tdm,pattern size=6pt,pattern thickness=0.75pt,pattern radius=0pt, pattern color={rgb, 255:red, 0; green, 0; blue, 0}]   (570.67,230.24) -- (220.67,230.25) -- (220.67,200.92) -- (570.67,200.91) ;
%Straight Lines [id:da6239287181023034] 
\draw [pattern=_xgzcfukxn,pattern size=6pt,pattern thickness=0.75pt,pattern radius=0pt, pattern color={rgb, 255:red, 0; green, 0; blue, 0}]   (570.67,280.24) -- (220.67,280.25) -- (220.67,250.92) -- (570.67,250.91) ;
%Straight Lines [id:da058937341640299956] 
\draw [pattern=_pegnvansc,pattern size=6pt,pattern thickness=0.75pt,pattern radius=0pt, pattern color={rgb, 255:red, 0; green, 0; blue, 0}]   (570.67,330.24) -- (220.67,330.25) -- (220.67,300.92) -- (570.67,300.91) ;
%Straight Lines [id:da5723099704828642] 
\draw [pattern=_6xgdbg32s,pattern size=6pt,pattern thickness=0.75pt,pattern radius=0pt, pattern color={rgb, 255:red, 0; green, 0; blue, 0}]   (570.67,380.24) -- (220.67,380.25) -- (220.67,350.92) -- (570.67,350.91) ;
%Straight Lines [id:da21582679997211862] 
\draw [pattern=_ar3mleoup,pattern size=6pt,pattern thickness=0.75pt,pattern radius=0pt, pattern color={rgb, 255:red, 0; green, 0; blue, 0}]   (570.67,430.24) -- (220.67,430.25) -- (220.67,400.92) -- (569.67,400.91) ;
%Straight Lines [id:da21719213370521284] 
\draw [pattern=_7ntc0297c,pattern size=6pt,pattern thickness=0.75pt,pattern radius=0pt, pattern color={rgb, 255:red, 155; green, 155; blue, 155}]   (290.67,480.58) -- (290.66,130.58) -- (319.99,130.58) -- (320,480.58) ;
%Straight Lines [id:da2845830142064236] 
\draw [pattern=_xgxpvldlz,pattern size=6pt,pattern thickness=0.75pt,pattern radius=0pt, pattern color={rgb, 255:red, 155; green, 155; blue, 155}]   (240.67,480.58) -- (240.66,130.58) -- (269.99,130.58) -- (270,480.58) ;
%Straight Lines [id:da7455635478461115] 
\draw [pattern=_0oj3xqw0w,pattern size=6pt,pattern thickness=0.75pt,pattern radius=0pt, pattern color={rgb, 255:red, 155; green, 155; blue, 155}]   (390.67,480.58) -- (390.66,130.58) -- (419.99,130.58) -- (420,480.58) ;
%Straight Lines [id:da7327190505756134] 
\draw [pattern=_nr9a77uzn,pattern size=6pt,pattern thickness=0.75pt,pattern radius=0pt, pattern color={rgb, 255:red, 155; green, 155; blue, 155}]   (340.67,480.58) -- (340.66,130.58) -- (369.99,130.58) -- (370,480.58) ;
%Straight Lines [id:da5454769564427788] 
\draw [pattern=_24e9ikrdj,pattern size=6pt,pattern thickness=0.75pt,pattern radius=0pt, pattern color={rgb, 255:red, 155; green, 155; blue, 155}]   (490.67,480.58) -- (490.66,130.58) -- (519.99,130.58) -- (520,480.58) ;
%Straight Lines [id:da9078799653345291] 
\draw [pattern=_5x8yqsjy0,pattern size=6pt,pattern thickness=0.75pt,pattern radius=0pt, pattern color={rgb, 255:red, 155; green, 155; blue, 155}]   (440.67,480.58) -- (440.66,130.58) -- (469.99,130.58) -- (470,480.58) ;
%Straight Lines [id:da883560593489499] 
\draw    (153.33,198.24) -- (231.47,167.64) ;
\draw [shift={(233.33,166.91)}, rotate = 158.61] [color={rgb, 255:red, 0; green, 0; blue, 0 }  ][line width=0.75]    (10.93,-3.29) .. controls (6.95,-1.4) and (3.31,-0.3) .. (0,0) .. controls (3.31,0.3) and (6.95,1.4) .. (10.93,3.29)   ;
%Straight Lines [id:da3455968820942169] 
\draw    (152.67,248.58) -- (230.8,217.97) ;
\draw [shift={(232.67,217.24)}, rotate = 158.61] [color={rgb, 255:red, 0; green, 0; blue, 0 }  ][line width=0.75]    (10.93,-3.29) .. controls (6.95,-1.4) and (3.31,-0.3) .. (0,0) .. controls (3.31,0.3) and (6.95,1.4) .. (10.93,3.29)   ;
%Straight Lines [id:da25815227068032487] 
\draw    (153,298.58) -- (231.14,267.97) ;
\draw [shift={(233,267.24)}, rotate = 158.61] [color={rgb, 255:red, 0; green, 0; blue, 0 }  ][line width=0.75]    (10.93,-3.29) .. controls (6.95,-1.4) and (3.31,-0.3) .. (0,0) .. controls (3.31,0.3) and (6.95,1.4) .. (10.93,3.29)   ;
%Straight Lines [id:da3589162256540678] 
\draw    (153.33,348.91) -- (231.47,318.31) ;
\draw [shift={(233.33,317.58)}, rotate = 158.61] [color={rgb, 255:red, 0; green, 0; blue, 0 }  ][line width=0.75]    (10.93,-3.29) .. controls (6.95,-1.4) and (3.31,-0.3) .. (0,0) .. controls (3.31,0.3) and (6.95,1.4) .. (10.93,3.29)   ;
%Straight Lines [id:da48275645356288577] 
\draw    (152.33,398.91) -- (230.47,368.31) ;
\draw [shift={(232.33,367.58)}, rotate = 158.61] [color={rgb, 255:red, 0; green, 0; blue, 0 }  ][line width=0.75]    (10.93,-3.29) .. controls (6.95,-1.4) and (3.31,-0.3) .. (0,0) .. controls (3.31,0.3) and (6.95,1.4) .. (10.93,3.29)   ;
%Straight Lines [id:da5951495169345221] 
\draw    (151.67,449.24) -- (229.8,418.64) ;
\draw [shift={(231.67,417.91)}, rotate = 158.61] [color={rgb, 255:red, 0; green, 0; blue, 0 }  ][line width=0.75]    (10.93,-3.29) .. controls (6.95,-1.4) and (3.31,-0.3) .. (0,0) .. controls (3.31,0.3) and (6.95,1.4) .. (10.93,3.29)   ;
%Straight Lines [id:da2741107063423531] 
\draw    (297.72,78.84) -- (257.11,140) ;
\draw [shift={(256,141.67)}, rotate = 303.58] [color={rgb, 255:red, 0; green, 0; blue, 0 }  ][line width=0.75]    (10.93,-3.29) .. controls (6.95,-1.4) and (3.31,-0.3) .. (0,0) .. controls (3.31,0.3) and (6.95,1.4) .. (10.93,3.29)   ;
%Straight Lines [id:da5481991998040879] 
\draw    (348.38,78.84) -- (307.77,140) ;
\draw [shift={(306.67,141.67)}, rotate = 303.58] [color={rgb, 255:red, 0; green, 0; blue, 0 }  ][line width=0.75]    (10.93,-3.29) .. controls (6.95,-1.4) and (3.31,-0.3) .. (0,0) .. controls (3.31,0.3) and (6.95,1.4) .. (10.93,3.29)   ;
%Straight Lines [id:da5908124398771761] 
\draw    (398.38,79.51) -- (357.77,140.67) ;
\draw [shift={(356.67,142.33)}, rotate = 303.58] [color={rgb, 255:red, 0; green, 0; blue, 0 }  ][line width=0.75]    (10.93,-3.29) .. controls (6.95,-1.4) and (3.31,-0.3) .. (0,0) .. controls (3.31,0.3) and (6.95,1.4) .. (10.93,3.29)   ;
%Straight Lines [id:da06546307782633776] 
\draw    (449.05,79.51) -- (408.44,140.67) ;
\draw [shift={(407.33,142.33)}, rotate = 303.58] [color={rgb, 255:red, 0; green, 0; blue, 0 }  ][line width=0.75]    (10.93,-3.29) .. controls (6.95,-1.4) and (3.31,-0.3) .. (0,0) .. controls (3.31,0.3) and (6.95,1.4) .. (10.93,3.29)   ;
%Straight Lines [id:da8351857928086602] 
\draw    (497.72,80.17) -- (457.11,141.33) ;
\draw [shift={(456,143)}, rotate = 303.58] [color={rgb, 255:red, 0; green, 0; blue, 0 }  ][line width=0.75]    (10.93,-3.29) .. controls (6.95,-1.4) and (3.31,-0.3) .. (0,0) .. controls (3.31,0.3) and (6.95,1.4) .. (10.93,3.29)   ;
%Straight Lines [id:da20124733702837672] 
\draw    (548.38,80.17) -- (507.77,141.33) ;
\draw [shift={(506.67,143)}, rotate = 303.58] [color={rgb, 255:red, 0; green, 0; blue, 0 }  ][line width=0.75]    (10.93,-3.29) .. controls (6.95,-1.4) and (3.31,-0.3) .. (0,0) .. controls (3.31,0.3) and (6.95,1.4) .. (10.93,3.29)   ;

% Text Node
\draw (190.53,164) node    {$F_{( 1,1)}$};
% Text Node
\draw (189.87,213.67) node    {$F_{( 2,1)}$};
% Text Node
\draw (137.52,253.67) node    {$x_{_{( 2,1)}}$};
% Text Node
\draw (190.2,264.33) node    {$F_{( 2,2)}$};
% Text Node
\draw (137.85,303.67) node    {$x_{_{( 2,2)}}$};
% Text Node
\draw (190.53,314.67) node    {$F_{( 3,1)}$};
% Text Node
\draw (138.18,354) node    {$x_{_{( 3,1)}}$};
% Text Node
\draw (189.53,364.67) node    {$F_{( 3,2)}$};
% Text Node
\draw (137.18,404) node    {$x_{_{( 3,2)}}$};
% Text Node
\draw (188.87,415) node    {$F_{( 3,3)}$};
% Text Node
\draw (136.52,454.33) node    {$x_{_{( 3,3)}}$};
% Text Node
\draw (254.61,109) node    {$G_{( 1,1)}$};
% Text Node
\draw (138.18,203.33) node    {$x_{_{( 1,1)}}$};
% Text Node
\draw (301.85,71.67) node    {$y_{_{( 1,1)}}$};
% Text Node
\draw (305.28,109) node    {$G_{( 1,2)}$};
% Text Node
\draw (352.52,71.67) node    {$y_{_{( 1,2)}}$};
% Text Node
\draw (355.28,109.67) node    {$G_{( 2,2)}$};
% Text Node
\draw (402.52,72.33) node    {$y_{_{( 2,2)}}$};
% Text Node
\draw (403.28,108.67) node    {$G_{( 1,3)}$};
% Text Node
\draw (453.18,72.33) node    {$y_{_{( 1,3)}}$};
% Text Node
\draw (453.95,110.33) node    {$G_{( 2,3)}$};
% Text Node
\draw (487.33,63.4) node [anchor=north west][inner sep=0.75pt]    {$y_{_{( 2,3)}}$};
% Text Node
\draw (505.28,110) node    {$G_{( 3,3)}$};
% Text Node
\draw (552.52,73) node    {$y_{_{( 3,3)}}$};
% Text Node
\draw (553.93,288.58) node  [font=\large]  {$\cdots $};
% Text Node
\draw (381.26,464.91) node  [font=\large,rotate=-90]  {$\cdots $};

\end{tikzpicture}
} 

The set $\{x_{2k+1,2l}\}_{k,l=1}^\infty$ stands for  the red squares in the figure below, and $\{y_{2l,2k+1}\}_{k,l=1}^\infty$ stands for the blue squares in the figure below.

\begin{center}
\scalebox{0.8}{ 
\tikzset{every picture/.style={line width=0.75pt}} %set default line width to 0.75pt        

\begin{tikzpicture}[x=0.75pt,y=0.75pt,yscale=-1,xscale=1]
%uncomment if require: \path (0,587); %set diagram left start at 0, and has height of 587

%Straight Lines [id:da5982246747328697] 
\draw [line width=1.5]    (220.44,89.8) -- (569.67,89.8) ;
%Straight Lines [id:da17960560551585125] 
\draw [line width=1.5]    (220.44,89.8) -- (220.43,421.15) ;
%Shape: Brace [id:dp7476766124589923] 
\draw   (220.11,110.42) .. controls (216.12,110.41) and (214.12,112.4) .. (214.11,116.39) -- (214.11,116.39) .. controls (214.1,122.08) and (212.1,124.92) .. (208.11,124.91) .. controls (212.1,124.92) and (214.08,127.78) .. (214.07,133.47)(214.08,130.91) -- (214.07,133.47) .. controls (214.06,137.46) and (216.05,139.46) .. (220.04,139.47) ;
%Shape: Brace [id:dp9931359129395855] 
\draw   (220.41,160.17) .. controls (216.38,160.12) and (214.35,162.11) .. (214.3,166.13) -- (214.3,166.13) .. controls (214.23,171.87) and (212.18,174.72) .. (208.16,174.67) .. controls (212.18,174.72) and (214.15,177.61) .. (214.08,183.36)(214.12,180.78) -- (214.08,183.36) .. controls (214.03,187.38) and (216.02,189.42) .. (220.04,189.47) ;
%Shape: Brace [id:dp2613797667242559] 
\draw   (220.41,210.17) .. controls (216.38,210.16) and (214.37,212.17) .. (214.36,216.19) -- (214.36,216.19) .. controls (214.35,221.94) and (212.33,224.8) .. (208.3,224.79) .. controls (212.33,224.8) and (214.33,227.68) .. (214.32,233.43)(214.32,230.84) -- (214.32,233.43) .. controls (214.31,237.45) and (216.31,239.46) .. (220.33,239.47) ;
%Shape: Brace [id:dp20123138528497453] 
\draw   (220.41,260.17) .. controls (216.38,260.17) and (214.37,262.18) .. (214.37,266.21) -- (214.37,266.21) .. controls (214.37,271.96) and (212.36,274.84) .. (208.33,274.84) .. controls (212.36,274.84) and (214.37,277.72) .. (214.37,283.47)(214.37,280.88) -- (214.37,283.47) .. controls (214.37,287.5) and (216.38,289.51) .. (220.41,289.51) ;
%Shape: Brace [id:dp4720563800156482] 
\draw   (220.78,310.42) .. controls (216.7,310.36) and (214.63,312.37) .. (214.57,316.45) -- (214.57,316.45) .. controls (214.48,322.28) and (212.4,325.16) .. (208.32,325.1) .. controls (212.4,325.16) and (214.4,328.1) .. (214.31,333.93)(214.35,331.31) -- (214.31,333.93) .. controls (214.25,338.01) and (216.26,340.08) .. (220.33,340.14) ;
%Shape: Brace [id:dp09829592935463916] 
\draw   (220.11,360.42) .. controls (216.12,360.41) and (214.12,362.4) .. (214.11,366.39) -- (214.11,366.39) .. controls (214.1,372.08) and (212.1,374.92) .. (208.11,374.91) .. controls (212.1,374.92) and (214.08,377.78) .. (214.07,383.47)(214.08,380.91) -- (214.07,383.47) .. controls (214.06,387.46) and (216.05,389.46) .. (220.04,389.47) ;
%Shape: Brace [id:dp4219363759736008] 
\draw   (270.02,89.96) .. controls (270.11,85.95) and (268.16,83.89) .. (264.15,83.8) -- (264.15,83.8) .. controls (258.42,83.66) and (255.6,81.59) .. (255.69,77.58) .. controls (255.6,81.59) and (252.68,83.52) .. (246.95,83.39)(249.53,83.45) -- (246.95,83.39) .. controls (242.94,83.29) and (240.89,85.25) .. (240.79,89.26) ;
%Shape: Brace [id:dp33643495165734083] 
\draw   (320.02,89.96) .. controls (320.04,85.99) and (318.07,84) .. (314.1,83.98) -- (314.1,83.98) .. controls (308.43,83.95) and (305.61,81.96) .. (305.64,77.99) .. controls (305.61,81.96) and (302.77,83.92) .. (297.11,83.89)(299.66,83.9) -- (297.11,83.89) .. controls (293.14,83.87) and (291.15,85.84) .. (291.13,89.81) ;
%Shape: Brace [id:dp7116045340155585] 
\draw   (370.02,89.96) .. controls (370.09,85.99) and (368.14,83.98) .. (364.17,83.91) -- (364.17,83.91) .. controls (358.5,83.82) and (355.7,81.79) .. (355.77,77.82) .. controls (355.7,81.79) and (352.84,83.72) .. (347.17,83.63)(349.72,83.67) -- (347.17,83.63) .. controls (343.21,83.56) and (341.2,85.5) .. (341.13,89.47) ;
%Shape: Brace [id:dp49438631090080776] 
\draw   (420.02,89.96) .. controls (420.09,85.99) and (418.14,83.98) .. (414.17,83.91) -- (414.17,83.91) .. controls (408.5,83.82) and (405.7,81.79) .. (405.77,77.82) .. controls (405.7,81.79) and (402.84,83.72) .. (397.17,83.63)(399.72,83.67) -- (397.17,83.63) .. controls (393.21,83.56) and (391.2,85.5) .. (391.13,89.47) ;
%Shape: Brace [id:dp5401968132851465] 
\draw   (469.8,89.96) .. controls (469.82,85.98) and (467.84,83.98) .. (463.86,83.96) -- (463.86,83.96) .. controls (458.17,83.93) and (455.34,81.92) .. (455.36,77.94) .. controls (455.34,81.92) and (452.49,83.9) .. (446.8,83.87)(449.36,83.88) -- (446.8,83.87) .. controls (442.81,83.85) and (440.81,85.83) .. (440.79,89.81) ;
%Shape: Brace [id:dp3162229615322666] 
\draw   (520.18,89.97) .. controls (520.25,85.89) and (518.24,83.82) .. (514.16,83.75) -- (514.16,83.75) .. controls (508.33,83.65) and (505.45,81.56) .. (505.52,77.49) .. controls (505.45,81.56) and (502.51,83.55) .. (496.68,83.46)(499.3,83.5) -- (496.68,83.46) .. controls (492.6,83.39) and (490.53,85.4) .. (490.46,89.48) ;
%Straight Lines [id:da08531370665324856] 
\draw  [dash pattern={on 4.5pt off 4.5pt}]  (520.33,140.35) -- (220.67,139.81) -- (220.67,110.47) -- (520.33,111.08) ;
%Straight Lines [id:da44707159043310385] 
\draw  [dash pattern={on 4.5pt off 4.5pt}]  (470.9,189.9) -- (220.67,189.81) -- (220.67,160.47) -- (470.9,160.63) ;
%Straight Lines [id:da6105282934264156] 
\draw  [dash pattern={on 4.5pt off 4.5pt}]  (520.08,240.02) -- (220.67,239.81) -- (220.67,210.47) -- (520.08,210.74) ;
%Straight Lines [id:da9504449883462697] 
\draw  [dash pattern={on 4.5pt off 4.5pt}]  (470.23,289.9) -- (220.67,289.81) -- (220.67,260.47) -- (470.23,260.63) ;
%Straight Lines [id:da9586075981786312] 
\draw  [dash pattern={on 4.5pt off 4.5pt}]  (521.23,339.9) -- (220.67,339.81) -- (220.67,310.47) -- (521.23,310.63) ;
%Straight Lines [id:da4593372659090399] 
\draw  [dash pattern={on 4.5pt off 4.5pt}]  (470.23,389.9) -- (220.67,389.81) -- (220.67,360.47) -- (470.23,360.63) ;
%Straight Lines [id:da019388457576953733] 
\draw  [dash pattern={on 4.5pt off 4.5pt}]  (290.29,339.9) -- (290.66,90.14) -- (319.99,90.14) -- (319.56,339.9) ;
%Straight Lines [id:da06401795857213655] 
\draw  [dash pattern={on 4.5pt off 4.5pt}]  (241,390.02) -- (240.66,90.14) -- (269.99,90.14) -- (270.33,389.35) ;
%Straight Lines [id:da6411259044823248] 
\draw  [dash pattern={on 4.5pt off 4.5pt}]  (390.29,339.9) -- (390.66,90.14) -- (419.99,90.14) -- (419.56,339.9) ;
%Straight Lines [id:da5308006520136684] 
\draw  [dash pattern={on 4.5pt off 4.5pt}]  (341.67,389.69) -- (340.66,90.14) -- (369.99,90.14) -- (369.56,388.9) ;
%Straight Lines [id:da7223505425133505] 
\draw  [dash pattern={on 4.5pt off 4.5pt}]  (490.95,339.9) -- (490.66,90.14) -- (519.99,90.14) -- (520.23,339.9) ;
%Straight Lines [id:da08663622841242158] 
\draw  [dash pattern={on 4.5pt off 4.5pt}]  (440.95,388.9) -- (440.66,90.14) -- (469.99,90.14) -- (470.23,389.9) ;
%Shape: Square [id:dp14989708150943815] 
\draw  [fill={rgb, 255:red, 74; green, 144; blue, 226 }  ,fill opacity=1 ] (390.62,110.63) -- (419.9,110.63) -- (419.9,139.9) -- (390.62,139.9) -- cycle ;
%Shape: Square [id:dp28158066371188273] 
\draw  [fill={rgb, 255:red, 74; green, 144; blue, 226 }  ,fill opacity=1 ] (290.95,110.63) -- (320.23,110.63) -- (320.23,139.9) -- (290.95,139.9) -- cycle ;
%Shape: Square [id:dp44184949966059994] 
\draw  [fill={rgb, 255:red, 74; green, 144; blue, 226 }  ,fill opacity=1 ] (290.5,210.74) -- (319.78,210.74) -- (319.78,240.02) -- (290.5,240.02) -- cycle ;
%Shape: Square [id:dp07559205594674034] 
\draw  [fill={rgb, 255:red, 74; green, 144; blue, 226 }  ,fill opacity=1 ] (390.29,310.63) -- (419.56,310.63) -- (419.56,339.9) -- (390.29,339.9) -- cycle ;
%Shape: Square [id:dp1925726088601144] 
\draw  [fill={rgb, 255:red, 74; green, 144; blue, 226 }  ,fill opacity=1 ] (390.47,210.74) -- (419.75,210.74) -- (419.75,240.02) -- (390.47,240.02) -- cycle ;
%Shape: Square [id:dp6128768535840142] 
\draw  [fill={rgb, 255:red, 74; green, 144; blue, 226 }  ,fill opacity=1 ] (290.29,310.63) -- (319.56,310.63) -- (319.56,339.9) -- (290.29,339.9) -- cycle ;
%Shape: Square [id:dp05578208615521729] 
\draw  [fill={rgb, 255:red, 74; green, 144; blue, 226 }  ,fill opacity=1 ] (490.95,310.63) -- (520.23,310.63) -- (520.23,339.9) -- (490.95,339.9) -- cycle ;
%Shape: Square [id:dp056546956345314814] 
\draw  [fill={rgb, 255:red, 74; green, 144; blue, 226 }  ,fill opacity=1 ] (490.81,210.74) -- (520.08,210.74) -- (520.08,240.02) -- (490.81,240.02) -- cycle ;
%Shape: Square [id:dp9861757911296564] 
\draw  [fill={rgb, 255:red, 74; green, 144; blue, 226 }  ,fill opacity=1 ] (491.06,111.08) -- (520.33,111.08) -- (520.33,140.35) -- (491.06,140.35) -- cycle ;
%Shape: Square [id:dp9926654888956047] 
\draw  [fill={rgb, 255:red, 208; green, 2; blue, 27 }  ,fill opacity=1 ] (240.95,160.63) -- (270.23,160.63) -- (270.23,189.9) -- (240.95,189.9) -- cycle ;
%Shape: Square [id:dp48314171516384075] 
\draw  [fill={rgb, 255:red, 208; green, 2; blue, 27 }  ,fill opacity=1 ] (240.95,260.63) -- (270.23,260.63) -- (270.23,289.9) -- (240.95,289.9) -- cycle ;
%Shape: Square [id:dp6438750747717561] 
\draw  [fill={rgb, 255:red, 208; green, 2; blue, 27 }  ,fill opacity=1 ] (241,360.75) -- (270.28,360.75) -- (270.28,390.02) -- (241,390.02) -- cycle ;
%Shape: Square [id:dp46238192538879974] 
\draw  [fill={rgb, 255:red, 208; green, 2; blue, 27 }  ,fill opacity=1 ] (341.67,360.41) -- (370.94,360.41) -- (370.94,389.69) -- (341.67,389.69) -- cycle ;
%Shape: Square [id:dp4257620523201352] 
\draw  [fill={rgb, 255:red, 208; green, 2; blue, 27 }  ,fill opacity=1 ] (440.95,360.63) -- (470.23,360.63) -- (470.23,389.9) -- (440.95,389.9) -- cycle ;
%Shape: Square [id:dp32141250260218] 
\draw  [fill={rgb, 255:red, 208; green, 2; blue, 27 }  ,fill opacity=1 ] (440.95,260.63) -- (470.23,260.63) -- (470.23,289.9) -- (440.95,289.9) -- cycle ;
%Shape: Square [id:dp11123790054649063] 
\draw  [fill={rgb, 255:red, 208; green, 2; blue, 27 }  ,fill opacity=1 ] (441.62,160.63) -- (470.9,160.63) -- (470.9,189.9) -- (441.62,189.9) -- cycle ;
%Shape: Square [id:dp7231834052718187] 
\draw  [fill={rgb, 255:red, 208; green, 2; blue, 27 }  ,fill opacity=1 ] (341.17,161.22) -- (370.45,161.22) -- (370.45,190.49) -- (341.17,190.49) -- cycle ;
%Shape: Square [id:dp45021402434529945] 
\draw  [fill={rgb, 255:red, 208; green, 2; blue, 27 }  ,fill opacity=1 ] (341.17,260.55) -- (370.45,260.55) -- (370.45,289.83) -- (341.17,289.83) -- cycle ;

% Text Node
\draw (198.53,124.22) node    {$K_{2}$};
% Text Node
\draw (199.87,173.22) node    {$K_{3}$};
% Text Node
\draw (200.2,223.89) node    {$K_{4}$};
% Text Node
\draw (198.53,274.22) node    {$K_{5}$};
% Text Node
\draw (198.53,325.22) node    {$K_{6}$};
% Text Node
\draw (198.87,374.56) node    {$K_{7}$};
% Text Node
\draw (254.61,70.56) node    {$L_{2}$};
% Text Node
\draw (305.28,70.56) node    {$L_{3}$};
% Text Node
\draw (355.28,71.22) node    {$L_{4}$};
% Text Node
\draw (403.28,70.22) node    {$L_{5}$};
% Text Node
\draw (453.95,70.89) node    {$L_{6}$};
% Text Node
\draw (505.28,70.56) node    {$L_{7}$};
% Text Node
\draw (255.59,176.27) node    {$x_{3,2}$};
% Text Node
\draw (255.59,276.27) node    {$x_{5,2}$};
% Text Node
\draw (255.64,376.38) node    {$x_{7,2}$};
% Text Node
\draw (356.3,377.05) node    {$x_{7,4}$};
% Text Node
\draw (455.59,376.27) node    {$x_{7,6}$};
% Text Node
\draw (355.81,276.19) node    {$x_{5,4}$};
% Text Node
\draw (355.81,176.86) node    {$x_{3,4}$};
% Text Node
\draw (456.26,176.27) node    {$x_{3,6}$};
% Text Node
\draw (455.59,276.27) node    {$x_{5,6}$};
% Text Node
\draw (305.59,126.27) node    {$y_{2,3}$};
% Text Node
\draw (305.14,226.38) node    {$y_{4,3}$};
% Text Node
\draw (304.92,326.27) node    {$y_{6,3}$};
% Text Node
\draw (405.26,126.27) node    {$y_{2,5}$};
% Text Node
\draw (505.7,126.72) node    {$y_{2,7}$};
% Text Node
\draw (405.11,226.38) node    {$y_{4,5}$};
% Text Node
\draw (505.44,226.38) node    {$y_{4,7}$};
% Text Node
\draw (404.92,326.27) node    {$y_{6,5}$};
% Text Node
\draw (505.59,326.27) node    {$y_{6,7}$};
% Text Node
\draw (543.26,125.41) node  [font=\large]  {$\cdots $};
% Text Node
\draw (543.26,225.41) node  [font=\large]  {$\cdots $};
% Text Node
\draw (542.26,325.08) node  [font=\large]  {$\cdots $};
% Text Node
\draw (255.93,412.75) node  [font=\large,rotate=-90]  {$\cdots $};
% Text Node
\draw (356.59,412.75) node  [font=\large,rotate=-90]  {$\cdots $};
% Text Node
\draw (457.26,412.75) node  [font=\large,rotate=-90]  {$\cdots $};
% Text Node
\draw (494.93,410.91) node  [font=\Large,rotate=-45]  {$\cdots $};
% Text Node
\draw (545.6,361.58) node  [font=\Large,rotate=-45]  {$\cdots $};
% Text Node
\draw (506.6,372.91) node  [font=\Large,rotate=-45]  {$\cdots $};

\end{tikzpicture}
}
\end{center}

\begin{corollary}\label{arazylemmaforisomorphism*}
Let $E,F$ be two quasi-Banach symmetric sequence spaces such that $c_0,\ell_2\not\hookrightarrow E$ and $\ell_2\not\hookrightarrow F$. 
If $T_F$ is isomorphic to a subspace of $C_E$, then there exists $b\in C_E$ such that the mapping  $W:C_F\to C_E$ defined by $$We_{i,j} = b\otimes e_{i,j},~i,j\ge 1,$$
is an isomorphism. 
\end{corollary}

\begin{proof}
Since $E$ is a quasi-Banach symmetric sequence space, it follows that  (see, e.g., \cite{S14})
 there exists a  positive number $C$ such that
\[
\left\Vert x+y\right\Vert_{C_E}\leq C\left(\left\Vert x\right\Vert_{C_E}+\left\Vert y\right\Vert_{C_E}\right)\;\;\;\;\;\;{\rm for\;any\;}x,y\in C_E,
\]
and
\[
\left\Vert x'+y'\right\Vert_{C_F}\leq C\left(\left\Vert x'\right\Vert_{C_F}+\left\Vert y' \right\Vert_{C_F}\right)\;\;\;\;\;\;{\rm for\;any\;}x',y'\in C_F.
\]
Let $T$ be an  isomorphic operator from $T_F$ into $C_E$.
Applying Theorem \ref{tftoce} to $T$, then there exist two closed subspaces $K$ and $L$ with $T_F(L,K)\subset T_F$, and two operators  $\tilde{R}_1,\tilde{R}_2\in\mathcal B\big(C_F(L,K),C_E\big)$ which satisfy the assertions in Theorem \ref{tftoce} and, in addition, $$R:=(\tilde{R}_1+\tilde{R}_2)|_{T_F(L,K)}$$ is an isomorphic embedding. Next, we will show that $\tilde{R}_2$ is an isomorphic embedding.

By  \eqref{x2k+12l} and \eqref{y2l2k+1}, $\tilde{R}_2$ is an isomorphic embedding if and only if $\tilde{R}_2|_{\left[e^{K,L}_{k,l}\right]_{k,l=N}^\infty}$ is an isomorphic embedding for some positive integer $N$. Note that $\tilde{R}_2|_{\left[e^{K,L}_{k,l}\right]_{k,l=N}^\infty}$ is an isomorphic embedding if and only if $\tilde{R}_2|_{{\rm span}\left(\left\{e^{K,L}_{k,l}\right\}_{k,l=N}^\infty\right)}$ is an isomorphic embedding.
Hence, if $\tilde{R}_2$ is not an isomorphic embedding, then, by induction, there exist an increasing sequence $\{N_n\}_{n=1}^\infty$ of positive integers and a sequence $\{a_n\}_{n=1}^\infty$ of elements in $T_F(L,K)$ such that
\[
\norm{\tilde{R}_2(a_n)}_{C_E}<\frac{1}{\Vert R^{-1}\Vert(2C)^{n+3}},
\]
where $a_n\in [e^{K,L}_{k,l}]_{N_n\leq l\leq k<N_{n+1}}$, for every $n\geq1$. This yields that $\norm{\tilde{R}_2|_{[a_n]_{n=1}^\infty}}\stackrel{\footnotemark}{\leq} C^2\norm{\tilde{R}_2|_{{\rm span}(\{a_n\}_{n=1}^\infty)}}<\frac{1}{2C\Vert R^{-1}\Vert}$\footnotetext{For any $b\in[a_n]_{n=1}^\infty$, choose $b_k\in{\rm span}(\{a_n\}_{n=1}^\infty)$ such that $b_k\to b$ as $k\to\infty$, and hence
\begin{align*}
\norm{\tilde{R}_2(b)}_{C_E}&\leq C\varliminf_{k\to\infty}\norm{ \tilde{R}_2(b_k)}_{C_E}\\
&\leq C\norm{ \tilde{R}_2|_{{\rm span}(\{a_n\}_{n=1}^\infty}}\left(\varliminf_{k\to\infty}\Vert b_k\Vert_{C_F}\right)\leq C^2\norm{ \tilde{R}_2|_{{\rm span}(\{a_n\}_{n=1}^\infty}}\norm{b}_{C_F}.
\end{align*}
}. Consequently,
\[
\left\Vert\tilde{R}_1(x)\right\Vert_{C_E}\geq C^{-1} \left\Vert R(x)\right\Vert_{C_E}-\left\Vert\tilde{R}_2(x)\right\Vert_{C_E}\geq\frac{ \left\Vert x\right\Vert_{C_F}}{2C\left\Vert R^{-1}\right\Vert}\;\;\;\;\;\;
\]
for any $x\in[a_n]_{n=1}^\infty$, and thus $\tilde{R}_1|_{[a_n]_{n=1}^\infty}$ is an isomorphic embedding. Since $\{a_n\}_{n=1}^\infty$ is a diagonal block sequence in $C_F$, it follows that  $[a_n]_{n=1}^\infty$ is isometrically isomorphic to a subspace of $F$. On the other hand, im$(\tilde{R}_1)$ is isomorphic to a Hilbert space, and thus $\ell_2\hookrightarrow F$, contradicting the hypothesis. Consequently, $\tilde{R}_2$ is an isomorphic embedding. Let $a,a'\in C_E$ be such that 
\[
\{x_{2k+1,2l}\}_{k,l=1}^\infty\stackrel{\mbox{{\rm\tiny b.s}}}{\boxplus}\{K_{2k+1}\}_{k=1}^\infty\otimes\{L_{2l}\}_{l=1}^\infty\curvearrowleft a
\]
and
\[
\{y_{2l,2k+1}\}_{k,l=1}^\infty\stackrel{\mbox{{\rm\tiny b.s}}}{\boxplus}\{K_{2l}\}_{l=1}^\infty\otimes\{L_{2k+1}\}_{k=1}^\infty\curvearrowleft a',
\]
then we have
$\big\{x_{2k+1,2l}+y_{2l,2k+1}\big\}_{k,l=1}^\infty$ is isometrically equivalent to
\[
\left\{
\left(
\begin{smallmatrix}
0&a'\\a&0
\end{smallmatrix}\right)\otimes e_{k,l}
\right\}_{k,l=1}^\infty
\]
in $C_E$. This completes the proof.
\end{proof}

The following corollary  treats Pelczynski's question concerning the uniqueness of symmetric structure for a wide class of quasi-Banach operator ideals, which extends Arazy's result\cite{Arazy81}. 
\begin{corollary}

    Let $E,F$ be two quasi-Banach symmetric sequence spaces such that $\ell_2\not\hookrightarrow E,F$, and $c_0\not\hookrightarrow E$ or $F$. If $C_E\hookrightarrow C_F$ and $C_F\hookrightarrow C_E$, then $(e^E_n)_{n\geq1}$ and $(e^F_n)_{n\geq1}$ are equivalent.
    
\end{corollary}
\begin{proof}
Without loss of generality, we may assume that $c_0\not\hookrightarrow E$. By Lemma \ref{arazylemmaforisomorphism}, there exists $b\in C_E$ such that the mapping  $W:C_F\to C_E$ defined by $$We_{i,j} = b\otimes e_{i,j},~i,j\ge 1,$$
is an isomorphism. Since $W|_{ \overline{[e_{i,i}] _{i=1}^\infty}^{\norm{\cdot}_{C_F}}}$ is an isomorphic embedding, it follows that $F$ is  isomorphic to a subspace $E$ and $F\subset E$. Hence, $c_0\not\hookrightarrow F$. Using the fact that $C_E\hookrightarrow C_F$,  the same argument yields that $E\subset F$. This completes the proof.
\end{proof}

  \section{An extension of Sch\"utt's theorem to $L_{W}^p(0,\infty)$}\label{sch-inf}
   Let $1\ge a_1 \ge a_2 \ge \cdots \ge a_{nm} \ge0$, where $n,m=1,2,\cdots$.
   Let $N=nm$.
In this section, we denote the $N$-dimensional weighted Lorentz sequence space by $d_n^m(a,p)$ equipped with the norm 
$$
\left\|x\right\|_{d_n^m(a,p)}=\left(\sum_{i=1}^{nm} a_i |x_i^*|^p\right)^{1/p}
$$
where $x_i^*$ is the decreasing rearrangement of $|x_i|$, $i=1,2,\cdots,nm$.

 For every $M \ge 1$, whenever $W$ is a nonincreasing, positive, continuous function defined on $(0,M)$ and satisfies 
$$
\int_{0}^{M} W(t)\,dt < \infty,
$$
Sch\"utt’s argument can be carried over mutatis mutandis to the weighted Lorentz space $L_{W}^{p}(0,M)$, whose union is dense in $L_W^p(0,\infty)$. 
This implies that $L_W^p(0,\infty)$ is isomorphic to a subspace of $L_p(0,1)$. 
However, the  conclusion of Sch\"utt's theorem  applies only to functions on the interval $(0,1)$.
For completeness, we provide a complete proof of Proposition~\ref{ifvos} in this section.
	
	%\begin{rem}\label{remark}
	%	Here, we present an alternative approach to the proof of (iii) $\Rightarrow$ (i).
		
	%	For $\forall n\in\mathbb N$, $L_{W}^p(0,n)\hookrightarrow L_{W}^{p}(0,\infty)$.
	%	Let $T_W^n(t)=W(nt),\enspace t\in (0,1)$.
	%	It is easy to see that $L_{W}^p(0,n)$ is isometric to $L_{T_W^n}^p(0,1)$.
	%	Indeed, we define a mapping $H_n:L_{W}^p(0,n) \rightarrow L_{T_W^n}^p(0,1)$ with $H_n(f(t))=n^{1/p}T_f(t)$ where $T_f(t)=f(nt),\enspace t\in (0,1)$.
		
	%	By the inequality~(\ref{ineq:ifvos}), we can get that for any $x \in (0,n)$,
	%	\begin{align*}
		%	\int_{0}^{x}T_W^n(\frac{t}{n}) t^{-p/2}\,dt \leq Cx^{-p/2} \int_{0}^{x}T_W^n(\frac{t}{n})\,dt,
	%	\end{align*}
	%	let $s={t/n}$ and $y={x/n}$, we have
	%	\begin{align*}
	%		\int_{0}^{y}T_W^n(s) (ns)^{-p/2} n\,ds \leq C(ny)^{-p/2} \int_{0}^{y}T_W^n(s) n\,ds,
	%	\end{align*}
	%	that is 
	%	\begin{align*}
	%		\int_{0}^{y}T_W^n(s) s^{-p/2} \,ds \leq Cy^{-p/2} \int_{0}^{y}T_W^n(s) \,ds, \quad y \in (0,1).
	%	\end{align*}
	%	This satisfies (iii) in Theorem~\ref{theor:1.1}, so the function space $L_{T_W^n}^{p}(0,1)$ is isomorphic to a subspace of $L_{p}(0,1)$ and then $L_{W}^p(0,n)$ is isomorphic to a subspace of $L_{p}(0,1)$.
	%	This means that for all positive integers n, $L_{W}^p(0,n)$ is uniformly isomorphic to some subspace of $L_{p}$.
	%	Taking $n\rightarrow \infty$, we can get the conclusion.\\
	%	\end{rem}

	\begin{proof}[Proof of Proposition~\ref{ifvos}]
		{\bf (iii)$\Rightarrow$(i):} Firstly, we claim that $L_{W}^p(0,\infty)$ is separable.
		Indeed, by~\cite[Preliminaries]{KamMal} (see also \cite[Theorem 6.3.6]{DPS2023}), it suffices to prove that 
		\begin{align*}
			\int_{0}^{\infty} W(t) \, dt = \infty.
		\end{align*}
		If $\int_{0}^{\infty} W(t) \, dt < \infty$, the right-hand side of inequality~(\ref{ineq:ifvos}) tends to zero when $x \rightarrow \infty$. However, the left-hand side of (\ref{ineq:ifvos}) is evidently strictly positive, which leads to a contradiction.
        
		For any positive integer $M\geq 1$, by \cite[lemma 1.4]{Schutt},  it suffices to show that the subspaces of $L_{W}^p(0,\infty)$ spanned by the vectors $\chi_{[{(i-1)/n},{i/n}]}$, $i=1,2,\ldots,nM$, denoted by $H_n^M$, are uniformly isomorphic to subspaces of $L_p(0,1)$.
	%	Note that if there exists a constant $m>0$ such that 
	%	$$
	%	\int_0^m W(t)\,dt=\infty ,
	%	$$
	%	then $L_{W}^p(0,\infty)=\{f\in L(0,\infty) |\enspace f(x)=0 \enspace a.e., \enspace \forall x\in (0,\infty)\}$.
		%Therefore, $L_{W}^p(0,\infty)\hookrightarrow L_p(0,1)$.
		
	%	Without loss of generality, we assume that 
	%	$$
	%	\int_0^1 W(t)\,dt<\infty.
	%	$$
	For any $x=(x_1,\cdots, x_{nM})\in \mathbb{R}^{nM}$, 	
        we have 
		\begin{align*}
			\left\|\sum_{i=1}^{nM}x_i\chi_{[{(i-1)/n},{i/n}]}\right\|_{L_{W}^p}=\left(\sum_{i=1}^{nM}|x_i^*|^{p}\int_{(i-1)/n}^{i/n}W(t)\,dt \right)^{1/p}.
		\end{align*}
		We put
		\begin{equation}
			a_i:=\int_{(i-1)/n}^{i/n}W(t)\,dt \left(\int_{0}^{1/n}W(t)\,dt \right)^{-1}, \enspace i=1,2,\ldots,nM.
			\label{ineq:3.1}
		\end{equation}
In particular,   $1\ge a_1\ge a_2\ge \cdots\ge a_{nM}>0$.
		For any $x=(x_1,x_2,\ldots,x_{nM}) \in d_n^M(a,p)$, its norm is given by 
		\begin{align*}
			\left\|x\right\|_{d_n^M(a,p)}&=\left(\sum_{i=1}^{nM}a_i|x_i^*|^{p}\right)^{1/p}\\
			&=\left(\sum_{i=1}^{nM}|x_i^*|^{p}\int_{(i-1)/n}^{i/n}W(t)\,dt \left(\int_{0}^{1/n}W(t)\,dt\right)^{-1}
			\right)^{1/p}.
		\end{align*}
        Moreover, $d_n^M(a,p)$  is isometric to $H^M_n$. 
		%We define a mapping $T_n^M:H_n^M \rightarrow d_n^M(a,p)$ by the rule 
		%\begin{align*}
		%	T_n^M\left(\sum_{i=1}^{nM}x_i\chi_{[{(i-1)/n},{i/n}]}\right)=\left(A_nx_1,A_nx_2,\ldots,A_nx_{nM}\right),
	%	\end{align*}
	%	while $A_n=\left(\int_{0}^{1/n}W(t)\,dt\right)^{1/p}$.
	%	And $T_n^M$ is an isometric isomorphism.
        
		For all   $1\leq k \leq nM$, by  (\ref{ineq:ifvos}),  we have 
		\begin{align*}
			\sum_{i=1}^{k}\int_{(i-1)/n}^{i/n}W(t)(\frac{i}{n})^{-p/2}\,dt&\leq
			\sum_{i=1}^{k}\int_{(i-1)/n}^{i/n}W(t)t^{-p/2}\,dt\\ &=
			\int_{0}^{k/n}W(t)t^{-p/2}\,dt\\ &\leq
			C(\frac{k}{n})^{-p/2}\int_{0}^{k/n}W(t)\,dt\\ &=C(\frac{k}{n})^{-p/2}\sum_{i=1}^{k}\int_{(i-1)/n}^{i/n}W(t)\,dt.
		\end{align*}
		That is, 
		\begin{align*}
			\sum_{i=1}^{k}{i}^{-p/2}\int_{(i-1)/n}^{i/n}W(t)\,dt&\leq C{k}^{-p/2}\sum_{i=1}^{k}\int_{(i-1)/n}^{i/n}W(t)\,dt.
		\end{align*}
		Multiplying both sides of the above inequality by $\left(\int_{0}^{1/n}W(t)\,dt \right)^{-1}$, we get
		\begin{equation}
			\sum_{i=1}^{k}a_i i^{-p/2} \leq Ck^{-p/2}\sum_{i=1}^{k}a_i.
			\label{ineq:3.2}
		\end{equation}
		For $1\leq i\leq k/{(2C)^{2/p}}$, we have  $i^{-p/2}\geq2Ck^{-p/2}$.
		Therefore,
		\begin{align*}
			2Ck^{-p/2}\sum_{1\le i\leq k/{(2C)^{2/p}}}a_i\leq \sum_{1\le i\leq k/{(2C)^{2/p}}}a_i i^{-p/2}.
		\end{align*}
		Therefore,
		\begin{align*}
			2Ck^{-p/2}\sum_{1\le i\leq k/{(2C)^{2/p}}}a_i \quad+ \sum_{1\le i\leq k/{(2C)^{2/p}}}a_i i^{-p/2} \leq 
			2\sum_{1\le i\leq k/{(2C)^{2/p}}}a_i i^{-p/2},
		\end{align*}
		that is, 
		\begin{align}\label{aii0p2}\begin{split}
	&		\sum_{1\le i\leq k/{(2C)^{2/p}}}a_i i^{-p/2}\\ &\leq 2\left( \sum_{1\le i\leq k/{(2C)^{2/p}}}a_i i^{-p/2} \quad - Ck^{-p/2}\sum_{1\le i\leq k/{(2C)^{2/p}}}a_i \right) \\ &=
			2\sum_{1\le i\leq k/{(2C)^{2/p}}} a_i \left(i^{-p/2} - C k^{-p/2}\right).\end{split}
		\end{align}
		On the other hand,
		\begin{align}\label{aii-p2}\begin{split}
			\sum_{1\le i\leq k/{(2C)^{2/p}}} a_i i^{-p/2} &\leq \sum_{i=1}^{k}a_i i^{-p/2} \\
			&\stackrel{(\ref{ineq:3.2})}{\leq}
			Ck^{-p/2}\sum_{i=1}^{k}a_i \\
			& =
			Ck^{-p/2} \left(\sum_{1\le i\leq k/{(2C)^{2/p}}} a_i+\sum_{i>{k/(2C)^{2/p}}}^{k}a_i\right)\end{split}
		\end{align}
		This means in particular that 
		\begin{align*}
			\frac{1}{2} \sum_{1\le i\leq k/{(2C)^{2/p}}} a_i i^{-p/2} \stackrel{\eqref{aii0p2}}{\leq} \sum_{1\le i\leq k/{(2C)^{2/p}}} a_i \left(i^{-p/2} - C k^{-p/2}\right)
			\stackrel{\eqref{aii-p2}}{\leq} Ck^{-p/2}\sum_{i> {k/(2C)^{2/p}}}^{k}a_i.
		\end{align*}
		Then, we obtain 
		\begin{align*}
			\frac{1}{2} \sum_{1\le i\leq k/{(2C)^{2/p}}} a_i i^{-p/2} &\leq Ck^{-p/2}\sum_{i> {k/(2C)^{2/p}}}^{k}a_i\\
			&\leq Ck^{-p/2}(k-[{k/(2C)^{2/p}}])a_{[{k/(2C)^{2/p}}]}\\
			&\leq Ck^{-p/2}\cdot k\cdot a_{[{k/(2C)^{2/p}}]} \\
			&\leq Ck^{(1-{p/2})}a_{[{k/(2C)^{2/p}}]},
		\end{align*}
        where $[x]:=\max\{n\in \mathbb{N}|\enspace n\le x\}$.
		We define $t:=[k/(2C)^{2/p}] $. In particular, $1 \leq t \leq \left(nM/(2C)^{2/p}\right)$.
	Then, we have 
		\begin{align*}
			\sum_{1\le i \leq t}a_i i^{-p/2}\leq 2C\left((2C)^{2/p}t\right)^{1-p/2}a_t.
		\end{align*}
		We denote $C':=(2C)^{2/p}$. Hence, we have  
		\begin{align}\label{mainineq}\begin{split}
			\sum_{i=1}^{t}a_i i^{-p/2}\leq C't^{1-p/2}a_t.
             \end{split}
		\end{align}
		%This implies  {\color{red}t}
		%\begin{align*}
		%	\sum_{i=1}^{k}a_i i^{-p/2} \leq C' k^{1-{p/2}}a_k,\quad 1\leq k \leq\frac{nM}{(2C)^{2/p}},
	%	\end{align*}
		and $C'$ depends only on $C$.
		This allows us to   establish  the inequality in the assumption of \cite[Proposition 1.2]{Schutt}. Indeed, 
		for $1\leq p <2$, we have  $0 < (1-p/2) \leq1/2 <1$ and therefore,   for any $t=1,2,\ldots,\left[\frac{nM}{(2C)^{2/p}}\right]$, we have 
		\begin{align*}
			a_1+\left(1-\frac{p}{2}\right)\sum_{i=1}^{t}a_i i^{-p/2}&\leq a_1+\sum_{i=1}^{t}a_i i^{-p/2}\\
			&\leq
			\sum_{i=1}^{t}a_i i^{-p/2}+\sum_{i=1}^{t}a_i i^{-p/2}\\
			& =
			2\sum_{i=1}^{t}a_i i^{-p/2}\\
			&\stackrel{\ref{mainineq}}{\leq}
			2C' t^{1-{p/2}}a_t\\
			& =
			C''t^{1-{p/2}}a_t   ,
		\end{align*}
        where $C'': = 2C'$.
	%	Therefore,
	%	\begin{equation}
	%		a_1+\left(1-\frac{p}{2}\right)\sum_{i=1}^{k}a_i i^{-p/2} \leq C''k^{1-{p/2}}a_k,\quad k=1,2,\ldots,\frac{nM}{(2C)^{2/p}},
		%	
	%	\end{equation}
	This shows that the assumption of  \cite[Proposition 1.2]{Schutt} is satisfied.
		Since  $n$ and $M$ are arbitrarily taken, it follows   the space $H_n^M$'s, $n,m\ge 1$, are uniformly isomorphic to subspaces of $L_p(0,1)$. Hence, $L _{W}^p(0,\infty)$ is isomorphic to a subspace of $L_p(0,1)$.

		{\bf (i)$\Rightarrow$(ii):} Assume that  $L_{W}^p(0,\infty)$ is isomorphic to a subspace of $L_p(0,1)$. Since $L_p(0,1)$ has 
         cotype $2$, it follows that $L_{W}^p(0,\infty)$ also has cotype $2$.
         This is equivalent to $L_{W}^p(0,\infty)$    being $2$-concave\cite[Theorem 1.f.16]{LT2}.

		{\bf (ii)$\Rightarrow$(iii):}    Suppose that $L_{W}^{p}(0,\infty)$ is $2$-concave.   Recall the fact that  a Banach lattice has cotype $2$ if and only if it is $2$-concave\cite[Theorem 1.f.16]{LT2}.
	Observe that    the constants for the $2$-concavity of the spaces $d_n^M(a,p)$ generated by (\ref{ineq:3.1}) are uniformly bounded.  
		We have 
		\begin{eqnarray*}
			\int_{\frac{1}{n}}^{\frac{k}{n}} W(t) t^{-p/2} \, dt
			&= &\sum_{i=2}^{k} \int_{\frac{i-1}{n}}^{\frac{i}{n}} W(t) t^{-p/2} \, dt \\
			&\leq &\sum_{i=2}^{k} \left(\frac{i-1}{n}\right)^{-p/2} \int_{\frac{i-1}{n}}^{\frac{i}{n}} W(t) \, dt\\
			&\stackrel{\eqref{ineq:3.1}}{ = }&\int_0^{\frac{1}{n}} W(t) \, dt \left( \sum_{i=2}^{k} \left(\frac{i-1}{n}\right)^{-p/2}\cdot a_i \right) \\
			&\le &\int_0^{\frac{1}{n}} W(t) \, dt \left(  2\sum_{i=2}^{k} \left(\frac{i}{n}\right)^{-p/2}\cdot a_i \right)\\
			&\leq& \int_0^{\frac{1}{n}} W(t) \, dt \left( 2\sum_{i=1}^{k} \left(\frac{i}{n}\right)^{-p/2}\cdot a_i \right) .
		\end{eqnarray*}
		By \cite[proposition 1.3]{Schutt}, we have 
        \begin{align}\label{prop1.3}
             2\sum_{i=1}^{k} \left(\frac{i}{n}\right)^{-p/2}\cdot a_i \le   C_p ({\frac{k}{n}})^{-p/2} \sum_{i=1}^{k}a_i
        \end{align}
        for all $1 \leq k \leq nM<\infty$, where $C_p$ is a constant depending on $p$ and the $2$-concavity constants of $d_n^M(a,p)$ which are uniformly bounded.
        Therefore, we have 
		\begin{align*}
        \int_{\frac{1}{n}}^{\frac{k}{n}} W(t) t^{-p/2} \, dt
			\stackrel{\eqref{prop1.3}}{ \le } \int_0^{\frac{1}{n}} W(t) \, dt \left( C_p ({\frac{k}{n}})^{-p/2} \sum_{i=1}^{k}a_i \right) \stackrel{\eqref{ineq:3.1}}{ = }  C_p ({\frac{k}{n}})^{-p/2} \int_{0}^{\frac{k}{n}} W(t) \,dt.
		\end{align*}
		The proof is complete by taking $n, M\rightarrow \infty$. 
	\end{proof}

Note that $L_q(0,1)\hookrightarrow L_r(0,1)$ for $1\le r\le q \le2$ \cite[Proposition 11.1.9]{AK}.
      We obtain the following corollary.
	\begin{corollary}\label{aoifvos}
		Let $1\leq r<2$ and $0<p,q\le \infty $.
		The Lorentz spaces $L_{p,q}(0,\infty)$ are isomorphic to subspaces of $L_r(0,1)$ if and only if 
		$1\leq r\leq q\leq p<2$ or $p=q=2$.
	\end{corollary} 
		\begin{proof}

			$(\Rightarrow)$: If $L_{p,q}(0,\infty)\hookrightarrow L_r(0,1)$, we have $$L_{p,q}(0,1)\hookrightarrow L_{p,q}(0,\infty)\hookrightarrow L_r(0,1)\hookrightarrow L_1(0,1).$$
            Note that $L_r(0,1)$ is separable, while $L_{\infty}(0,\infty)$ and $L_{p,\infty}(0,\infty)$ are non-separable.
            Therefore, it suffices to consider the case when $0<p,q<\infty$.
            By \cite[Theorem 5.1]{Maligranda}, it satisfies $1\le q\le p<\infty$ or $1<p<q<\infty$.
		By \cite[Corollary~1.5]{Schutt}, we get that $$ 1\leq r\le q\leq p<2~\mbox{ or }~p=q=2.$$
			
			$(\Leftarrow)$: If $1\leq r\le q\leq p<2$,  then by \cite[Theorem 3.5]{Creekmore}, $L_{p,q}(0,\infty)$ is of  cotype $2$ .
            Let $$W(s):=\frac{q}{p}s^{q/p-1} ,$$ which is nonincreasing, positive and continuous with $\int_0^{\infty} W(s)\,ds=\infty$.
			In particular, $L_W^q(0,\infty)=L_{p,q}(0,\infty)$.
			By $(i)$ and $(ii)$ in Proposition~\ref{ifvos}, we have $$L_{p,q}(0,\infty)\hookrightarrow L_q(0,1) .$$
            Therefore, $L_{p,q}(0,\infty)\hookrightarrow L_q(0,1)\hookrightarrow L_r(0,1)$.
			If $p=q=2$, then $L_2(0,\infty)\approx L_2(0,1)\hookrightarrow L_r(0,1)$~\cite[Proposition 11.1.9]{AK}.
		\end{proof}

\section{$d(w,p)$ with exactly two nonequivalent symmetric basic sequences}

In this section, we give a necessary and sufficient condition that $d(w,p)$, $0<p<\infty$, has exactly two nonequivalent symmetric basic sequences.
This generalizes the work of Altshuler, Casazza, and Lin \cite[Theorem 6]{ACL} and follows their established line of reasoning. 

\begin{definition}
    Let $\{s_n\}_{n\ge 1}$ and $\{t_n\}_{n\ge 1}$ be two sequences of nonnegative numbers. 
    We say that $\{t_n\}_{n\ge 1}$ dominates $\{s_n\}_{n\ge 1}$, denoted by $t_n>s_n$, if there exists a positive number $A$ such that $s_n\le At_n$, $n=1,2,\cdots$.
    We say that $\{s_n\}_{n\ge 1}$ is equivalent to $\{t_n\}_{n\ge 1}$, and write $s_n\sim t_n$, if $s_n>t_n$ and $t_n>s_n$.
\end{definition}

\begin{prop}\label{comequiva}\cite[Proposition 7]{ACL}\cite[p. 56, Proposition 3.6]{BS}
    Let $\{v_i\}_{i\ge 1}$ and $\{w_i\}_{i\ge1}$ be sequences of nonnegative numbers and let $s_n=\sum_{i=1}^n v_i$, $t_n=\sum_{i=1}^n w_i$, $n=1,2,\cdots$.
    Then $t_n >s_n$ (for a positive number $A$) if and only if   $\sum_{i=1}^{\infty} \beta_i v_i\le A\sum_{i=1}^{\infty} \beta_i w_i$ for all nonincreasing sequences $\{\beta_i\}_{i\ge 1}$ of nonnegative numbers.
\end{prop}
%\begin{proof}
%    $(\Rightarrow)$: If $t_n >s_n$, then there exists $A>0$ such that $s_n\le At_n$, $n=1,2,\cdots$.
%    Let $s_0=0$.
%    For any nonincreasing sequences $\{\beta_i\}_{i\ge 1}$ of nonnegative numbers and $N=1,2,\cdots$, 
%    \begin{align*}
%            \sum_{i=1}^N \beta_i v_i
%            &=
%            \sum_{i=1}^N \beta_i (s_{i}-s_{i-1})\\
%            &=
%            \beta_1 s_1+\beta_2(s_2-s_1)+\cdots+\beta_N(s_N-s_{N-1})\\
%            &=
%            s_1(\beta_1-\beta_2)+s_2(\beta_2-\beta_3)+\cdots+s_{N-1}(\beta_{N-1}-\beta_N)+s_N\beta_N\\
%            &=
%            \sum_{i=1}^{N-1} s_i (\beta_i-\beta_{i+1})+s_N\beta_N.
%     \end{align*}
%     Hence, 
%        \begin{align*}
%        \sum_{i=1}^N \beta_i v_i
%            \le
%            A\left(\sum_{i=1}^{N-1} t_i (\beta_i-\beta_{i+1})+t_N\beta_N\right)
%            =
%            A\sum_{i=1}^N \beta_i w_i.
%    \end{align*}
%Let $N\rightarrow\infty$, we have  $\sum_{i=1}^{\infty} \beta_i v_i\le A\sum_{i=1}^{\infty} \beta_i w_i$.
%
%$(\Leftarrow)$: For $n=1,2,\cdots$, let $\beta_i=1$ for $1\le i \le n$ and all remaining entries are zero.
%This completes the proof.
%\end{proof}

%It is worth noting that the number $A$ appearing in Definition~\ref{sequequiva}
%coincides with the number $A$ 
%appearing in Proposition~\ref{comequiva}.
Below,  we extend \cite[Lemma 2]{ACL} to the case $0<p<\infty $, which shows that two weighted Lorentz sequence spaces are isomorphic if and only if the partial sums of their respective weight sequences are equivalent.
\begin{lemma}\label{isomorequiva}
    Let $0<p<\infty$.
    Let $w_1=\{w_n^1\}_{n\geq 1}$ and $w_2=\{w_n^2\}_{n\ge 1}$ be nonincreasing sequences of positive numbers.
    For each $n=1,2,\cdots$, let $t_n=\sum_{i=1}^n w_i^1$ and $s_n=\sum_{i=1}^n w_i^2$.
    Then, $d(w_1,p)$ is isomorphic to $d(w_2,p)$ if and only if $t_n\sim s_n$.
\end{lemma}
\begin{proof}
    Let $\{x_n\}_{n\ge 1}$ and $\{y_n\}_{n\ge }$ be the unit vector basis of $d(w_1,p)$ and $d(w_2,p)$, respectively.
    By Proposition~\ref{comequiva}, $t_n\sim s_n$ if and only if there exist $A>0$ ,$B>0$ such that $A\sum_{i=1}^{\infty} \alpha_i^p w_i^2\le \sum_{i=1}^{\infty} \alpha_i^p w_i^1 \le B\sum_{i=1}^{\infty} \alpha_i^p w_i^2 $ for every nonincreasing and nonnegative sequence $\{\alpha_i\}_{i\ge 1}$.
    That is, 
$$    A \norm{\sum_{i=1}^\infty \alpha_i y _i }_{d(w_2,p)}^p\le   \norm{\sum_{i=1}^\infty \alpha_i x _i }_{d(w_1,p)}^p  \le B\norm{\sum_{i=1}^\infty \alpha_i y_i }_{d(w_2,p)}^p. $$
    Since $\{x_n\}_{n\ge 1}$ and $\{y_n\}_{n\ge }$ are symmetric, this means that $t_n\sim s_n$ if and only if $\{x_n\}_{n\ge 1} \sim \{y_n\}_{n\ge 1}$.
   Moreover, by \cite[Theorem 4]{ACL} and \cite[Proposition 3]{NO}, $\{x_n\}_{n\ge 1}$ and $\{y_n\}_{n\ge }$ are unique, up to equivalence, symmetric bases in $d(w_1,p)$ and $d(w_2,p)$ respectively.
   That is,   $d(w_1, p)$ is isomorphic to $d(w_2,p)$ if and only if $\{x_n\}_{n\ge 1} \sim \{y_n\}_{n\ge 1}$, i.e., if and only if $t_n\sim s_n$.
\end{proof}

 Let $n\ge 1$. The dilation operator  $\sigma_n$ is defined by  
$$
\sigma_n\alpha=(\underbrace{\alpha_1,\cdots,\alpha_1}_{\text{n elements}},\underbrace{\alpha_2,\cdots\alpha_2}_{\text{n elements}},\cdots,\underbrace{\alpha_k,\cdots,\alpha_k}_{\text{n elements}},\cdots), ~\forall \alpha=(\alpha_1,\alpha_2,\cdots)\in \ell_\infty .
$$

\begin{definition}\label{definition m}
    Let $0\not=\alpha=\sum_{i=1}^{\infty} \alpha_ie_i^{d(w,p)}\in d(w,p)$ with $\alpha_1\ge \alpha_2 \ge \ldots \ge 0$.
    For each $n=1,2,\cdots$, let $s_n=\sum_{i=1}^n w_i$, $s_n^{(\alpha)}=\sum_{i=1}^{\infty} \alpha_i^p(s_{ni}-s_{n(i-1)})
    =\norm{ \sigma_n \alpha  }_{d(w,p)}^p 
    $, where $\sigma_n$ is the dilation operator,  $s_0=0$, $m_n^{(\alpha)}=s_n^{(\alpha)}-s_{n-1}^{(\alpha)}$ and $s_0^{(\alpha)}=0$.
\end{definition}

\begin{proposition}\label{decreastozero}
      Let $0<p<\infty$ and $d(w,p)$ be the   weighted Lorentz sequence space where $w=\{w_n\}_{n\ge 1}$ is strictly decreasing to $0$.
    For $0\not=\alpha=\sum_{i=1}^{\infty} \alpha_i e_i^{d(w,p)}\in d(w,p)$, let $v_i=\sum_{j=1}^{\infty} \alpha_j^p w_{i,j}$, where $\{w_{i,j}\}_{j=1,2,\ldots}$ (respectively, $\{w_{i,j}\}_{i=1,2,\ldots}$) is a subsequence of $\{w_n\}_{n\ge 1}$ for $i=1,2,\cdots$ (respectively, $j=1,2,\cdots$).
    Then $\{v_i\}_{i\ge 1}$ is decreasing to zero.
\end{proposition}

\begin{proof}
    It follows from the hypothesis that $w_{i,j}\le w_{i,h}$, $w_{j,i}\le w_{h,i}$ for $j\ge h$ and $w_{i,j}\le w_i$ for $i=1,2,\cdots$.
    Therefore, $\{v_i\}_{i\ge 1}$ is decreasing.
    For any $\epsilon>0$, choose $N$ and $M$ such that $\sum_{j=N+1}^{\infty} \alpha_j^p w_j<\epsilon/2$ and $w_M \sum_{j=1}^N \alpha_j^p <\epsilon/2$.
    Then for every $i\ge M$, 
    $$v_i=\sum_{j=1}^N \alpha_j^p w_{i,j} + \sum_{j=N+1}^{\infty} \alpha_j^p w_{i,j}
    \le
    w_M\sum_{j=1}^N \alpha_j^p + \sum_{j=N+1}^{\infty} \alpha_j^p w_j
    <\epsilon
    $$
    which means $v_i\rightarrow 0$ as $i\rightarrow \infty$.
\end{proof}

\begin{proposition}\label{weightfunc}
      Let $0<p<\infty$ and $d(w,p)$ be the   weighted Lorentz sequence space where $w=\{w_n\}_{n\ge 1}$ is strictly decreasing to $0$.
    For every $\alpha\in d(w,p)$,  we have $s_n^{(\alpha)}\le n\norm{\alpha}_{d(w,p)}^p$, $\sum_{n=1}^{\infty} m_n^{(\alpha)}=+\infty$ and $\{m_n^{(\alpha)}\}_{n\ge 1}$ is a sequence decreasing to zero, where $s_n^{(\alpha)}$ and $m_n^{(\alpha)}$ are defined in Definition~\ref{definition m}.
\end{proposition}
\begin{proof}
Obviously, $$s_n^{(\alpha)}=\sum_{i=1}^{\infty} \alpha_i^p(s_{ni}-s_{n(i-1)})\le n\norm{\alpha}_{d(w,p)}^p$$ and $$\sum_{n=1}^{\infty} m_n^{(\alpha)}=\lim_{n\to \infty }s_n^{(\alpha)} = \infty.$$
For a fixed $n$, we have $s_{(n+1)k}\ge s_{nk}$ and $2s_{nk}
%=\sum_{i=1}^{(n-1)k} w_i +\sum_{i=(n-1)k+1}^{nk} %w_i +\sum_{i=1}^{nk}w_i
\ge 
%s_{(n-1)k}+\sum_{i=nk+1}^{(n+1)k} w_i+ \sum_{i=1}^{nk} w_i=
s_{(n-1)k}+s_{(n+1)k}$ for $k=1,2,\cdots$.
Thus, 
$$\sum_{i=1}^k \left(s_{(n+1)i}-s_{(n+1)(i-1)}\right)=s_{(n+1)k}\ge s_{nk}=\sum_{i=1}^k \left( s_{ni}-s_{n(i-1)}\right)$$ 
and 
\begin{align*}
2\sum_{i=1}^k \left( s_{ni}-s_{n(i-1)}\right)
&=
2s_{nk}\\
&\ge
s_{(n-1)k}+
s_{(n+1)k}\\
&=
\sum_{i=1}^k \left(s_{(n-1)i}-s_{(n-1)(i-1)}\right)+
\sum_{i=1}^k \left(s_{(n+1)i}-s_{(n+1)(i-1)}\right)\\
\end{align*}
By Proposition~\ref{comequiva} with $\beta_i=\alpha_i^p $, we have $s_{n+1}^{(\alpha)}\ge s_n^{(\alpha)}$  and $2s_n^{(\alpha)}\ge s_{n-1}^{(\alpha)}+s_{n+1}^{(\alpha)}$.
Hence, $\{m_n^{(\alpha)}\}_{n\ge 1}$ is a decreasing sequence of nonnegative numbers.
Now,
\begin{align*}
    m_n^{(\alpha)}
    &=s_n^{(\alpha)}-s_{n-1}^{(\alpha)}\\
    &=\sum_{i=1}^{\infty} \alpha_i^p\left(s_{ni}-s_{n(i-1)}\right)-\sum_{i=1}^{\infty} \alpha_i^p\left(s_{(n-1)i}-s_{(n-1)(i-1)}\right)\\
    &=\sum_{i=1}^{\infty} \alpha_i^p \left(\sum_{j=1}^n w_{n(i-1)+j}\right)-\sum_{i=1}^{\infty} \alpha_i^p \left(\sum_{j=1}^{n-1} w_{(n-1)(i-1)+1}\right)\\
   & =\sum_{i=1}^{\infty} \alpha_i^p \left(\sum_{j=1}^n w_{n(i-1)+j}-\sum_{j=1}^{n-1} w_{(n-1)(i-1)+1}\right)\\
    &\le
    \sum_{i=1}^{\infty}\alpha_i^p w_{ni} =v_n \stackrel{\tiny \rm Prop.~\ref{decreastozero}}{\to} 0 
\end{align*}
     as $n\rightarrow \infty$.
\end{proof}

Let $m^{(\alpha)}:=\left(m_1^{(\alpha)}, m_2^{(\alpha)}, \cdots\right)$. 
By Proposition~\ref{weightfunc}, $d\left(m^{(\alpha)},p\right)$ is a weighted Lorentz space for $0<p<\infty$.
By Lemma~\ref{isomorequiva}, we  provide a necessary and sufficient condition for the isomorphism of $d\left(m^{(\alpha)},p\right)$ and $d(w,p)$.
\begin{proposition}\label{mequiva}
      Let $0<p<\infty$ and $d(w,p)$ be the   weighted Lorentz sequence space where $w=\{w_n\}_{n\ge 1}$ is strictly decreasing to $0$.
     Then $d\left(m^{(\alpha)},p\right)$ is isomorphic to $d(w,p)$ if and only if $s_n^{(\alpha)} \sim s_n$.
\end{proposition}

\begin{lemma}\label{lessinequi}
         Let $0<p<\infty$ and $d(w,p)$ be the   weighted Lorentz sequence space where $w=\{w_n\}_{n\ge 1}$ is strictly decreasing to $0$.
    If for every $\alpha \in d(w,p)$ with $\alpha_1\ge \alpha_2 \ge \cdots \ge 0$ and $\norm{\alpha}_{d(w,p)}=1$, there exists $B_{\alpha}>0$ such that $s_n^{(\alpha)}\le B_{\alpha} s_n$ for $n=1,2,\cdots$, then there exists $B>0$ such that for all $\alpha \in d(w,p)$ with $\norm{\alpha}_{d(w,p)}=1$, $s_n^{(\alpha)}\le B s_n$, $n=1,2,\cdots$.
\end{lemma}
\begin{proof}
    For every fixed $n$, let $w_k^{(n)}:=\frac{s_{nk}-s_{(n-1)k}}{s_n}$, $k=1,2,\cdots$.
    Then $w^{(n)}:=\left(w_1^{(n)}, w_2^{(n)}, \cdots,w_k^{(n)}, \cdots\right)\in c_0\setminus \ell_1$.
    Let $d: =\left(\oplus_{n=1}^{\infty} d\left(w^{(n)},p\right)\right)_{c_0}$ and let $\{e_i^{(n)}\}_{i\ge 1}$ and $\{e_i^{d(w,p)}\}_{i\ge 1}$ be the unit vector basis of $d\left(w^{(n)},p\right)$ and $d(w,p)$, respectively.
    Define $T_n:d(w,p)\rightarrow d$ by
    $$
    T_n\left(\sum_{i=1}^{\infty} \alpha_ie_i^{d(w,p)}\right):= (0,\cdots,0,\underbrace{\sum_{i=1}^{\infty} \alpha_i e_i^{(n)}}_{\text{n'th place}},0,\cdots)
    $$
    for all $\alpha=\sum_{i=1}^{\infty} \alpha_ie_i^{d(w,p)} \in d(w,p)$.
    If $\alpha_1\ge \alpha_2 \ge \cdots \ge 0$, then we have 
    \begin{align}\label{norm T}
    \norm{T_n(\alpha)}_{d}^p=\sum_{i=1}^{\infty} \frac{\alpha_i^p\left(s_{nk}-s_{(n-1)k}\right)}{s_n}=\frac{s_n^{(\alpha)}}{s_n}\stackrel{\rm \tiny Prop. \ref{weightfunc}}{\le} \frac{n}{s_n}\norm{\alpha}_{d(w,p)}^p.
    \end{align}
    and so $\norm{T_n}\le (n/s_n)^{1/p}$.
    Now, for each $\alpha \in d(w,p)$, by the hypothesis that $s_n^{(\alpha)}\le B_{\alpha} s_n$, we have  $\norm{T_n(\alpha)}_d=\left(s_n^{(\alpha)}/s_n\right)^{1/p}\le (B_{\alpha})^{1/p}$ for all $n=1,2,\cdots$.
    By the uniform boundedness principle\cite{KPR}, there exists  $B>0$ such that $\sup_{n}\norm{T_n}\le B^{1/p}$.
    Thus, for every $\alpha \in d(w,p)$ with $\norm{\alpha}_{d(w,p)}=1$, by \eqref{norm T},
    we get $s_n^{(\alpha)}\le B s_n$, $n=1,2,\cdots$.
\end{proof}

The following theorem provides a necessary and sufficient condition for $\{ s_n^{(\alpha)}\}_{n\ge 1}$ to be equivalent with $\{ s_n\}_{n\ge 1}$, which extends \cite[Theorem 5]{ACL}.
\begin{theorem}\label{supequiva}
     Let $0<p<\infty$ and $d(w,p)$ be the   weighted Lorentz sequence space where $w=\{w_n\}_{n\ge 1}$ is strictly decreasing to $0$.
     Then $\sup_{1\le n,k<\infty} s_{nk}/s_ns_k <\infty$ if and only if for every $\alpha \in d(w,p)$ with $\alpha_1\ge \alpha_2 \ge \cdots \ge 0$ and $\norm{\alpha}_{d(w,p)}=1$, $s_n^{(\alpha)} \sim s_n$.
\end{theorem}
\begin{proof}
   $(\Rightarrow)$: Let $\alpha=\sum_{i=1}^{\infty} \alpha_ie_i^{d(w,p)}\in d(w,p)$ such that $\norm{\alpha}_{d(w,p)}=1$.
By definition,  we always have $s_n^{(\alpha)}\ge \alpha_1^ps_n$.
    Suppose that $\sup_{1\le n,k<\infty} s_{nk}/s_ns_k =B<\infty$.
    Then, $s_{nk}\le Bs_ns_k$ for all $1\le n,k,<\infty $, that is $$\sum_{i=1}^{k} \left(s_{ni}-s_{n(i-1)}\right) = s_{nk}\le Bs_ns_k= Bs_n\left(\sum_{i=1}^k w_i\right)$$ for all $n,k=1,2,\cdots$.
    For any fixed  $n$,
    by Proposition~\ref{comequiva}, we get $$Bs_n=Bs_n\left(\sum_{i=1}^{\infty}\alpha_i^pw_i\right)\ge \sum_{i=1}^{\infty} \alpha_i^p \left( s_{ni}-s_{n(i-1)}\right)=s_n^{(\alpha)}.$$
    Hence, $s_n^{(\alpha)} \sim s_n$.

    $(\Leftarrow)$: Suppose $s_n^{(\alpha)} \sim s_n$ for all $\alpha \in d(w,p)$ with $\alpha_1\ge \alpha_2 \ge \cdots \ge 0$ and $\norm{\alpha}_{d(w,p)}=1$.
    By Lemma~\ref{lessinequi}, there exists $B>0$ such that $s_n^{(\alpha)}\le Bs_n$, $n=1,2,\cdots$.
    For each $k$, let $\gamma_i=(1/s_k)^{1/p}$ if $i\le k$ and $\gamma_i=0$ if $i>k$.
    Let $\gamma=\sum_{i=1}^{\infty} \gamma_ie_i^{d(w,p)}$.
    Then, \begin{align*}
  \norm{\gamma}_{d(w,p)}=\left(\frac{\sum_{i=1}^k w_i}{s_k}\right)^{1/p}=1 \enspace\text{and} \enspace s_n^{(\gamma)}=\sum_{i=1}^{\infty} \gamma_i^p(s_{ni}-s_{n(i-1)})=s_{nk}/s_k.
    \end{align*} 
    Hence, $s_{nk}\le Bs_ns_k$, $n,k=1,2,\cdots$.
\end{proof}
The following   property is well-known for symmetric sequence spaces, see \cite[Proposition 4]{ACL}. The proof for quasi-Banach symmetric sequence spaces follows from the same argument. Therefore, the proof of the following lemma is omitted.  
\begin{lemma}\label{dominatelemma}
    Let $\{e_n\}_{n\ge1}$ be a  symmetric basis in a quasi-Banach symmetric sequence  space $X$.
    If $y_n=\sum_{i=p_n+1}^{p_{n+1}}\alpha_ie_i$, $n=1,2,\cdots$, is a bounded block sequence of $\{e_n\}_{n\ge1}$ such that $\inf_{1\le n<\infty}{\sup_{p_n+1\le i\le p_{n+1}}|\alpha_i|}>0$, then $\{y_n\}_{n\ge1}$ dominates $\{e_n\}_{n\ge 1}$.
\end{lemma}
%\begin{proof}
%    Since $\{e_n\}_{n\ge1}$ is symmetric, we may assume that there exists $\epsilon>0$ and $0\le k_n\le p_{n+1}-p_n$ such that $\alpha_{p_n+k_n}\ge\epsilon$ for $n=1,2,\cdots$.
%    Suppose that $\sum_{n=1}^{\infty}a_ny_n$ converges in $X$.
%    Then
%    $$
%    \norm{\sum_{i=1}^n a_ie_{p_i+k_i}}_{X}\le \frac{1}{\epsilon}\norm{\sum_{i=1}^n a_i \alpha_{p_i+k_i}e_{p_i+k_i}}_{X}\le\frac{1}{\epsilon} \norm{\sum_{i=1}^n a_iy_i}_{X}.
%    $$
%    Thus $\sum_{i=1}^{\infty}a_ie_{p_i+k_i}$ converges in $X$, so that $\sum_{i=1}^n a_ie_i$ converges in $X$.
%\end{proof}
Let $\{e_n\}_{n\ge1}$ be a symmetric  basis in a  $p$-Banach symmetric sequence space $X$.
If $0\not=\alpha=\sum_{i=1}^{\infty} \alpha_i e_i\in X$ and if $\{p_n\}_{n\ge1}$ is an increasing sequence of positive integers, let $$y_n^{(\alpha)}:=\sum_{i=p_n+1}^{p_{n+1}} \alpha_{i-p_n} e_i, ~n=1,2,\cdots.$$
Then $\{y_n^{(\alpha)}\}_{n\ge1}$ is a bounded block basic sequence of $\{e_n\}_{n\ge1}$.
And we call it a \textit{a block basic sequence of type I of $\{e_n\}_{n\ge1}$}.
Moreover, by Lemma~\ref{dominatelemma}, it is evident that  $\{y_n^{(\alpha)}\}_{n\ge1}$ dominates $\{e_n\}_{n\ge 1}$.
We denote by $\left[\{y_n\}_{n\ge 1}\right]$ the closed linear span of the sequence $\{y_n\}_{n\ge 1}$.

By \cite[Proposition 3.a.3]{LT1}, every symmetric basic sequence in Banach space is subsymmetric.
It also applies to the quasi-Banach spaces.
Consequently, we establish necessary and sufficient conditions for $\left[\{y_n\}_{n\ge 1}\right]$ to be isomorphic to $d(w,p)$.

\begin{lemma}\label{isomoequiva2}
    Let $0<p<\infty$ and $d(w,p)$ be the   weighted Lorentz sequence space,  where $w=\{w_n\}_{n\ge 1}$ is strictly decreasing to $0$.
     Let $\alpha=\sum_{i=1}^{\infty} \alpha_i e_i^{d(w,p)} \in d(w,p)$ with $\alpha_1\ge \alpha_2 \ge \cdots \ge 0$ and $\norm{\alpha}_{d(w,p)}=1$.
     If the block sequence $y_n=\sum_{i=p_n+1}^{p_{n+1}} \alpha_{i-p_n} e_i^{d(w,p)}$, $n=1,2,\cdots$ is symmetric, then $\left[\{y_n\}_{n\ge 1}\right]$ is isomorphic to $d(w,p)$ (equivalently, $\{y_n\}_{n\ge 1} \sim  \{e_n^{d(w,p)}\}_{n\ge 1}$) if and only if $s_n^{(\alpha)} \sim s_n$.
\end{lemma}
\begin{proof}
    Let $\{N_i\}_{i\ge 1}$ be 
    a collection of 
   subsets of natural numbers, $\mathbb{N}$, such that $\mathbb{N}=\bigcup_{i=1}^{\infty} N_i$, $N_i\cap N_j=\emptyset$ for all $i\not=j$, and the cardinality of  $N_i$ equals that of  $\mathbb{N}$, $n=1,2,\cdots$.
    For each $i$, we denote $N_i=\{k_{i,j}\}_{j\ge 1}$. 
    Let $u_i^{(\alpha)}=\sum_{j=1}^{\infty} \alpha_j e_{i,j}^{d(w,p)}$, where $e^{d(w,p)}_{i,j}:= e^{d(w,p)}_{   k_{i,j} } $.
    It is clear that $\{u_i^{(\alpha)}\}_{i\ge 1}$ is a symmetric basic sequence in $d(w,p)$ and $$\norm{u_i^{(\alpha)}}_{d(w,p)}=\norm{\alpha}_{d(w,p)}$$ for all $i\ge1$.

   Recall that   $d(w,p)$  is  a $q$-Banach space, where  $q=\min\{1,p\}$.
    For any $\epsilon>0$, we can find an increasing sequence $\{n_i\}_{i\ge 1}$ such that $$\norm{\sum_{j=p_{n_{i+1}}-p_{n_i}+1}^{\infty} \alpha_je_{i,j}^{d(w,p)}}_{d(w,p)}^q<\epsilon/2^i\enspace, \enspace i=1,2,\cdots.$$
    Let $z_i=\sum_{j=1}^{p_{n_{i+1}}-p_{n_i}} \alpha_je_{i,j}^{d(w,p)}$, $i=1,2,\cdots$.
    Then,
    $$
    \sum_{i=1}^{\infty} \norm{u_i^{(\alpha)}-z_i}_{d(w,p)}^q=\sum_{i=1}^{\infty}\norm{\sum_{j=p_{n_{i+1}}-p_{n_i}+1}^{\infty} \alpha_je_{i,j}^{d(w,p)}}_{d(w,p)}^q<\sum_{i=1}^{\infty} \epsilon/2^i=\epsilon.
    $$
    By the Principle of small perturbations\cite[Theorem 1.1]{KR}, $\{u_i^{(\alpha)}\}_{i\ge 1} \sim \{z_i\}_{i\ge 1}$.
    It is clear
    that $\{z_i\}_{i\ge 1} \sim \{y_{n_i}\}_{i\ge 1}$.
    Hence, $\{u_i^{(\alpha)}\}_{i\ge 1} \sim \{y_{n_i}\}_{i\ge 1}$ and $\{y_{n_i}\}_{i\ge 1}$ is symmetric.
     Since $\{y_{n}\}_{n \ge 1}$ is symmetric, it follows that  $\{y_{n_i}\}_{i\ge 1} \sim \{y_{n}\}_{i\ge 1}$.
    Hence, $\{u_i^{(\alpha)}\}_{i\ge 1} \sim \{y_{i}\}_{i\ge 1}$.
  Moreover,   $\norm{\sum_{i=1}^n u_i^{(\alpha)}}_{d(w,p)}^p=s_n^{(\alpha)}$.

    $(\Rightarrow)$: Suppose that $\left[\{y_n\}_{n\ge 1}\right]$ is isomorphic to $d(w,p)$.
    Then $\left[\{u_n^{(\alpha)}\}_{n\ge 1}\right]$ is isomorphic to $d(w,p)$, and by the uniqueness of symmetric basis of $d(w,p)$ (see \cite[Theorem 4]{ACL} and \cite[Proposition 3]{NO}), all symmetric bases in $d(w,p)$ are equivalent, which means $\{u_n^{(\alpha)}\}_{n\ge 1}$ is equivalent to $\{e_n^{d(w,p)}\}_{n\ge1}$.
    Thus, $$(s_n^{(\alpha)})^{1/p } = \norm{\sum_{i=1}^n u_i^{(\alpha)}}_{d(w,p)} \sim \norm{\sum_{i=1}^n e_i^{d(w,p)}}_{d(w,p)}=(s_n)^{1/p } ,$$which means $s_n^{(\alpha)} \sim s_n$.

    $(\Leftarrow)$: Suppose $s_n^{(\alpha)} \sim s_n$.
Recall that      $m_1^{(\alpha)}=s_1^{(\alpha)}$, $m_n^{(\alpha)}= s_{n+1}^{(\alpha)}-s_n^{(\alpha)}$, $n=2,3,\cdots$ and $m^{(\alpha)}=\left(m_1^{(\alpha)}, m_2^{(\alpha)}, \cdots\right)$ (see Definition \ref{definition m} above).
    By Proposition~\ref{mequiva}, $d\left(m^{(\alpha)},p\right)$ is isomorphic to $d(w,p)$.
    Let $\{v_n^{(\alpha)}\}_{n\ge1}$ be the unit vector basis of  $d\left(m
    ^{(\alpha)}, p\right)$.
    By the uniqueness of symmetric basis of $d(w,p)$ (see \cite[Theorem 4]{ACL} and \cite[Proposition 3]{NO}), $\{e_n^{d(w,p)}\}_{n\ge 1} \sim \{v_n^{(\alpha)}\}_{n\ge 1}$.
    Let $\{\beta_n\}_{n\ge 1}$ be an arbitrary  decreasing sequence of nonnegative numbers.
    Then
    $$
    \norm{\sum_{i=1}^N \beta_iu_i^{(\alpha)}}_{d(w,p)}^p=\sum_{i=1}^N \beta_i^p\left(\sum_{j=1}^{\infty} \alpha_j^p w_{i,j}\right)
    $$
    where for every $i=1,2,\cdots, N$ (respectively, for every $j$), $\{w_{i,j}\}_{j\ge 1}$ is a subsequence of $\{w_n\}_{n\ge 1}$ such that $\{w_{i,j}\}_{j\ge 1} \cap  \{w_{i' ,j}\}_{j\ge 1 }   =\emptyset $ when $i\ne i'$. 
    Now, for every $l$ and $k$, 
    $$
    \sum_{i=1}^k\left(\sum_{j=1}^l w_{i,j}\right)\le s_{kl}= \sum_{j=1}^l \left(s_{kj}-s_{k(j-1)}\right).
    $$
  By Proposition~\ref{comequiva},  for each fixed $k=1,2,\cdots, N$,  we have 
    $$
    \sum_{i=1}^k\left(\sum_{j=1}^{\infty} \alpha_j^p w_{i,j}\right) \le
    \sum_{j=1}^{\infty} \alpha_j^p \left(s_{kj}-s_{k(j-1)}\right)\stackrel{\tiny \rm Def.~\ref{definition m}}{=}s_k^{(\alpha)}=\sum_{i=1}^k m_i^{(\alpha)}.
    $$
    Since  $\{\beta_n\}_{n\ge 1}$ is decreasing, it follows from  Proposition~\ref{comequiva} that  $$\norm{\sum_{i=1}^N \beta_i u_i^{(\alpha)}}_{d(w,p)}^p\le \sum_{i=1}^N \beta_i^p m_i^{(\alpha)}.$$
    Hence,  $$\{v_n^{(\alpha)}\}_{n\ge 1} > \{u_n^{(\alpha)}\}_{n\ge 1}.$$ 
    Since $\{e_n^{d(w,p)}\}_{n\ge 1} \sim \{v_n^{(\alpha)}\}_{n\ge 1} $ and $\{y_n\}_{n\ge 1} \sim \{u_n^{(\alpha)}\}_{n\ge 1}$, we get $\{e_n^{d(w,p)}\}_{n\ge1}>\{y_n\}_{n\ge 1}.$
    By Lemma~\ref{dominatelemma}, $$\{y_n\}_{n\ge 1}>\{e_n^{d(w,p)}\}_{n\ge 1}$$
Thus,  $\{y_n\}_{n\ge 1}\sim\{e_n^{d(w,p)}\}_{n\ge 1}$ which means $\left[\{y_n\}_{n\ge 1}\right]$ is isomorphic to $d(w,p)$.
This completes the proof.
\end{proof}
   Below, we characterize the weighted Lorentz sequence spaces $d(w,p)$ for $0<p<\infty$ which  have only two nonequivalent symmetric basic sequences.
\begin{theorem}\label{fourmainthm}
    Let $0<p<\infty$ and $d(w,p)$ be a    weighted Lorentz sequence space,  where $w=\{w_n\}_{n\ge 1}$ is strictly decreasing to $0$.
    Then,  there are exactly two nonequivalent symmetric basic sequences in $d(w,p)$ (namely, the unit vector basis of $\ell_p$ and that of $d(w,p)$ itself) if and only if $\sup_{1\le n
    ,k<+\infty} s_{nk}/s_ns_k<+\infty$.
\end{theorem}
\begin{proof}
    $(\Leftarrow)$: Let $\{y_n\}_{n\ge 1}$ be a symmetric basic sequence in $d(w,p)$.
    By \cite[Lemma 2.1]{KR07}, any basic sequence in $d(w,p)$ is equivalent to a block basic sequence. 
    By
    \cite[Theorem 3]{ACL} and \cite[Theorem 3.2]{Popa}, $\{y_n\}_{n\ge 1}$ is equivalent either to the unit vector basis of $\ell_p$ or to a block basic sequence of type I of $\{e_n^{d(w,p)}\}_{n\ge 1}$, 
    where $\{e_n^{d(w,p)}\}_{n\ge 1}$ is the unit vector basis of $d(w,p)$.
   If $\{y_n\}_{n\ge 1}$ is not equivalent to the unit vector basis of $\ell_p$, then there exists  $\alpha=\sum_{i=1}^{\infty} \alpha_i e_i^{d(w,p)} \in d(w,p)$ with $\alpha_1\ge \alpha_2 \ge \cdots \ge 0$ and $\norm{\alpha}_{d(w,p)}=1$
     such that $y_n=\sum_{i=p_n+1}^{p_{n+1}} \alpha_{i-p_n} e_i^{d(w,p)}$, $n=1,2,\cdots$.
    Since $\sup_{1\le n
    ,k<+\infty} s_{nk}/s_ns_k<+\infty$, 
    it follows from  Theorem~\ref{supequiva} that  $s_n^{(\alpha)} \sim s_n$. 
   % Since $\{y_n\}_{n\ge 1}$ is symmetric, by Lemma~\ref{isomoequiva2}, $\left[\{y_n\}_{n\ge 1}\right]$ is isomorphic to $d(w,p)$.
   %  Moreover, by \cite[Theorem 4]{ACL} and \cite[Proposition 3]{NO},  $\left\{e_n^{d(w,p)}\right\}_{n\ge 1}$ is the  unique symmetric basis in $d(w,p)$  (up to equivalence).
   Since $\{y_n\}_{n\ge1}$ is symmetric, by Lemma~\ref{isomoequiva2}, $\{y_n\}_{n\ge 1}$ is equivalent to the unit vector basis of $d(w,p)$.

    $(\Rightarrow)$: Now, assume by contradiction that  $\sup_{1\le n
    ,k<+\infty} s_{nk}/s_ns_k=+\infty$.
    By Theorem~\ref{supequiva}, there exists  $\alpha \in d(w,p)$ with $\alpha_1\ge \alpha_2 \ge \cdots \ge 0$ and $\norm{\alpha}_{d(w,p)}=1$ such that $s_n^{(\alpha)} \not\sim s_n$.
    Similar to the proof of Lemma~\ref{isomoequiva2}, we can construct a block symmetric basic sequence $\left\{z_n^{(\alpha)}\right\}_{n\ge 1}$ of type I of $\left\{e_n^{d(w,p)}\right\}_{n
    \ge 1}$,  where $z_n^{(\alpha)}=\sum_{i=p_n+1}^{p_{n+1}}\alpha_{i-p_n}e_i^{d(w,p)}$ for $i=1,2,\cdots$.
    Since $\left\{z_n^{(\alpha)}\right\}_{n\ge 1}$ is symmetric and $s_n^{(\alpha)} \not\sim s_n$, by Lemma~\ref{isomoequiva2}, $\left\{z_n^{(\alpha)}\right\}_{n\ge 1}$ is not equivalent to
    $\left\{e_n^{d(w,p)}\right\}_{n
    \ge 1}$.
    Now, it is sufficient to prove that $\left\{z_n^{(\alpha)}\right\}_{n\ge 1}$ is not equivalent to the unit vector basis of $\ell_p$.
    Since $\lim_{i\rightarrow \infty} \norm{\sum_{n=1}^{\infty} \alpha_{n+i}e_n^{d(w,p)}}_{d(w,p)}=0$, for any $\epsilon>0$, there exists a positive integer $N_{\epsilon}$ such that $$\norm{\sum_{n=1}^{\infty} \alpha_{N_{\epsilon}+n} e_n^{d(w,p)}}_{d(w,p)}^p<\left(\alpha_1 \epsilon\right)^p/2.$$
    Since $\lim_{n\rightarrow \infty} w_n=0$, we have that $\lim_{n\rightarrow \infty} \left(\sum_{i=1}^{nN_{\epsilon}} w_i\right)/n=0$. 
    Choose $n$ such that $\left(\sum_{i=1}^{nN_{\epsilon}} w_i\right)/n<\epsilon^p/2$.
    Then,
    \begin{align*}
        \norm{\sum_{i=1}^{n} z_{N_\epsilon+i}^{(\alpha)}}_{d(w,p)}^p
        &=
        s_n^{(\alpha)}\\
        &=\sum_{i=1}^{\infty} \alpha_i^p \left(\sum_{j=1}^n w_{n(i-1)+j}\right)\\
        &\le
        \sum_{i=1}^{N_\epsilon} \alpha_i^p\left(\sum_{j=(i-1)n+1}^{in} w_j\right)+\sum_{i=N_\epsilon+1}^{\infty} \alpha_i^p (nw_i)\\
        &=
        \sum_{i=1}^{N_\epsilon} \alpha_i^p\left(\sum_{j=(i-1)n+1}^{in} w_j\right)+n\norm{\sum_{i=1}^{\infty}\alpha_{N_\epsilon+i}e_i^{d(w,p)}}_{d(w,p)}^p\\
        &<
        \alpha_1^p\sum_{i=1}^{nN_\epsilon}w_i+n\left(\frac{\alpha_1^p\epsilon^p}{2}\right)\\
  &      <n\left(\frac{\alpha_1^p\epsilon^p}{2}\right)+n\left(\frac{\alpha_1^p\epsilon^p}{2}\right)\\
&        =n\alpha_1^p\epsilon^p
    \end{align*}
    Hence, $\norm{\sum_{i=1}^n z_{N_\epsilon+i}^{(\alpha)}}_{d(w,p)}<n^{1/p}\alpha_1\epsilon$.
    However,  $n^{1/p}=\norm{\sum_{i=1}^n e_i^{\ell_p}}_{\ell_p}$,  where $\{e_n^{\ell_p}\}_{n\ge 1}$ stands for  the unit vector basis of $\ell_p$.
    Thus $\left\{z_n^{(\alpha)}\right\}_{n\ge 1}$ is not equivalent to $\left\{e_n^{\ell_p}\right\}_{n\ge 1}$\footnote{Indeed, arguing similarly, one can show that no subsequence of $\left\{z_n^{(\alpha)}\right\}_{n\ge 1}$ is equivalent to $\left\{e_n^{\ell_p}\right\}_{n\ge 1}$.}.
    Thus, there are more than two nonequivalent symmetric basic sequences in $d(w,p)$.
\end{proof}

The following corollary is a direct result of Theorem~\ref{fourmainthm}.
\begin{corollary}\label{ddequiva}
    Let $0<p_1, p_2<\infty$.
    Let $d(w_1,p_1)$ and $d(w_2,p_2)$ be the weighted Lorentz sequence space,  where $w_1=\{w_n^1\}_{n\geq 1}$ is a nonincreasing  sequence 
    and $w_2=\{w_n^2\}_{n\ge 1}$ is a  sequences strictly decreasing to $0$.
    Let $s_n:=\sum_{i=1}^n w_i^2$.
    If $\sup_{1\le n
    ,k<+\infty} s_{nk}/s_ns_k<+\infty$, then $d(w_1,p_1)\hookrightarrow d(w_2,p_2)$ implies that  $$\mbox{either~}\{e_n^{d(w_1,p_1)}\}_{n\ge 1}\sim\{e_n^{\ell_{p_2}}\}_{n\ge1}\mbox{ or }\{e_n^{d(w_1,p_1)}\}_{n\ge 1}\sim\{e_n^{d(w_2,p_2)}\}_{n\ge1}.$$
\end{corollary}

\begin{rem}\label{equivapp}
Let   $w_1=\{w_n^1\}_{n\geq 1}$ and $w_2=\{w_n^2\}_{n\ge 1}$  be two  nonincreasing  sequences.
    If $d(w_1,p_1)\hookrightarrow d(w_2,p_2)$, then $p_1=p_2$.
    Indeed, by \cite[Lemma 2.1]{KR07}, the standard basis of $d(w_1,p_1)$ is equivalent to a block basis $(y_n)_{n \ge 1 }$ in $d(w_2,p_2)$.  
    Let $$y_n=\sum_{i= k_{n}+1 }^{k_{n+1 } } \alpha_{  i} e_{ i}^{d(w_2,p_2)} .$$ 
    Since $d(w_2,p_2)$ is symmetric, we may assume that $\alpha_{(n)}: =(\alpha_{ i})_{i= k_{n}+1 }^{k_{n+1 } }$ is decreasing. 
    Let $T: d(w_1,p_1) \rightarrow \overline{\text{span}\{y_n\}}$ be an isomorphism,   where $\overline{\text{span}\{y_n\}}$ is the closure of   $\text{span}\{y_n\}$  in  $d(w_2,p_2)$.
    By Lemma~\ref{Lorequivalp}, 
		we obtain that $$\left( \frac{\sum_{i=2^n +1}^{2^{n+1}} e_i ^{d(w_1,p_1)}}{\norm{\sum_{i=2^n +1}^{2^{n+1}} e_i ^{d(w_1,p_1)}}_{d(w_1,p_1)}}  \right)_{n\ge 1}$$ contains a subsequence which is equivalent to the vector basis of $\ell_{p_1}$.
		Since $T$ is an isomorphism, passing to a further subsequence, 
		$$\left( \frac{\sum_{i=2^n +1}^{2^{n+1}} T(e_i ^{d(w_1,p_1)})}{\norm{\sum_{i=2^n +1}^{2^{n+1}} T(e_i ^{d(w_1,p_1)})}_{d(w_2,p_2)}}  \right)_{n\ge 1}$$ contains a subsequence which is equivalent to the vector basis of $\ell_{p_2}$. 
        In other words, $\ell_{p_1}\hookrightarrow \ell_{p_2}$. 
        Therefore, we  have  $p_1=p_2$ (see \cite[Proposition 2.9]{KPR} or Lemma \ref{classical} above).
    %   {\color{red}
    %    \cite[Corollary 14]{CL} and \cite[Theorem 4.8]{Popa} imply that if $d(w_1,p)$ is isomorphic to $d(w_2,p)$ for $0<p<\infty$, then the embedding: $d(w_1,p)\hookrightarrow d(w_2,p)$ must be complemented.}
\end{rem}
       
%{\color{red}This motivates the following question:
%\begin{quest}\label{myquestion}
%		 Does there exist an isomorphic embedding $T:d(w_1,p)\hookrightarrow d(w_2,p)$ such that $T$ is not complemented when $d(w_1,p)\not=\ell_p$ and $\sup_{1\le n
 %   ,k<+\infty} s_{nk}/s_ns_k=+\infty$?
%		\end{quest}
%If we have a negative answer to Problem~\ref{myquestion}, then Remark~\ref{equivapp} shows us that if $d(w_1,p_1)\hookrightarrow d(w_2,p_2)$, then either $d(w_1,p_1)\approx\ell_{p_2}$ or $d(w_1,p_1)\approx d(w_2,p_2)$.
%}

\end{document}